\documentclass[a4paper,5p,number,sort&compress,times]{elsarticle}
\usepackage{amsmath,amsfonts,amssymb}
\usepackage{graphicx,psfrag,subfigure}

\newtheorem{thm}{Theorem}

\newdefinition{defi}{Definition}
\newdefinition{remark}{Remark}
\newdefinition{example}{Example}

\newcommand{\Order}{\mathop{\rm O}\nolimits}

\newcommand{\rmd}{{\rm d}}
\newcommand{\rme}{{\rm e}}
\newcommand{\rmh}{{\rm h}}
\newcommand{\rmi}{\mathrm{i}}
\newcommand{\rmp}{{\rm p}}
\newcommand{\rmr}{{\rm r}}
\newcommand{\rms}{{\rm s}}
\newcommand{\rmu}{{\rm u}}
\newcommand{\RF}{{\rm RF}}
\newcommand{\acos}{\mathop{\rm acos}\nolimits}

\newcommand{\Trace}{\mathop{\rm tr}\nolimits}
\newcommand{\FixedSet}{\mathop{\rm Fix}\nolimits}
\newcommand{\Identity}{{\rm I}}

\newcommand{\Nset}{\mathbb{N}}
\newcommand{\Qset}{\mathbb{Q}}
\newcommand{\Rset}{\mathbb{R}}
\newcommand{\Tset}{\mathbb{T}}
\newcommand{\Zset}{\mathbb{Z}}

\newif\iffigures
\figurestrue

\begin{document}

\begin{frontmatter}

\title{Stability of the phase motion in race-track microtons\tnoteref{t1}}

\author[iet]{Yu.~A.~Kubyshin}
\ead{Iouri.Koubychine@upc.edu}

\author[mazu]{O.~Larreal}
\ead{olarreal@luz.edu.ve}

\author[ma]{R.~Ram{\'\i}rez-Ros\corref{cor}}
\ead{Rafael.Ramirez@upc.edu}

\author[ma]{T. M. Seara}
\ead{Tere.M.Seara@upc.edu}

\cortext[cor]{Corresponding author}

\address[iet]{Institut de T\`{e}cniques Energ\`{e}tiques,
             Universitat Polit\`ecnica de Catalunya,
             Diagonal 647, 08028 Barcelona, Spain}

\address[mazu]{Departamento de Matem\'aticas,
               Universidad del Zulia,
               Maracaibo, Venezuela}

\address[ma]{Departament de Matem\`{a}tiques,
             Universitat Polit\`ecnica de Catalunya,
             Diagonal 647, 08028 Barcelona, Spain} 
             
\tnotetext[t1]{Research supported in part by
MINECO-FEDER grant MTM2015-65715-P(Spain)
and CUR-DIUE grants 2014SGR504 and 2014SGR846 (Catalonia).
T. M. Seara is also supported by the Russian Scientific Foundation
grant 14-41-00044 and the European Marie Curie Action
FP7-PEOPLE-2012-IRSES: BREUDS. 
We acknowledge the use of the UPC Applied Math cluster system for
research computing.
Useful conversations with Amadeu Delshams, Vasiliy Shvedunov,
and Arturo Vieiro are gratefully acknowledged.}

\begin{abstract}
We model the phase oscillations of electrons in race-track microtrons 
by means of an area preserving map with a fixed point at the origin, which
represents the \emph{synchronous trajectory} of a reference particle in the
beam. We study the nonlinear stability of the origin in terms of the
\emph{synchronous phase} ---the phase of the synchronous particle
at the injection.
We estimate the size and shape of the stability domain around the origin,
whose main connected component is enclosed by the last rotational invariant
curve.
We describe the evolution of the stability domain as the synchronous phase varies.
Besides, we approximate some rotational invariant curves by level sets of certain
Hamiltonians. Finally, we clarify the role of the stable and unstable invariant
curves of some hyperbolic (fixed or periodic) points.
\end{abstract}


\end{frontmatter}

\section{Introduction}
\label{Sec:Intro}

Race-track microtron (RTM) is a specific type of electron accelerator
with beam recirculation combining properties of the linear accelerator
and a circular machine~\cite{KapitzaMelekhin1978,Rand1984}.
For applications in which a modest beam power at a relatively
high beam energy is required the RTM turns out to be
the most optimal source of electron beams. 
This is the case of applications like the cargo inspection or
Intraoperative Radiation Therapy for which RTMs allow
to get pulsed and continuous beams in a quite cost and energy effective
way with the most optimal dimensions of the machine.
An example of such accelerator is a compact 12 MeV RTM
which is under construction at the Technical University of Catalonia
in collaboration with the Moscow State University and CIEMAT
(Spain)~\cite{Aloev_etal2010,Kubyshin_etal2009,Vladimirov_etal2014}.

In the design of a particle accelerator its main parameters are optimized
for some reference particle, usually referred to as synchronous particle.
Real particles of the beam perform transverse and longitudinal oscillations
with respect to this synchronous particle.
One of the important issues of the RTM design is to assure the stability
of these oscillations,
in particular to avoid an uncontrolled growth of their amplitude that
leads to the loss of the beam.
In the RTM beam dynamics the main role is played by longitudinal oscillations
of individual electrons with respect to the synchronous
particle~\cite{Veksler1945,Lidbjork2001}.
Such oscillations are usually referred to as phase motion,
and the region of initial states in the phase space which give rise to
stable oscillations is called acceptance.

Stability of the phase oscillations of particles in the beam is a
matter of primary concern both at the stage of the accelerator design
and during its operation. 
In particular,
the size of the acceptance determines the efficiency of caption into
acceleration of particles injected from a source of electrons,
usually an electron gun. 
The shape and size of the acceptance depend on the phase of
the synchronous particle, also called synchronous phase and
denoted as $\phi_\rms$ along this paper. 
The definition of these notions will be given in Section~\ref{Sec:RTMModel}.
A specific feature of the RTM is that the range of values of $\phi_\rms$ for
which the acceptance exists is quite narrow and therefore this parameter
must be carefully chosen and controlled during the RTM operation.
The main criteria here is to keep the synchronous phase equal or close
to the value for which the size of the acceptance is maximal and to
avoid resonant values of $\phi_\rms$. 
The latter is important
because resonant phase oscillations lead to a buildup of beam instabilities,
rapid growth of the amplitude of the phase oscillations and eventually
to a loss of the beam. 
Therefore, a good understanding of the longitudinal
dynamics of the beam in RTMs in general and the phase oscillations
of electrons in a particular machine is important.

The phase oscillations of the beam in an RTM are described by a system
of nonlinear difference equations~\cite{Henderson_etal1953,Rand1984} that
will be derived in Section~\ref{Sec:RTMModel}.
It gives rise to an analytic area preserving map which is the object of study
in the present paper. 
These difference equations cannot be approximated by differential
equations without loss of accuracy, since RTMs have large energy gain
per turn and high frequency of phase oscillations. 
The description of the acceptance in the linear approximation
is well known~\cite{Rand1984}.
The change of the acceptance under the variation of the synchronous phase,
the appearance of stable regions,
and the emergence of the stochastic regime is studied
in~\cite{Melekhin1972} by looking at the normal form of
the difference equations up to order four. 
For instance, it was shown that there is a shift of the center of
oscillations and that the oscillation frequency depends on the amplitude. 
We will reproduce these known results at the beginning of our analysis.

In the present paper we will study the phase motion of the beam in
RTMs using modern Dynamical Systems methods. 
Our goal is to characterize the acceptance,
called stability domain in Dynamical Systems,
as a function of the synchronous phase. 
In particular,
we will establish the intervals of values of $\phi_\rms$ for which
the acceptance and its central connected component exist,
analyze their geometry, shape and size, and calculate their area for
the whole range of $\phi_\rms$.
Our main theoretical tools are the Moser twist theorem and several
Sim\'o's stability results~\cite{Simo1982}. 
We will also approximate some rotational invariant curves of
the RTM beam longitudinal dynamics by level curves of suitable Hamiltonians. 
Finally, we will carry out a global study of the acceptance and
its connected component using results and algorithms, like the orbit method,
developed in~\cite{LuqueVillanueva2013,Vieiro2009}.

The article is organized as follows. 
We give a short introduction into the RTM beam longitudinal dynamics,
define the notion of synchronous trajectory and synchronous phase
and derive the map describing the longitudinal beam dynamics that
we will refer to as RTM map in Section~\ref{Sec:RTMModel}.
We introduce the notion of the acceptance and its connected component,
we find fixed points of the RTM map and analyze their type in the linear
approximation in Section~\ref{Sec:MainResult}. 
We also formulate a theorem about the local stability of
the synchronous trajectory and give a summary of further results on
the behavior of the acceptance as a function of the synchronous phase. 
The proof of the theorem is given in Section~\ref{Sec:LocalStability}.
Section~\ref{Sec:HamiltonianApproximations} is devoted to the approximation
of the RTM map by the integrable dynamics of certain Hamiltonians
and correspondingly the approximation of the invariant curves around
the synchronous trajectory by means of Hamiltonian level curves. 
The global stability of the synchronous trajectory is studied
in Section~\ref{Sec:GlobalStability},
also a detailed characterization of the acceptance is given there.
These results are obtained numerically by using the orbit method.
The invariant curves of hyperbolic points are analyzed
in Section~\ref{Sec:InvariantCurves}. 
A summary and a discussion of the obtained results
are given in Section~\ref{Sec:Conclusions}.

\section{The RTM model}
\label{Sec:RTMModel}

\begin{figure}
\includegraphics*[width=85mm]{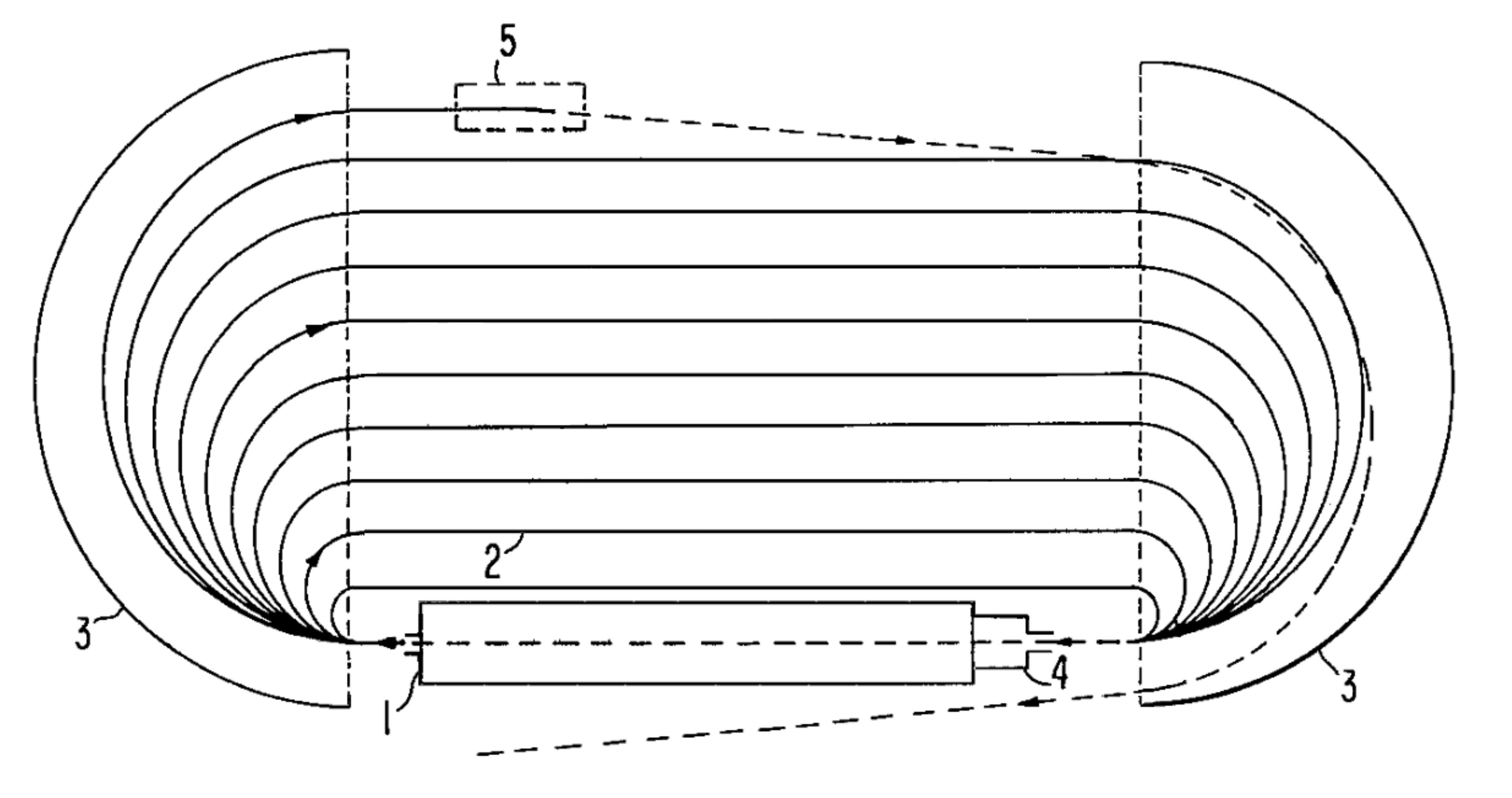}
\caption{Schematic view of our RTM model:
1) Accelerating structure,
2) Drift space,
3) End magnets,
4) Electron gun, and
5) Extraction magnet.}
\label{fig:MicrotronScheme}
\end{figure}

The operation of a race-track microtron (RTM) is illustrated
in Fig.~\ref{fig:MicrotronScheme}.
The initial beam is injected from an electron source
(electron gun or external pre-accelerator)
into an accelerating structure (AS) consisting of a few resonant cavities.
The longitudinal electric component of a high frequency electromagnetic
wave (usually a standing wave) accelerates the electrons in the AS.
Then the beam is bent by the magnetic field of a $180^{\circ}$ bending magnet,
called \emph{end magnet}, travels through a free space,
usually referred to as \emph{drift space}, follows along a circular trajectory
inside the second end magnet and returns to the AS. 
In this way the beam makes a number of recirculations through the RTM with
the energy being increased at each pass through the AS.
Once the beam gets the final design energy,
it is deflected by an extraction magnet and is
directed towards the accelerator beam exit.

We approximate the AS by a zero length accelerating gap and consider
a stationary regime with a constant amplitude of the accelerating field. 
Then the energy gain of an electron passing through the AS is equal to
$\Delta_{\max} \cos\phi$, where $\Delta_{\max}$ is the maximum energy gain
in the AS and $\phi$ is the phase of the accelerating field. 
In accelerator physics this parameter is usually referred to
as particle phase and is used to characterize the longitudinal position
of the particle along the orbit with respect to the accelerating gap
(in a real machine with respect to, say, the exit, of the last cavity
of the AS in the direction of the beam motion).
Another idealization in our study is that the end magnets will be considered
as hard-edge dipole magnets so that the fringe-field effects
are not taken into account.

Let us consider the longitudinal (phase) motion of electrons in an
RTM with magnetic field induction $B$ in the end magnets
and separation $l$ between the magnets (straight section length).
We assume that the injected electrons are already ultra-relativistic,
so that in the formulas below the velocity of the particles in the
beam is equal to the velocity of light $c$.

Our dynamical variables are the full particle energy $E_n$ and
its phase $\phi_n$ at the $n$-th turn at some point of the orbit.
For the sake of convenience, we choose this point to be the exit of the AS. 
Let $\phi_0$ and $E_0$ be the particle phase and energy at the injection;
that is, just before its first passage through the AS. 
We would like to note that, in some pulsed RTM,
the electrons reverse their trajectory in special reverse field magnets
after the injection and first acceleration at the
AS~\cite{Kubyshin_etal2009,Vladimirov_etal2014}.
In that case,
$E_0$ is the energy after the first acceleration and reflection of the beam
in the end magnet, before the second passage through the AS,
but we will still use the term injection energy.

For the RTM model described above,
the duration of the $n$-th revolution of the beam is 
\begin{equation}\label{eq:TimeEnergy}
T_n =
\frac{2 l + 2 \pi r_n}{c} =
\frac{2l}{c} + \frac{2 \pi E_n}{e c^2 B},
\end{equation}
where $e$ is the elementary charge and $r_n$ is the beam trajectory radius
in the end magnets. 
We have used that $E_n = e c B r_n$ in bending
magnets for ultrarelativistic beams~\cite{Rand1984}. 
Relation~(\ref{eq:TimeEnergy}) allows us to work indistinctly with
time variables $T_n$ or energy variables $E_n$. 
The beam dynamics in the phase-energy variables, for a zero-length AS,
is governed by the difference equations
\begin{equation}\label{eq:Dynamics_in_phi_E}
\phi_{n+1} = \phi_n + 2 \pi T_n/T_{\RF}, \qquad
E_{n+1} = E_n + \Delta_{\max} \cos\phi_{n+1},
\end{equation}
$T_{\RF}$ being the period of the accelerating electromagnetic field,
usually called \emph{radiofrequency (RF) field} in the context
of accelerator beam dynamics.
See~\cite{Rand1984}.

The design of a particle accelerator requires a proper choice of physical
and technical parameters in order to guarantee the existence of a
reference trajectory, called \emph{synchronous trajectory},
so that an ideal particle following this trajectory is accelerated
in the most optimal way.
Namely, it can be tuned in resonance with the accelerating field.
Let us explain the choice of parameters in case of an RTM. 
We fix two positive integers $m$ and $k$.
The resonance conditions in the case of an RTM are
\begin{equation}\label{eq:ResonanceConditionsInTime}
\frac{T_{1,\rms}}{T_\RF} = m, \qquad
\frac{T_{n+1,\rms}-T_{n,\rms}}{T_\RF} = k.
\end{equation}
Henceforth, the subindex ``$\rms$'' indicates quantities related
to the synchronous trajectory. 
The integer $m$ is the number of RF field periods during the first turn
and defines the synchronicity condition at this turn,
whereas $k$ is the increase of the multiplicity factor due to the increase
of the  period of revolution of the reference particle in each turn.
The ratio 
\[
j_n := T_{n,\rms}/T_\RF = m + (n-1)k
\]
is called \emph{harmonic number}. 
RTMs are accelerators with variable harmonic numbers,
since $j_n$ depends on the $n$-th turn.

If we rewrite the resonance conditions~(\ref{eq:ResonanceConditionsInTime})
using the energy variables $E_{n,\rms}$ and relation~(\ref{eq:TimeEnergy}),
we get that $E_{n,\rms} = E_\rms + n \Delta_\rms$, where 
\[
E_\rms= \left( \frac{m}{k} - 1 - \frac{2l}{k\lambda} \right) \Delta_\rms, \qquad
\Delta_\rms = \frac{e c B \lambda k}{2 \pi},
\]
and $\lambda = c T_\RF$ is the wavelength of the RF field.
Finally, it is straightforward to check that if the injection
phase for the synchronous trajectory
$\phi_{0,\rms} = \phi_{1,\rms} = \phi_\rms$ satisfies relation 
\[
\Delta_\rms = \Delta_{\max} \cos\phi_\rms,
\]
then we get a particular solution $(\phi_{n,\rms},E_{n,\rms})$ of
equations~(\ref{eq:TimeEnergy})--(\ref{eq:Dynamics_in_phi_E}), whose
energy undergoes a constant gain $\Delta_\rms$ at each turn: 
\begin{equation}\label{eq:SynchronousDynamics}
\phi_{n,\rms} = \phi_\rms + 2\pi i_n,\qquad
E_{n,\rms} = E_\rms + n \Delta_\rms,
\end{equation}
with $i_n = (n-1)m + (n-1)(n-2)k/2$.
We note that $i_n \in \Zset$,
so this synchronous particle passes through the AS in the same phase
of the RF field at each turn:
$\phi_{n,\rms} = \phi_\rms$ (mod $2\pi$) for all $n$.
We say that $\phi_\rms$ is the \emph{synchronous phase}.

Once we realize that
the synchronous trajectory~(\ref{eq:SynchronousDynamics}) exists,
we wonder whether the oscillations of other trajectories around it are stable.
This depends on the synchronous phase $\phi_\rms$.

We study this dependence by introducing the variables 
\begin{equation}\label{eq:psi-w-defn}
\psi_n = \phi_n-\phi_\rms, \qquad
w_n = 2 \pi k (E_n-E_{n,s})/\Delta_\rms,
\end{equation}
that describe the phase and energy deviation of an arbitrary trajectory
from the synchronous one.
Then the beam dynamics is modeled by the map
$(\psi_{n+1},w_{n+1})=f(\psi_n,w_n)$, where
\[
\left\{
\begin{array}{ccl}
\psi_{n+1} & = & \psi_n+w_n, \\
w_{n+1} & = & w_n + 2 \pi k \big( \cos(\psi_{n+1} + \phi_\rms )/\cos\phi_\rms-1 \big).
\end{array}
\right.
\]
For simplicity, we will assume that $k=1$ in the rest of the paper.
The study of other values can be carried out in a similar way.

We end this section by stressing that the model of longitudinal oscillations
around the synchronous trajectory for non-ultra-relativistic beams
is more complicated~\cite{Kubyshin_etal2008}.

\section{Terminology and main results}
\label{Sec:MainResult}

Under the assumption $k = 1$, our model for the nonlinear oscillations around
the synchronous trajectory~(\ref{eq:SynchronousDynamics}) is
the analytic area preserving diffeomorphism $(\psi_1,w_1) = f(\psi,w)$, defined by
\begin{equation}\label{eq:f}
\left\{
\begin{array}{ccl}
\psi_1 & = & \psi + w,\\
w_1    & = & w + 2\pi (\cos \psi_1 -1) - \mu \sin \psi_1.
\end{array}
\right.
\end{equation}
Here, $\mu = 2\pi \tan \phi_\rms$ is a parameter, and
$\psi$ (respectively, $w$) is the deviation of
the phase (respectively, energy) of an arbitrary particle from the
phase (respectively, energy) of the synchronous trajectory,
which corresponds to the fixed point
\[
p_\rms = (\psi_\rms,w_\rms) = (0,0).
\]
The angular coordinate $\psi$ is defined modulo $2\pi$,
so our phase space is $\Tset \times \Rset$, with $\Tset = \Rset/2\pi\Zset$.
For brevity, we will write $p = (\psi,w)$ and
$p_n = (\psi_n,w_n) = f^n(\psi,w) = f^n(p)$.
We look for initial conditions that give rise to particles with
a bounded energy deviation from the synchronous trajectory.
More precisely, we will estimate the size and the shape of
\[
\mathcal{A} =
\left\{ p \in \Tset \times \Rset : \mbox{$(w_n)_{n \in \Zset}$ is bounded} \right\}.
\]
This domain is called \emph{(longitudinal) acceptance} in
Accelerator Physics, and \emph{stability domain} in Dynamical Systems.
We will also study its connected component $\mathcal{D} \subset \mathcal{A}$
that contains $p_\rms$.

\begin{remark}\label{rem:PhaseDeviation}
We have numerically checked that
\[
\mathcal{A} \subset [-0.45,0.35] \times [-0.8,0.8],\qquad
\forall \mu > 0.
\]
Hence, if $(\psi_n,w_n) \in \Tset \times \Rset$ is the phase-energy
deviation at the $n$-th turn of a trajectory contained in $\mathcal{A}$
and $(\tilde{\psi}_n,w_n) \in \Rset \times \Rset$ denotes the lift of
this trajectory determined by $-\pi \le \tilde{\psi}_0 < \pi$,
then $-\pi \le \tilde{\psi}_n < \pi$ for all $n$.
This means that trajectories with bounded energy deviation,
have also bounded phase deviation respect to the synchronous one.
\end{remark}

The map $f$, its acceptance $\mathcal{A}$, and the connected component
$\mathcal{D}$ depend on $\mu$, but we will frequently omit this dependence.
Otherwise, we will write $f_\mu$, $\mathcal{A}_\mu$,
and $\mathcal{D}_\mu$, respectively.

We say that the synchronous trajectory is:
\begin{itemize}
\item
\emph{Globally stable} when $\mathcal{D}$ is a neighborhood of $p_\rms$;
\item
\emph{Locally stable} when $\forall \epsilon > 0$ there exists
$\delta > 0$ such that
\[
\| p-p_\rms \| < \delta \Rightarrow
\|p_n - p_\rms \| < \epsilon ,\ \forall n \in \Zset;
\]
\item
\emph{Linearly stable} when the sequence $(M_\rms^n)_{n \in \Zset}$
is bounded, where $M_\rms$ is the linear part of the map $f$
at $p_\rms$; that is, when the linearized map $p  \mapsto M_\rms p$
is locally stable.
\end{itemize}
Otherwise, it is \emph{globally}, \emph{locally}, or \emph{linearly unstable}.

\begin{figure*}
\centering
\subfigure[ICs and two unbounded orbits for $\mu = -0.1$]{
\label{fig:BeforeParabolic0}
\iffigures
\includegraphics[height=1.62in]{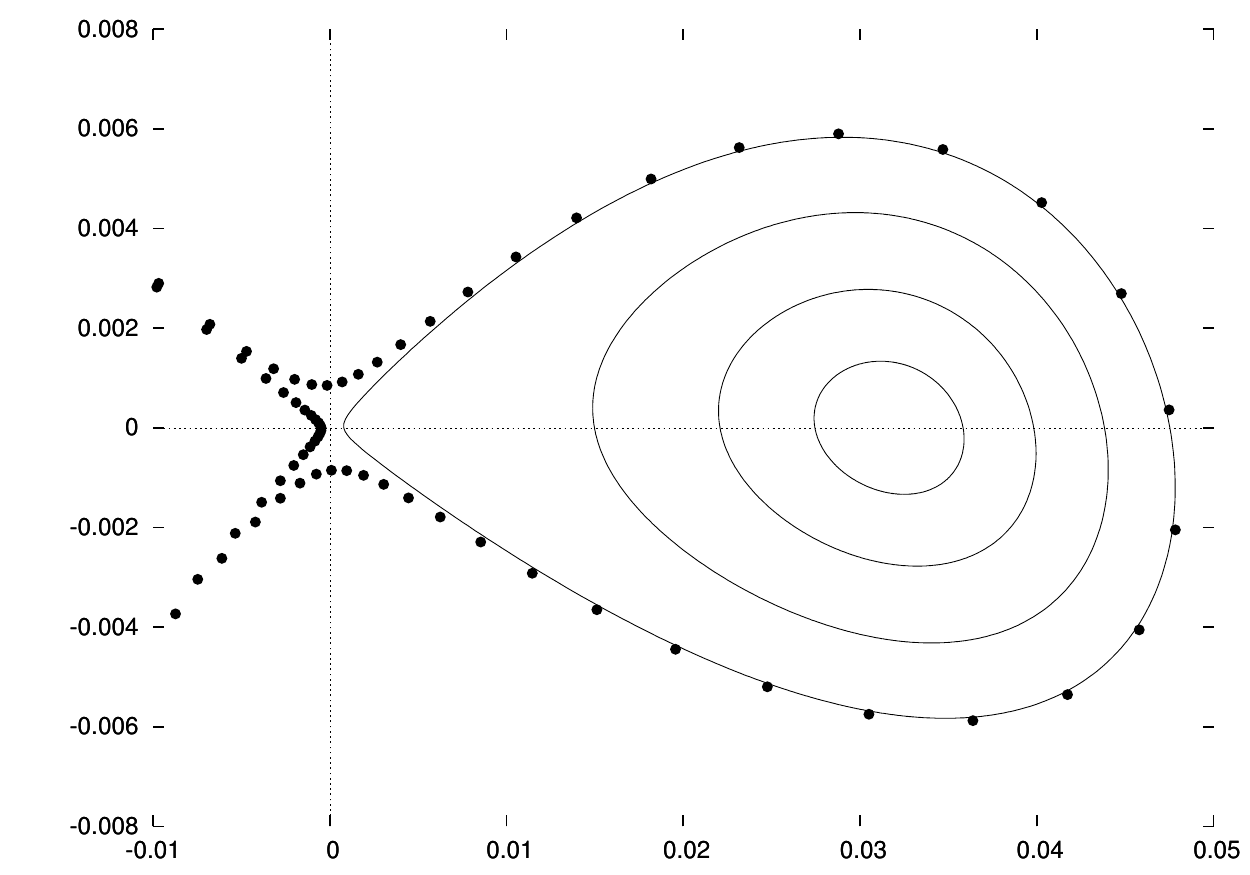}
\else
\vspace{1.62in}
\fi
}
\hfill
\subfigure[Five unbounded orbits for $\mu = 0$]{
\label{fig:Parabolic0}
\iffigures
\includegraphics[height=1.62in]{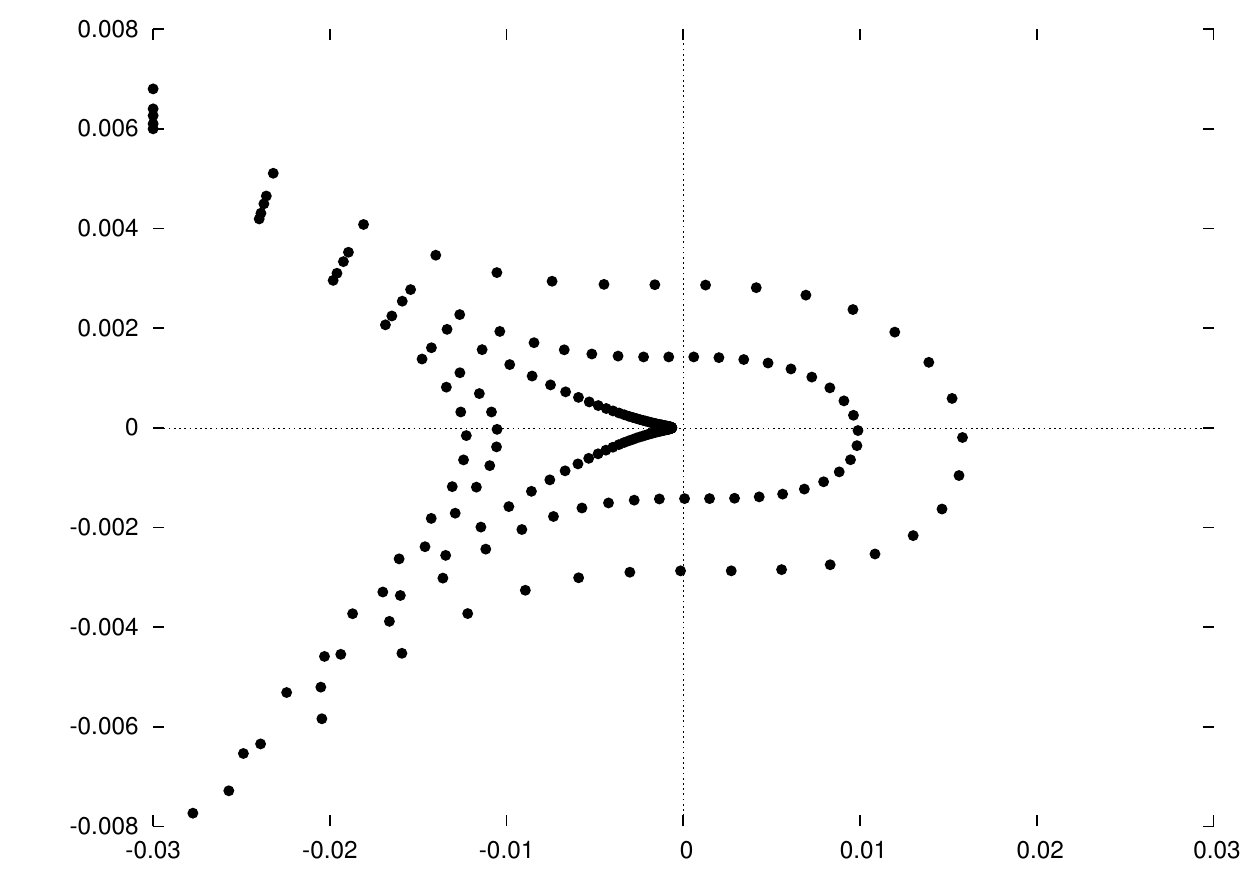}
\else
\vspace{1.62in}
\fi
}
\hfill
\subfigure[RICs and one unbounded orbit for $\mu = 0.1$]{
\label{fig:AfterParabolic0}
\iffigures
\includegraphics[height=1.62in]{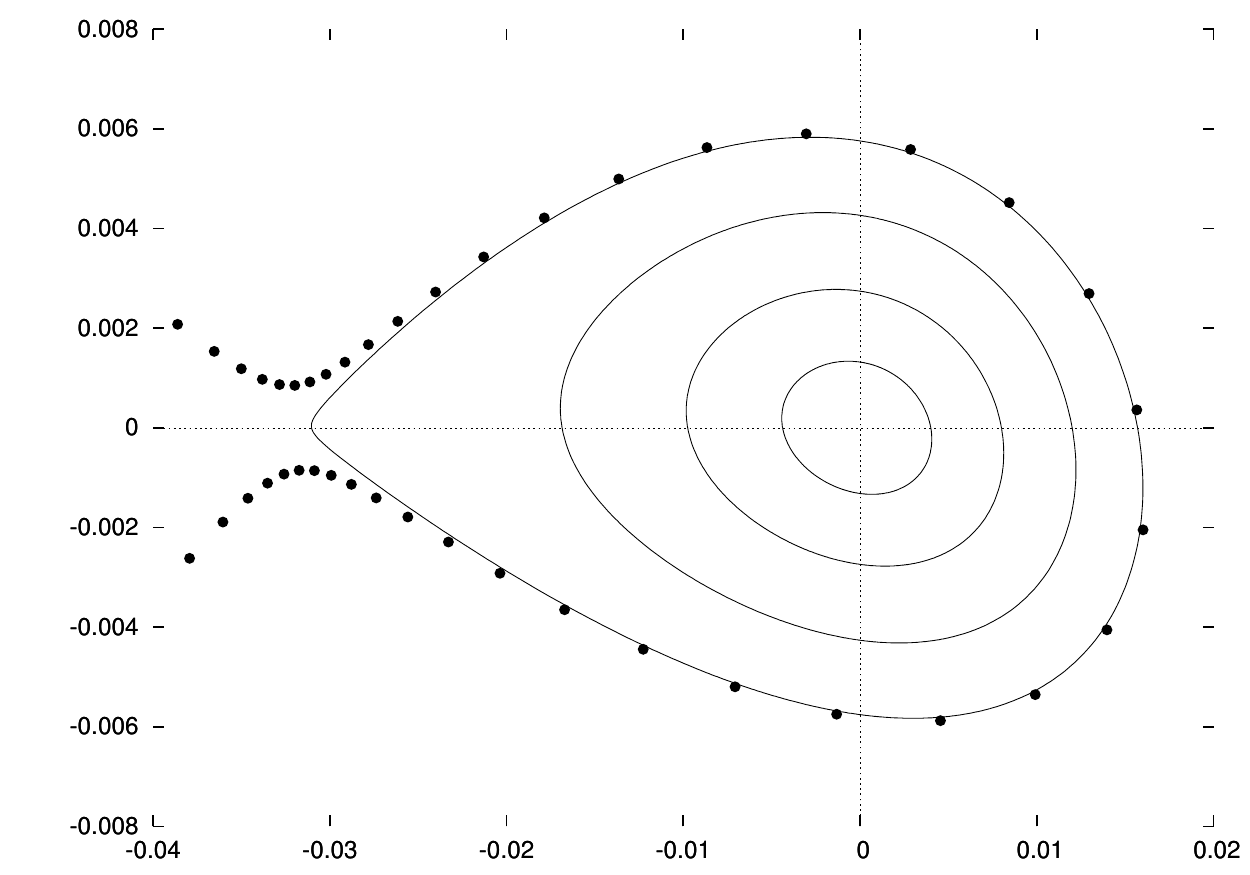}
\else
\vspace{1.62in}
\fi
}
\subfigure[RICs and one unbounded orbit for $\mu = 1$]{
\label{fig:mu_1}
\iffigures
\includegraphics[height=1.62in]{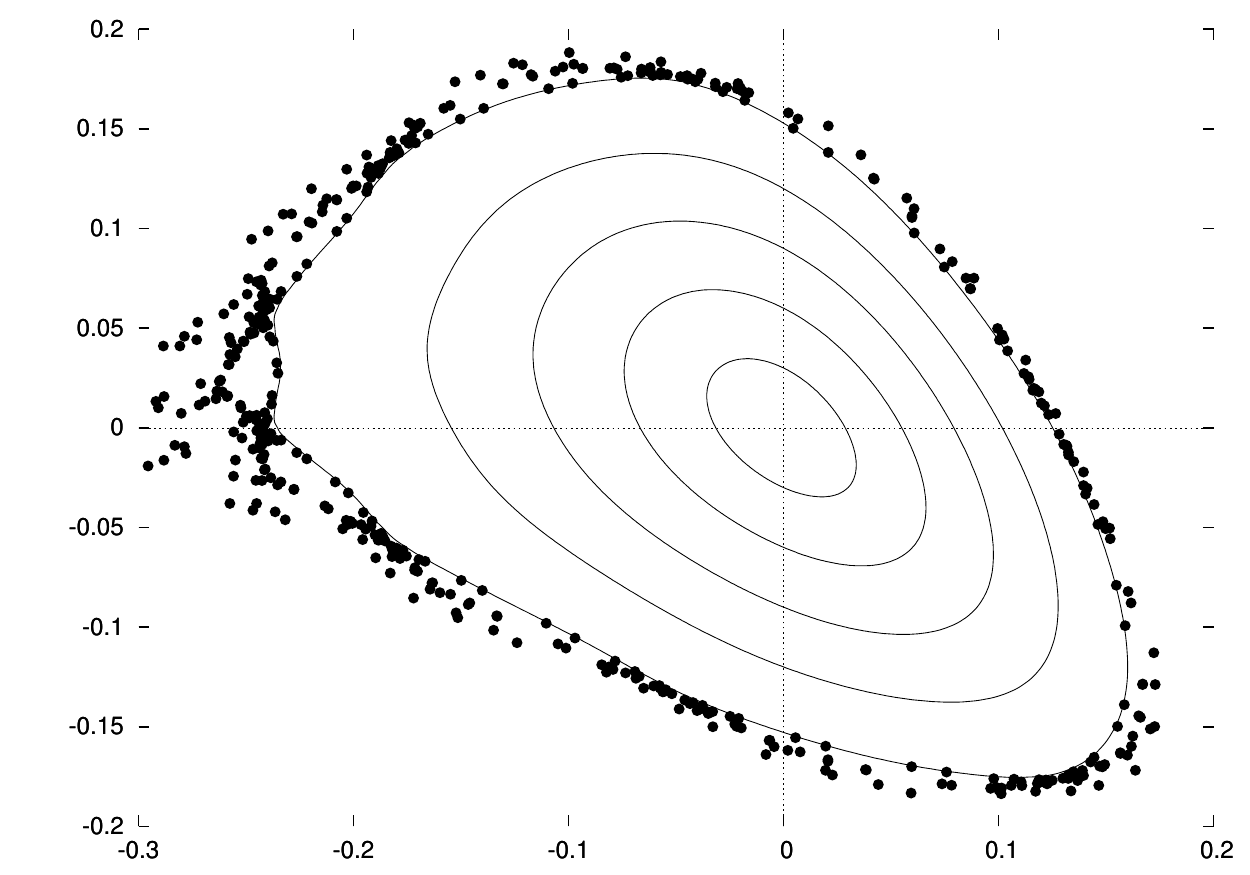}
\else
\vspace{1.62in}
\fi
}
\hfill
\subfigure[RICs and one unbounded orbit for $\mu = 2$]{
\label{fig:FourthOrderResonance}
\iffigures
\includegraphics[height=1.62in]{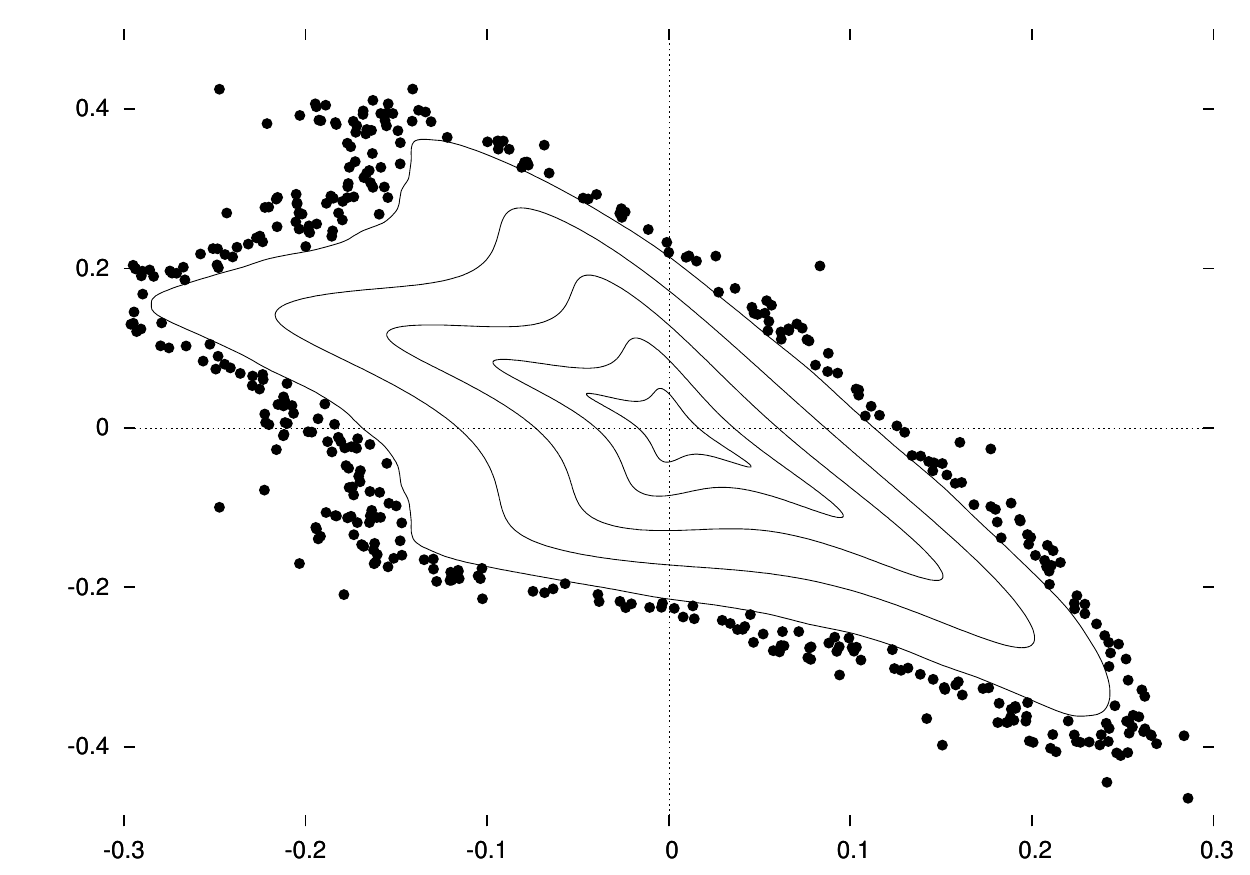}
\else
\vspace{1.62in}
\fi
}
\hfill
\subfigure[RICs and one unbounded orbit for $\mu = \mu_\rmr \simeq 2.53$]{
\label{fig:RootOfTheTwistCoefficient}
\iffigures
\includegraphics[height=1.62in]{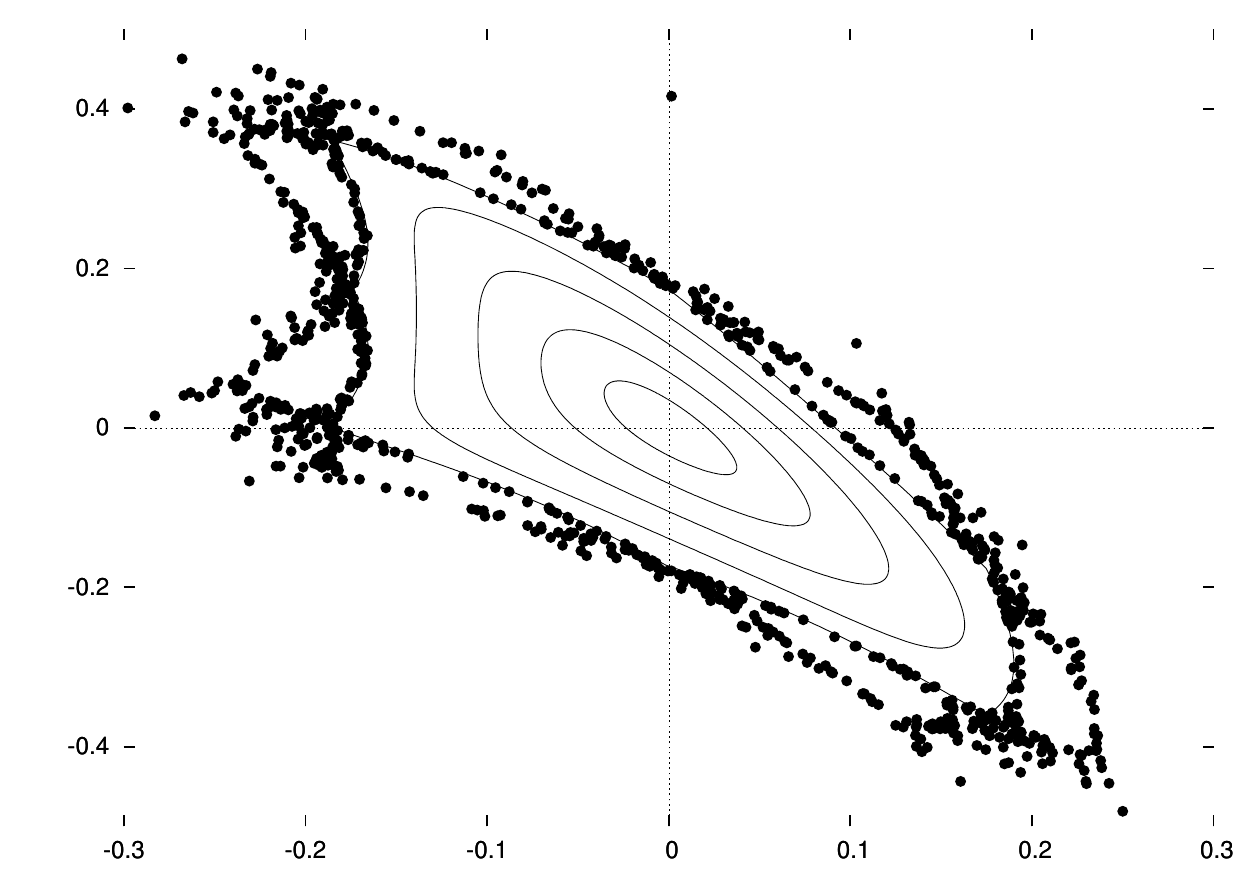}
\else
\vspace{1.62in}
\fi
}
\subfigure[RICs and one unbounded orbit for $\mu = 2.9$]{
\label{fig:BeforeThirdOrderResonance}
\iffigures
\includegraphics[height=1.62in]{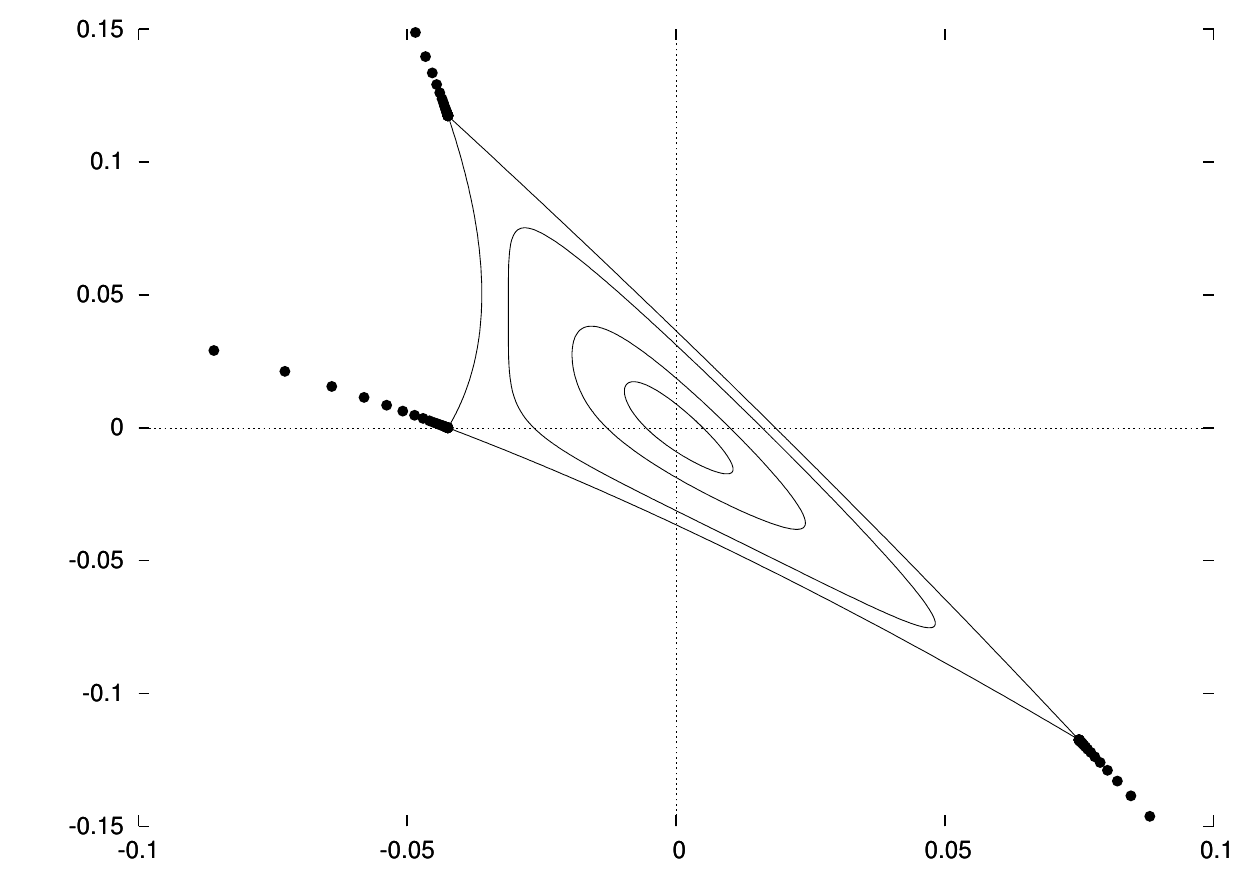}
\else
\vspace{1.62in}
\fi
}
\hfill
\subfigure[One unbounded orbit for $\mu = 3$]{
\label{fig:ThirdOrderResonance}
\iffigures
\includegraphics[height=1.62in]{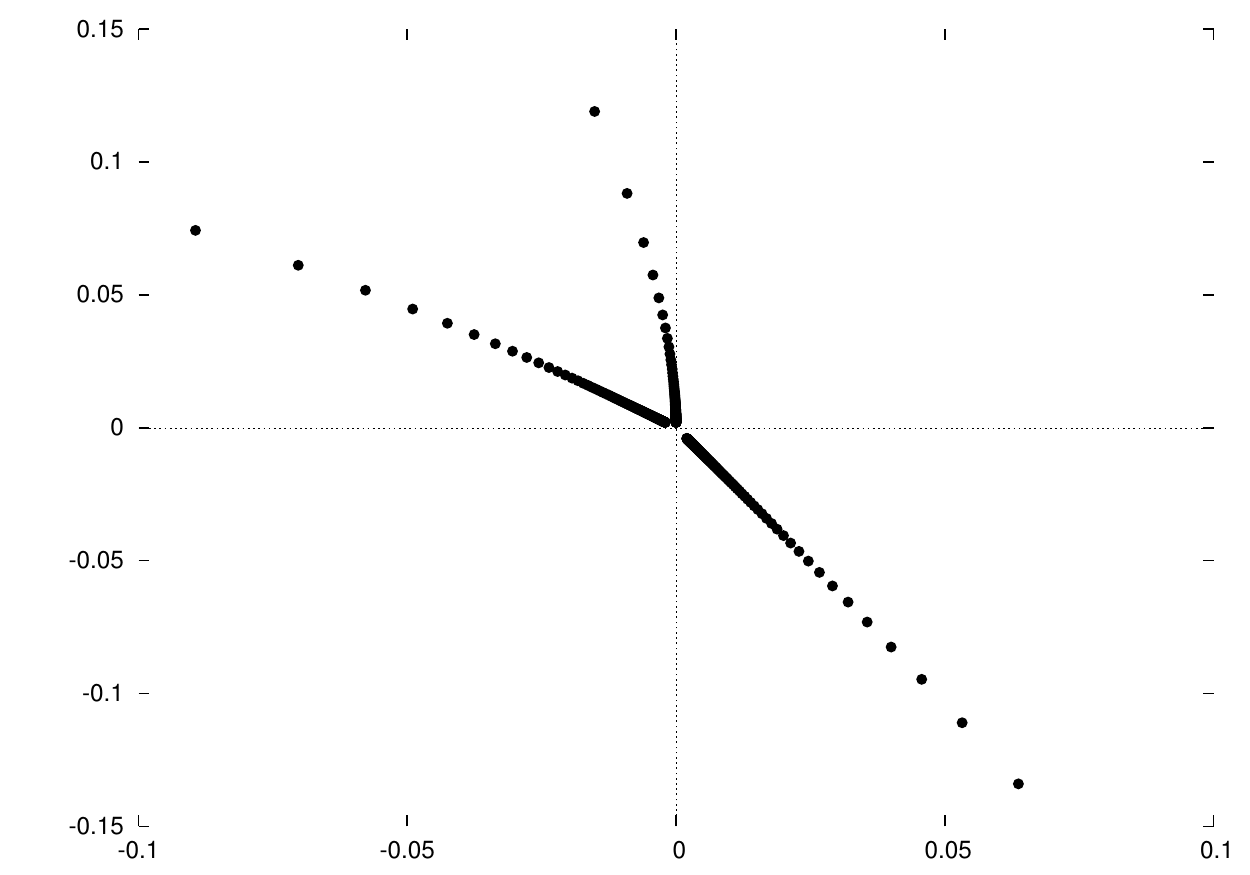}
\else
\vspace{1.62in}
\fi
}
\hfill
\subfigure[RICs and one unbounded orbit for $\mu = 3.1$]{
\label{fig:AfterThirdOrderResonance}
\iffigures
\includegraphics[height=1.62in]{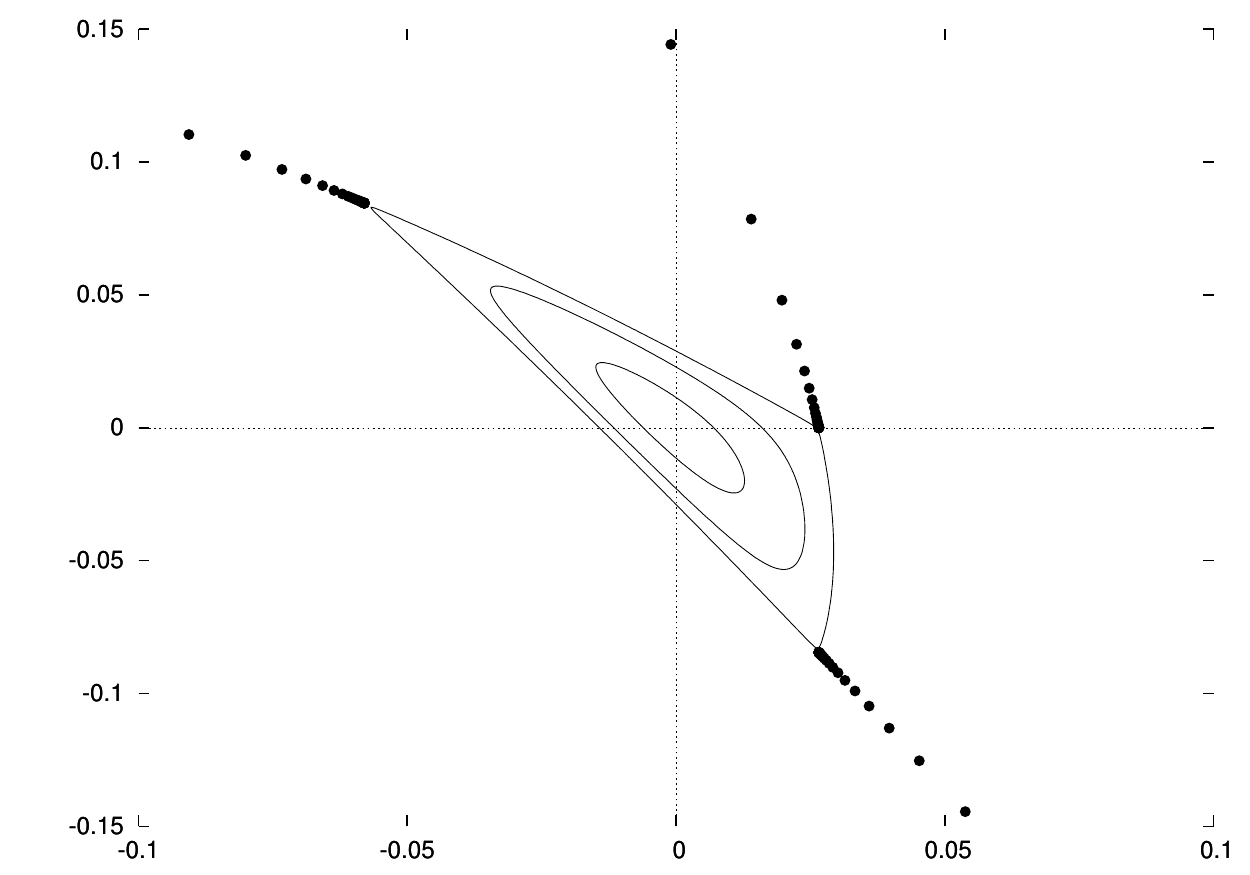}
\else
\vspace{1.62in}
\fi
}
\subfigure[RICs and one unbounded orbit for $\mu = 3.9$]{
\label{fig:BeforeParabolic4}
\iffigures
\includegraphics[height=1.62in]{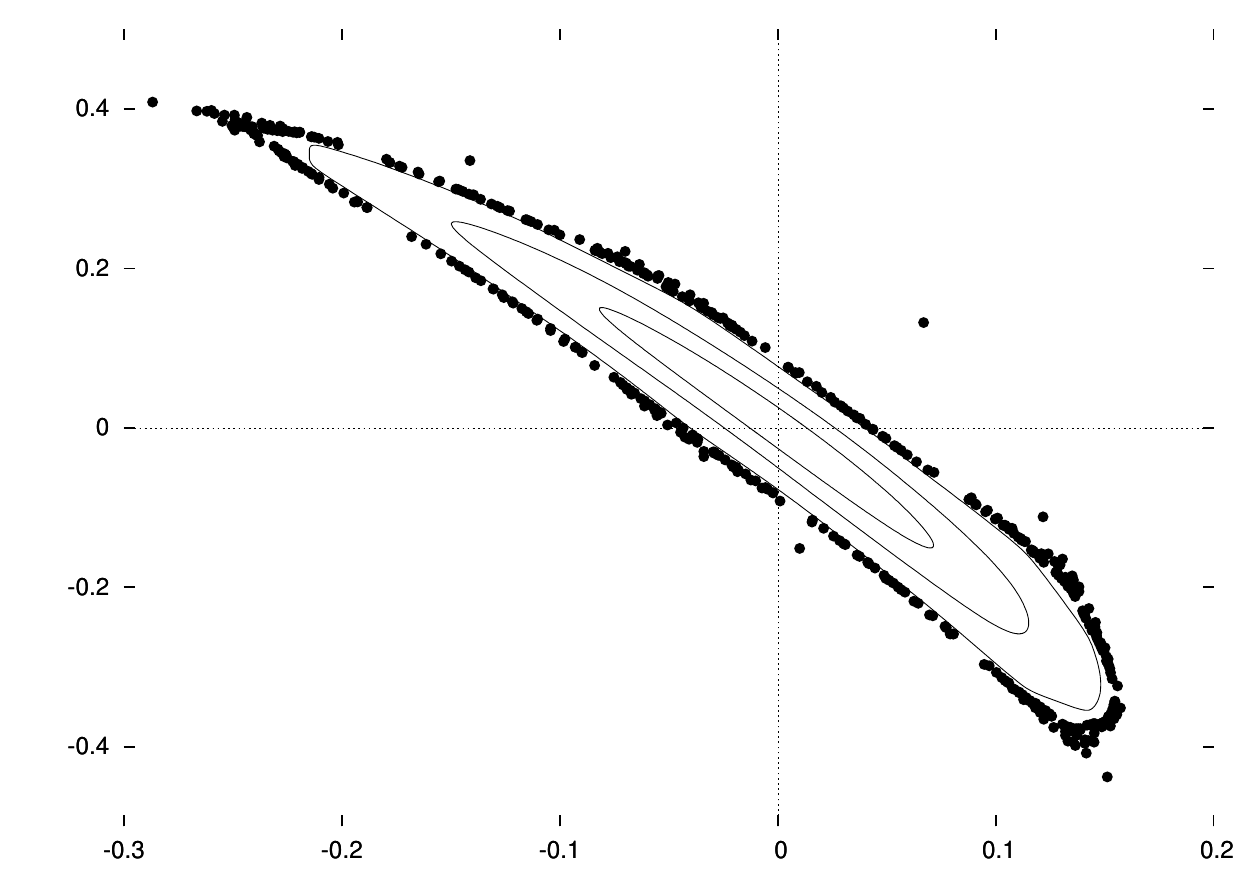}
\else
\vspace{1.62in}
\fi
}
\hfill
\subfigure[RICs and one unbounded orbit for $\mu = 4$]{
\label{fig:Parabolic4}
\iffigures
\includegraphics[height=1.62in]{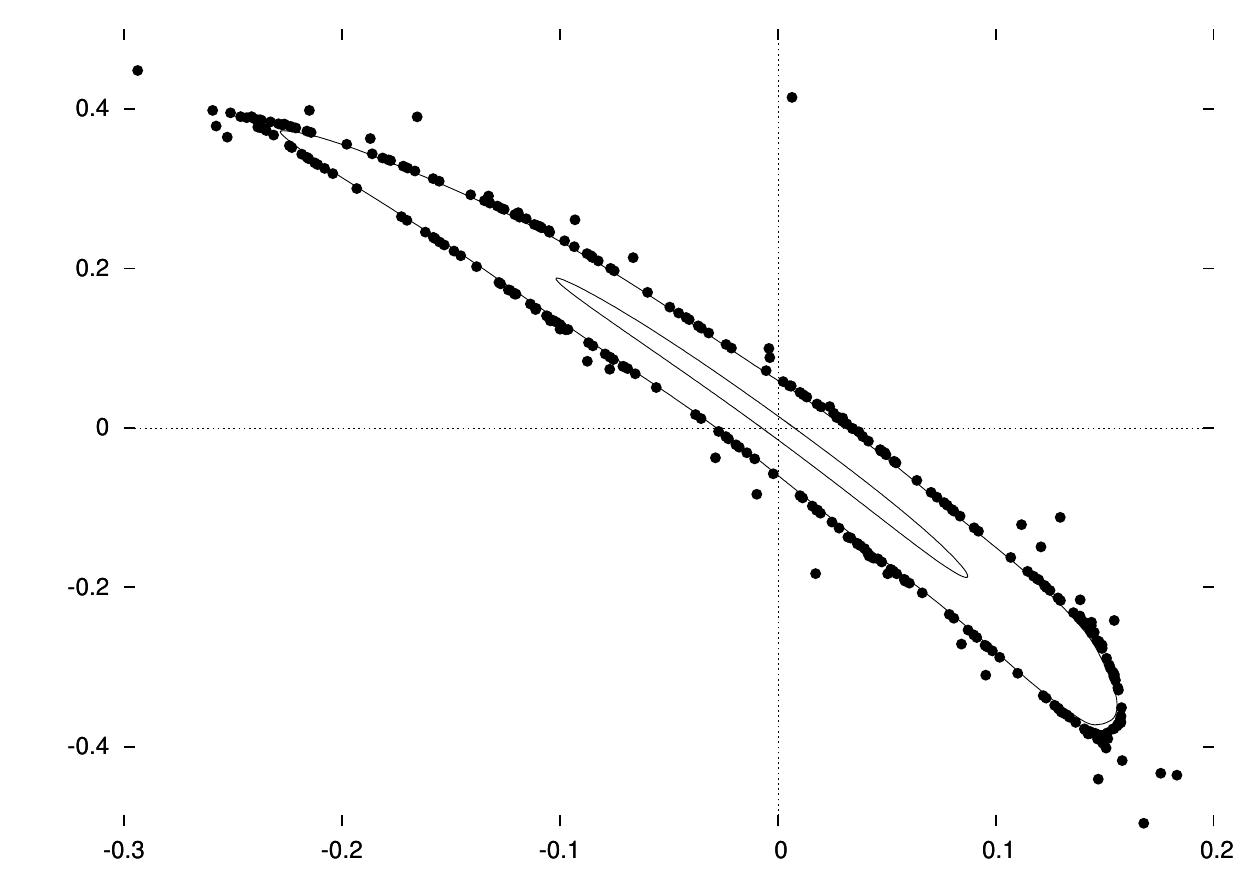}
\else
\vspace{1.62in}
\fi
}
\hfill
\subfigure[One RIC, one unbounded orbit and one
homoclinic (to the origin) orbit for $\mu = 4.05$]{
\label{fig:AfterParabolic4}
\iffigures
\includegraphics[height=1.62in]{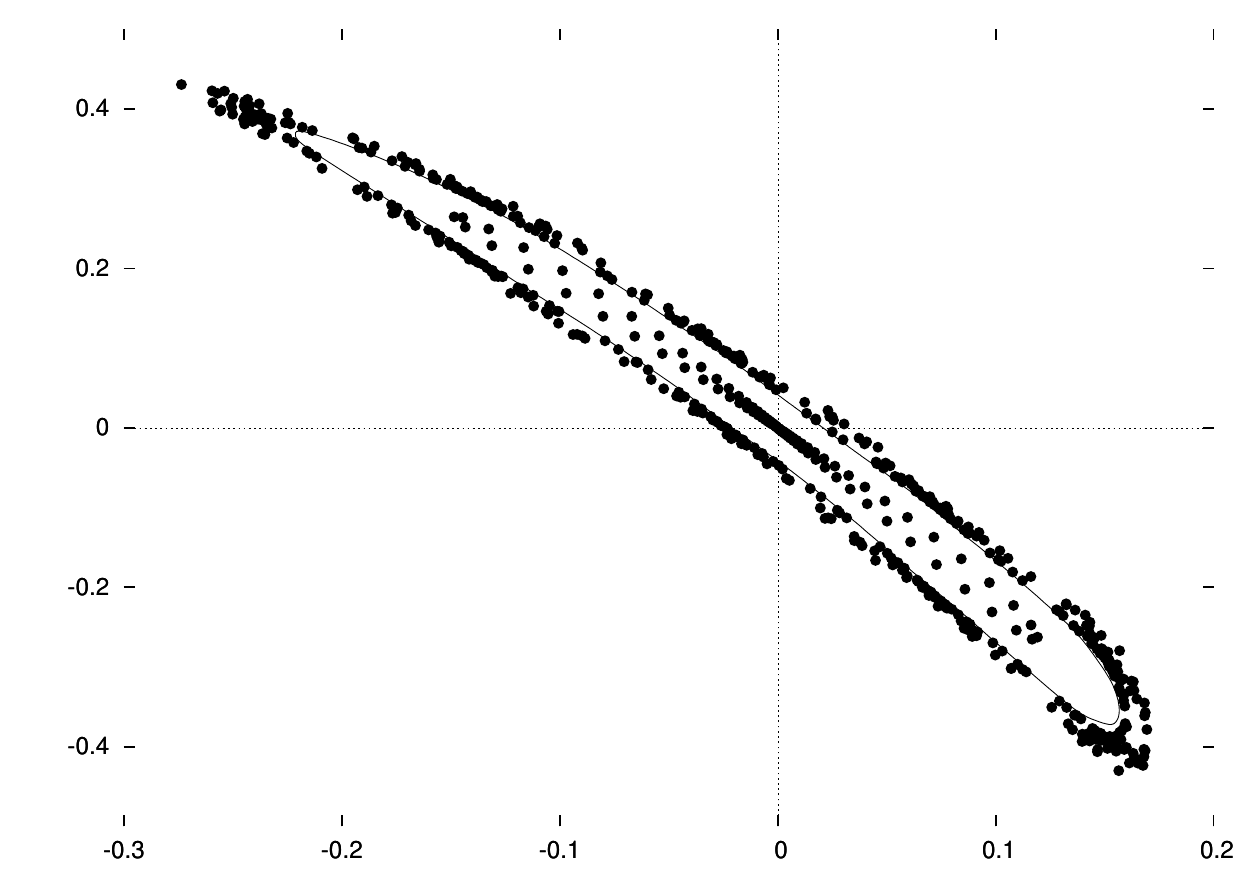}
\else
\vspace{1.62in}
\fi
}
\caption{Race-track microtron phase dynamics in $(\psi,w)$-coordinates
for several values of the parameter $\mu = 2\pi \tan \phi_\rms$.
In the electronic version one can magnify the plots to check
details of the invariant curves: proximity to the unbounded orbits,
regularity at their ``corners'', etcetera.
This also applies to many other figures.}
\label{fig:SummaryOfDynamics}
\end{figure*}

Local stability implies global stability,
but local instability does not imply global instability.
The case $\mu \gtrsim 4$ is a sample of this claim.
See Fig.~\ref{fig:AfterParabolic4}.
Our map~(\ref{eq:f}) shows the four possible combinations of
local stability/instability and linear stability/instability.

The synchronous trajectory is \emph{hyperbolic},
\emph{parabolic}, or \emph{elliptic} when the eigenvalues of
the matrix $M_\rms$ are real of modulus different from one, real of modulus
equal to one, or non-real of modulus equal to one, respectively.
We note that $\det[M_\rms] = 1$, since $f$ is an area preserving map.
Therefore, the behavior in the linear approximation of
the synchronous trajectory only depends
on the trace $T_\rms = \Trace[M_\rms]$.
Concretely, it is hyperbolic, parabolic or elliptic if and only if
$|T_\rms| > 2$, $|T_\rms| = 2$, or $|T_\rms| < 2$, respectively.

The linear type and the local stability are related as follows.
The hyperbolic type implies local instability,
whereas the elliptic type is generically locally stable,
but local instability may take place in degenerate cases~\cite{SiegelMoser1995}.
The parabolic type is the hardest one, but it can be
studied using results from~\cite{LeviCivita1901,Simo1980,Simo1982}.

The linear part of the map~(\ref{eq:f}) at the fixed point $p_\rms$ is
\[
M_\rms =
\left.
\frac{\partial(\psi_1,w_1)}{\partial(\psi,w)}
\right|_{(\psi,w) = (\psi_\rms,w_\rms)} =
\left(
\begin{array}{rc}
1 & 1 \\ -\mu & 1 - \mu
\end{array}
\right),
\]
so $T_\rms = \Trace[M_\rms] = 2 - \mu = 2 - 2\pi \tan \phi_\rms$.
Thus, the linear type of $p_\rms$ depends on the parameter $\mu$
as follows.
It is hyperbolic, parabolic or elliptic if and only if
$\mu \not \in [0,4]$, $\mu \in \{0,4\}$, and $\mu \in (0,4)$,
respectively.

If $\mu \in (0,4)$,
then the eigenvalues of $M_\rms$ have the form
\[
\lambda_\rms = \rme^{\rmi \theta},\qquad
\lambda_\rms^{-1} = \bar{\lambda_\rms} = \rme^{-\rmi \theta},
\]
for some $\theta \in (0,\pi)$ such that
\begin{equation}\label{eq:RotationNumber}
\cos \theta =
(\lambda_\rms + \bar{\lambda_\rms})/2 =
\Trace[M_\rms] /2 = 1 -\mu/2 =
1 - \pi \tan \phi_\rms.
\end{equation}
This means that the linear dynamics around $p_\rms$ is
conjugated to a rotation by angle $\theta$.

If $p_n = f^n(p)$ is an orbit of the map,
then we say that it is \emph{unbounded} when $(p_n)_{n \in \Zset}$
is unbounded in $\Tset \times \Rset$,
and we say that it is \emph{homoclinic} to the fixed
point $p_\rms$ when $\lim_{n \to \pm \infty} p_n = p_\rms$,
but $p_n \neq p_\rms$.
Any unbounded orbit is placed outside the acceptance.

The map~(\ref{eq:f}) has two fixed points:
$p_\rms = (0,0)$ and
\[
p_\rmh = (\psi_\rmh,w_\rmh) := (-2\phi_\rms,0).
\]
The linear part of the map~(\ref{eq:f}) at the fixed point $p_\rmh$ is
\[
M_\rmh =
\left.
\frac{\partial(\psi_1,w_1)}{\partial(\psi,w)}
\right|_{(\psi,w) = (\psi_\rmh,w_\rmh)} =
\left(
\begin{array}{cc}
1 & 1 \\ \mu & 1 + \mu
\end{array}
\right),
\]
so $T_\rmh = \Trace[M_\rmh] = 2 + \mu$.
Hence, $p_\rmh$ is hyperbolic for any $\mu > 0$.
Indeed, if $\mu > 0$, then the eigenvalues of $M_\rmh$ have the form
\[
\lambda_\rmh = \rme^h, \qquad \lambda^{-1}_\rmh = \rme^{-h},
\]
for some $h > 0$ such that
\[
\cosh h = (\lambda_\rmh + \lambda^{-1}_\rmh)/2 =
\Trace[M_\rmh]/2 = 1+ \mu/2 =
1 + \pi \tan \phi_\rms.
\]
This means that the linear dynamics around $p_\rmh$
expands the unstable direction by a factor $\lambda_\rmh$ and
contracts the stable direction by a factor $\lambda^{-1}_\rmh$.

The quantity $\theta/2\pi$ is called \emph{tune} in
Accelerator Physics~\cite{Rand1984}.
Nevertheless, following a standard terminology in Dynamical Systems,
we will say that $\theta/2\pi$ is the \emph{rotation number} of
the elliptic fixed point $p_\rms$.
The elliptic point $p_\rms$ is called \emph{$(m,n)$-resonant} when
$\theta = 2\pi m/n$ for some relatively prime integers $m$ and $n$
such that $1 \le m \le n/2$, which implies that $\lambda_\rms^n = 1$.
Besides, $n$ is the \emph{order} of the resonance.
The quantity $h$ is the \emph{characteristic exponent} of
the hyperbolic fixed point $p_\rmh$.

The quantities $\theta$ and $h$ carry the main local information around
the fixed points when $\mu \in (0,4)$.

\emph{Moser twist theorem} is the standard tool to prove
that fixed points of analytic area preserving maps are locally
stable~\cite{SiegelMoser1995}.
If the dynamics around the fixed point satisfies some
\emph{twist condition}, Moser twist theorem implies that
any neighborhood of the fixed point contains infinitely many
closed invariant curves.
We say that such curves are \emph{rotational
invariant curves (RICs)} since they surround the fixed point and
their internal dynamics is conjugated to a rigid rotation.
The local stability follows from the fact that the domain
enclosed by any RIC is invariant.

We will prove the following characterization
in Section~\ref{Sec:LocalStability}.

\begin{thm}\label{thm:LocalStability}
The synchronous trajectory is locally stable
if and only if $\mu \in (0,4] \setminus \{3\}$.
\end{thm}

The range of local stability in terms of $\theta$ and $h$
is given in Table~\ref{tab:Relations}.

\begin{table}
\centering
\begin{tabular}{ll}
\hline
Relation & Local stability \\
\hline
$\mu = 2\pi \tan \phi_\rms$ & $\mu \in (0,4]\setminus\{3\}$ \\
$\cos \theta = 1 - \pi \tan \phi_\rms$ & $\theta \in (0,\pi]\setminus \{ 2\pi/3\}$ \\
$\cosh h = 1 + \pi \tan \phi_\rms$ & $h \in (0,h_\rmp] \setminus \{ h_\rmu \}$ \\
\hline
\end{tabular}
\caption{\label{tab:Relations}Relations among
the synchronous phase $\phi_\rms \in \Tset$,
the parameter $\mu \in \Rset$,
the rotation number $\theta/2\pi$,
and the characteristic exponent $h>0$.
Here, the values $h_\rmp$ and $h_\rmu$ are
given by $\cosh h_\rmp = 3$ and $\cosh h_\rmu = 5/2$.}
\end{table}

Once we have a clear understanding of the local dynamics,
we focus on the global picture, which rises several questions.
Which is the ``last'' RIC (LRIC)?
How can we compute it?
Does it coincide with the border of the connected component $\mathcal{D}$?

If the border $\mathcal{C} = \partial \mathcal{D}$ is an analytic curve
and its internal dynamics is conjugated to a rigid rotation,
then $\mathcal{C}$ is the LRIC that we are looking for and we can numerically
compute it by using the algorithm described in~\cite{SimoTreschev1998}.
These two hypotheses about $\mathcal{C}$ are not so restrictive
as they may look.
Indeed, all RICs obtained from the Moser's twist theorem satisfy them.
This means that even when $\mathcal{C} = \partial \mathcal{D}$ does not
satisfy them, probably there are many Moser-like RICs close to $\mathcal{C}$,
which are the objects that our computations will find.

A summary of the phenomena that take place when $\mu$
moves is displayed in Fig.~\ref{fig:SummaryOfDynamics}.
If there exist RICs around the origin,
we plot some of them, the last one (the LRIC),
and one unbounded orbit very close to the LRIC.
If the origin is globally unstable, we plot suitable unbounded orbits
to highlight it.
The invariant curves (ICs) displayed in Fig.~\ref{fig:BeforeParabolic0}
are no RICs for the synchronous trajectory,
since they do not surround the origin.

Fig.~\ref{fig:SummaryOfDynamics} shows the basic skeleton from which
an expert reader can deduce the main dynamical changes.
Let us describe them.

The elliptic fixed point $p_\rms$ and the hyperbolic fixed point $p_\rmh$
merge in a parabolic point as $\mu \to 0^+$;
see Figs.~\ref{fig:BeforeParabolic0}--\ref{fig:AfterParabolic0}.
This scenario corresponds to a
saddle-center bifurcation~\cite{Broer_etal1996,Gelfreich2000}.
The stable and unstable invariant curves (called \emph{separatrices})
of the hyperbolic fixed point seem to form a small loop
around the elliptic fixed point when $0 < \mu \ll 1$.
The LRIC around the elliptic point is exponentially close
(in the parameter $\mu$) to the separatrices.
See Fig.~\ref{fig:AfterParabolic0}.
We will use this fact to see that
\begin{equation}\label{eq:AsymptoticSaddleCenterBifurcation}
|\mathcal{A}_\mu|,|\mathcal{D}_\mu| = 6 \mu^{5/2}/5\pi^2 + \Order(\mu^3),
\qquad \mbox{as $\mu \to 0^+$}.
\end{equation}
Therefore, both regions are almost the same
when $\mu \to 0 ^+$.
(If $\mathcal{R}$ is a subset of $\Tset \times \Rset$,
 then $|\mathcal{R}|$ denotes its area.)

The LRIC grows in size, changes its shape, and moves away from
the separatrices as $\mu$ increases; see Fig.~\ref{fig:mu_1}.
We reach the fourth order resonance $\theta = \pi/2$ at $\mu = 2$.
The elliptic point is locally stable at this resonance,
and the RICs near it look like a Latin cross,
whose arms have a width of the order of the square of its length;
see Fig.~\ref{fig:FourthOrderResonance}.
The H\'enon map displays exactly the same behavior at the
fourth order resonance~\cite{Simo1980}.

The twist coefficient (also called first Birkhoff coefficient)
goes from negative to positive at $\mu = \mu_\rmr \simeq 2.537706$
and from positive to negative at $\mu = 3$,
where the third order resonance $\theta = 2\pi/3$ takes place.
This is one of the two typical behaviors for the twist
coefficient described in~\cite{Moeckel1990}.
It means that the orbits around the origin rotate slower
(resp., faster) as they go away from the origin
when $\mu \in (0,\mu_\rmr)\cup (3,4)$ (resp., $\mu \in (\mu_\rmr,3)$).
Figs.~\ref{fig:BeforeThirdOrderResonance}--\ref{fig:AfterThirdOrderResonance}
show that behavior,
since we see a three-periodic orbit near the elliptic point
for both $\mu \lesssim 3$ and $\mu \gtrsim 3$.

The elliptic point is locally and globally unstable at
the third order resonance; see Fig.~\ref{fig:ThirdOrderResonance}.
If $\mu = 3 + \epsilon$ with $0 < |\epsilon| \ll 1$,
the LRIC is approximately a triangle of vertices
\begin{equation}\label{eq:TriangleVertices}
\frac{\epsilon}{\pi}(1,0),\qquad \frac{\epsilon}{\pi}(1,-3),\qquad
\frac{\epsilon}{\pi}(-2,3).
\end{equation}
In particular,
\begin{equation}\label{eq:AsymptoticThirdOrderResonace}
|\mathcal{D}_{3+\epsilon}| = 9\epsilon^2/2\pi^2 + \Order(\epsilon^3),
\qquad \mbox{as $\epsilon \to 0$}.
\end{equation}

The LRIC grows in size and its triangular shape changes to
a ``banana-like'' shape as the parameter moves from the three order resonance
at $\mu = 3$ to the second order resonance at $\mu = 4$.
The fixed point $p_\rms$ is locally stable at this last resonance,
and the RICs near it looks like ``bananas'' whose width is of
the order of the square of its length. See Fig.~\ref{fig:Parabolic4}.

Finally, a period-doubling bifurcation takes place after
the second order resonance.
To be precise, the point $p_\rms$ becomes a saddle for $\mu > 4$,
being parabolic with reflection at $\mu = 4$.
One then expects to find an elliptic two-periodic orbit for $\mu > 4$.
The \emph{separatrices} of the hyperbolic fixed point $p_\rms$
form two small loops around the two elliptic two-periodic points.
Nevertheless, the saddle $p_\rms$ remains globally stable when $\mu \gtrsim 4$,
since some of the RICs around $p_\rms$ still persist.
See Fig.~\ref{fig:AfterParabolic4}.

In spite of all the above comments, the expulsion of some
resonance is a more relevant phenomenon than any resonance at the elliptic
fixed point, but the third order resonance.
The reason is that only the first phenomenon changes
drastically $\mathcal{A}$ and $\mathcal{D}$,
since the second one takes place inside $\mathcal{D}$.
We will check that $|\mathcal{A}|$ and $|\mathcal{D}|$ experiment a
jump for each primary resonance that is thrown away.
Besides,
these primary jumps are smaller in $|\mathcal{A}|$ than in $|\mathcal{D}|$.
However, $|\mathcal{A}|$ displays many more jumps, the secondary ones,
which are related to resonances inside elliptic islands
that are already outside $\mathcal{D}$.
We will visualize such properties in Section~\ref{Sec:GlobalStability}.

\section{Local stability of the synchronous trajectory}
\label{Sec:LocalStability}

This section has a purely local character.
It contains the proof of Theorem~\ref{thm:LocalStability}.
The hyperbolic type always implies local instability,
so we only study the range $0 \le \mu \le 4$.

\subsection{Local instability in the saddle-center bifurcation}
\label{Ssec:SaddleCenterBifurcation}

Set $\mu = 0$.
Then the map~(\ref{eq:f}) has the form
\[
\left\{
\begin{array}{ccl}
w_1    & = & w - \pi (\psi+w)^2 + \Order_4(\psi,w), \\
\psi_1 & = & \psi + w.
\end{array}
\right.
\]
On the other hand,
the Levi-Civita criterion~\cite{LeviCivita1901} implies that
the origin is unstable under any analytic map of the form
\[
\left\{
\begin{array}{ccl}
w_1 & = & w + \Order_2(\psi,w), \\
\psi_1 & = & \psi + w + \Order_2(\psi,w),
\end{array}
\right.
\]
with $\frac{\partial^2 w_1}{\partial \psi^2}(0,0) \neq 0$.
In our concrete map,
$\frac{\partial^2 w_1}{\partial \psi^2} (0,0) = -2\pi \neq 0$.
Thus, the synchronous trajectory is locally unstable if $\mu = 0$.

\subsection{Local stability in the elliptic case: Generic values}
\label{Ssec:LocalStabilityElliptic}

Let $\mu \in (0,4)$.
Then the point $p_\rms$ is elliptic with rotation number $\theta/2\pi$.
The relation between $\theta$ and $\mu$ is given in Table~\ref{tab:Relations}.

First, we bring the linear part of map~(\ref{eq:f}) around
the elliptic fixed point $p_\rms$ into its real Jordan normal form
\[
R_\theta = 
\left( \begin{array}{rr}
\cos \theta & -\sin \theta \\ \sin \theta & \cos \theta
\end{array} \right) 
\]
by means of an area preserving linear change $(\psi,w) \mapsto (x,y)$.
We write our map in these new variables as
\[
\left( \begin{array}{c} x_1 \\ y_1 \end{array} \right) =
R_\theta 
\left( \begin{array}{c} x \\ y \end{array} \right) + \Order_2(x,y).
\]

Then we introduce the complex variables $z = x + \rmi y$
and $\bar{z} = x - \rmi y$.
The map reads in these variables as
\begin{equation}\label{eq:ComplexVariable}
z_1 = \lambda_\rms z + \Order(|z|^2) = 
\lambda_\rms \left(
z + \sum_{\substack{2 \le j+k \le 3 \\ j,k \ge 0}} c_{jk} z^j \bar{z}^k
\right) + \Order(|z|^4),
\end{equation}
for some complex coefficients $c_{jk}$,
which depend on the original parameter $\mu$.
We recall that $\lambda_\rms = \rme^{\rmi \theta}$.

If $\lambda_\rms^n = 1$ for some integer $n$ such that $1 \le n \le 4$,
then we say that the elliptic point $p_\rms$ is \emph{strongly resonant}.
It is known~\cite{SiegelMoser1995} that if the elliptic point $p_\rms$
is not strongly resonant, then there exists an area preserving polynomial
change of variables which brings the map~(\ref{eq:ComplexVariable})
into its third order \emph{Birkhoff normal form}
\begin{equation}\label{eq:BirkhoffNormalForm}
z_1 = \lambda_\rms \left(z + \rmi \tau z^2 \bar{z} \right) + \Order(|z|^4)
= \rme^{\rmi (\theta + \tau |z|^2)} z + \Order(|z|^4),
\end{equation}
for some real coefficient $\tau$,
called \emph{first Birkhoff coefficient} or \emph{twist coefficient}.
It turns out that
\[
\tau = \Im c_{21} +
\frac{2(4 \cos \theta +1)\sin\theta}
     {(\cos \theta -1)(2\cos \theta + 1)}|c_{20}|^2,
\]
where $\Im c_{21}$ denotes the imaginary part of $c_{21}$,
see~\cite{Moeckel1990,KamphorstPinto2005}.
The coefficients $c_{21}$ and $c_{20}$ also depend on $\theta$,
via the parameter $\mu$.
After some tedious computations~\cite{Larreal2012},
one gets the expression
\begin{equation}\label{eq:TwistCoefficient}
\tau(\theta) =
\frac{2 \cos^3 \theta - 3 \cos^2 \theta + 4\pi^2 \cos \theta + 1 + \pi^2}
{(\cos \theta - 1)(2\cos \theta + 1)\sin^2 \theta}.
\end{equation}
We observe that $\lim_{\theta \to 0^+} \theta^4 \tau(\theta)$,
$\lim_{\theta \to 2\pi/3} (\theta-2\pi/3)\tau(\theta)$, and
$\lim_{\theta \to \pi^-} (\theta-\pi)^2 \tau(\theta)$ exist and are negative.
That is, $\tau(\theta)$ has a negative fourth-order pole at $\theta = 0$,
a negative simple pole at $\theta = 2\pi/3$,
and a negative second-order pole at $\theta = \pi$.
These properties agree with the generic behavior of twist coefficients
of families of area-preserving maps described in~\cite{Moeckel1990}.

\begin{figure}
\iffigures
\centering
\psfrag{0}{\footnotesize{$0$}}
\psfrag{pi}{\footnotesize{$\pi$}}
\psfrag{pi/3}{\footnotesize{$\pi/3$}}
\psfrag{2*pi/3}{\footnotesize{$2\pi/3$}}
\includegraphics[height=2.4in]{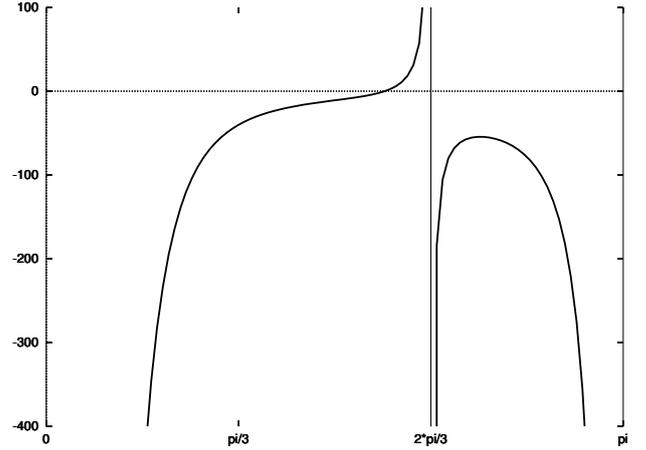}
\else
\vspace{2.4in}
\fi
\caption{The twist coefficient $\tau$ versus $\theta$.}
\label{fig:twist_coefficient}
\end{figure}

Consequently, we deduce that $\tau(\theta)$ has at least one root
in the interval $(0,2\pi/3)$.
Let us prove that it has no more real roots.
The numerator in~(\ref{eq:TwistCoefficient}) is the polynomial
$a \zeta^3 + b \zeta^2 + c \zeta + d$ in the variable $\zeta = \cos \theta$,
with $a = 2$, $b = -3$, $c = 4\pi^2$, and $d = 1 + \pi^2$.
This cubic polynomial has one real root and two complex conjugated roots,
because its discriminant is negative:
\[
\Delta =
b^2 c^2 - 4 a c^3 - 4 b^3 d - 27 a^2 d^2 + 18 a b c d \simeq
-536135.
\]
This means that $\tau(\theta)$ has just one root $\theta_\rmr$
in the interval $(0,\pi)$.
It was numerically computed in~\cite{Larreal2012} that
\[
\theta_\rmr \simeq 1.842998343412199023198246043 \in (0,2\pi/3).
\]
We have plotted the twist coefficient $\tau$ versus $\theta$
in Fig.~\ref{fig:twist_coefficient}. 

Finally, we recall that the Moser twist theorem~\cite{SiegelMoser1995}
implies that the origin is a stable elliptic fixed point of any analytic map
of the form~(\ref{eq:BirkhoffNormalForm}) when $\tau \neq 0$.
Hence,
the synchronous trajectory is locally stable
when $\theta \in (0,\pi) \setminus \{\pi/2,\theta_\rmr,2\pi/3\}$;
or, equivalently, when $\mu \in (0,4) \setminus \{2,\mu_\rmr,3\}$, where
\[
\mu_\rmr = 2-2\cos \theta_\rmr \simeq  2.537706055658189018165133406.
\]

\subsection{Local stability in the fourth order resonance}
\label{Ssec:LocalStabilityFourthOrder}

Set $\mu = 2$.
Then the map~(\ref{eq:f}) takes the form
\[
\left\{
\begin{array}{ccl}
\psi_1 & = & \psi + w,\\
w_1    & = & -2\psi -w - a(\psi_1),
\end{array}
\right.
\]
where $a(\psi_1) = 2\pi(1 - \cos \psi_1) - 2(\psi_1 - \sin \psi_1) =
\Order(\psi_1^2)$.

We consider the area preserving linear change of variables
\begin{equation}\label{eq:CoordinatesFourthOrder}
x = \psi + w,\qquad y = \psi.
\end{equation}
This change brings the linear part of the map at the origin
into its real Jordan normal form $R_{\pi/2}$.
Indeed, the map becomes
\begin{equation}\label{eq:Map_At_4Resonance}
\left( \begin{array}{c} x_1 \\ y_1 \end{array} \right) =
\left( \begin{array}{c} -y - a(x) \\ x \end{array} \right) =
R_{\pi/2}
\left( \begin{array}{c} x \\ y + a(x) \end{array} \right).
\end{equation}
Let $\sum_{j \ge 2} a_j x^j$ be the Taylor expansion of $a(x)$.
We know from Corollary 4.2 in~\cite{Simo1982} that if $a_2 \neq 0$
and $a_3 \neq 0$, then the origin is locally unstable
under the analytic map~(\ref{eq:Map_At_4Resonance}) if and only if
\[
0 < a_3 \le a_2^2.
\]
In our concrete map, $a_2 = \pi$ and $a_3 = -1/3$.
Therefore, the synchronous trajectory is locally stable when $\mu = 2$,
and we will analytically study the shape of its RICs
in Section~\ref{Ssec:HamiltonianFourthOrder}.

\subsection{Local stability for the root of the twist coefficient}
\label{Ssec:LocalStabilityRoot}

If $\mu = \mu_\rmr$, then the third order Birkhoff normal
form~(\ref{eq:BirkhoffNormalForm}) does not provide much information
because the first Birkhoff coefficient $\tau_1 = \tau$ vanishes.
Thus, we should compute the fifth order Birkhoff normal form
\[
z_1 =  \rme^{\rmi (\theta_\rmr + \tau_1 |z|^2 + \tau_2 |z|^4)} z + \Order(|z|^6),
\]
where $\tau_1 = 0$ and
$\tau_2 \in \Rset$ is the \emph{second Birkhoff coefficient}.
The analytical computation of $\tau_2$ is cumbersome,
so we have taken a simpler numerical approach.
First, we have computed the rotation number $\rho(p)$
of the points of the form $p = (\psi,0)$ using the algorithm described
in subsection~\ref{Ssec:RotationNumber}.
Next, we have checked that
\begin{equation}\label{eq:FlatBehavior}
\rho(\psi,0) = \theta_\rmr/2\pi + \rho_2 \psi^4 + \Order(\psi^5),
\end{equation}
for some non-zero coefficient $\rho_2 \approx -200$,
which implies that $\tau_2 \neq 0$.
Then the Moser twist theorem implies that the origin is a
stable elliptic fixed point~\cite{SiegelMoser1995}.

\begin{figure}
\iffigures
\centering
\includegraphics[height=2.4in]{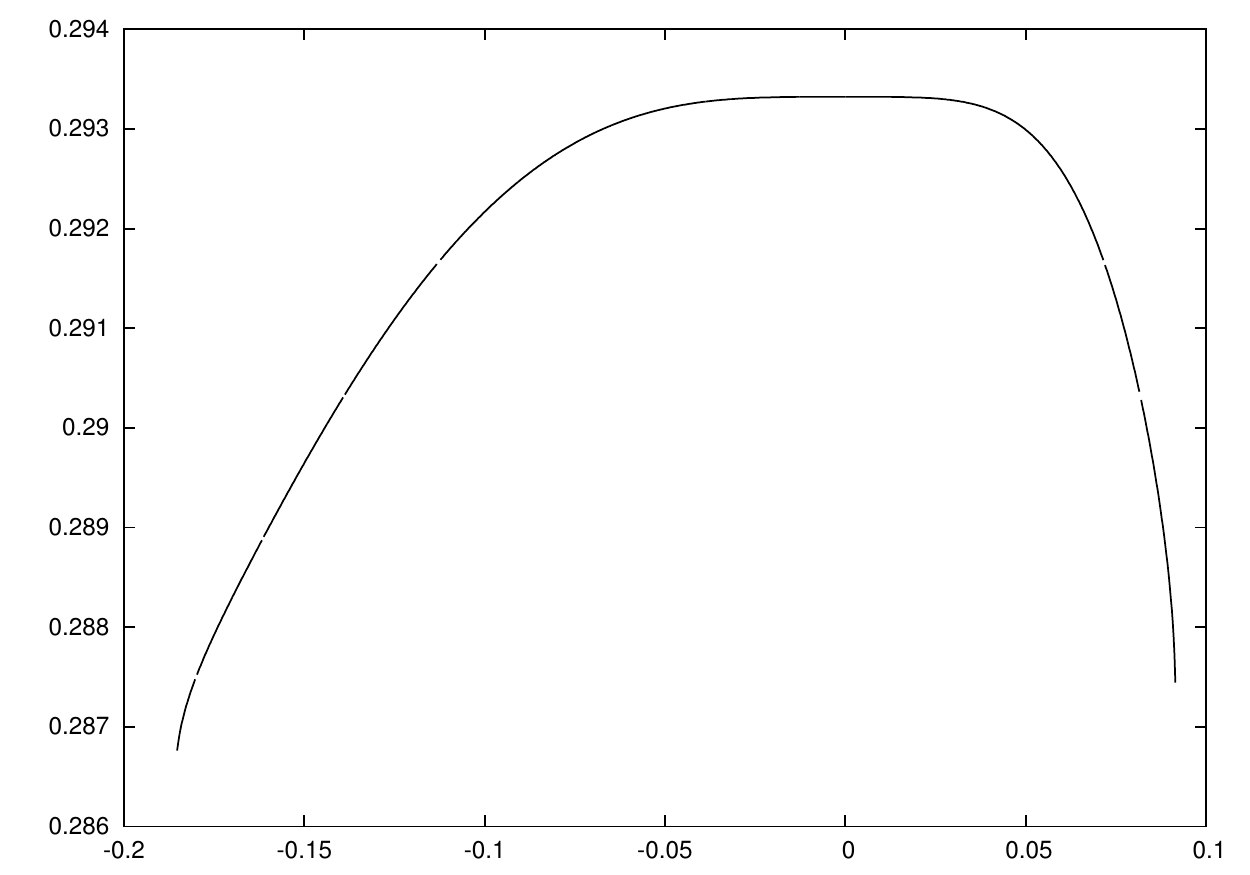}
\else
\vspace{2.4in}
\fi
\caption{The rotation number $\rho(\psi,0)$ versus $\psi$ for $\mu = \mu_\rmr$.}
\label{fig:RotationNumber_mu_root}
\end{figure}

The flat behavior~(\ref{eq:FlatBehavior}) is clearly observed
in Fig.~\ref{fig:RotationNumber_mu_root}.
We will explain the meaning of the small gaps that appear in
Fig.~\ref{fig:RotationNumber_mu_root} at the end of
subsection~\ref{Ssec:RotationNumber}.

\subsection{Local instability in the third order resonance}
\label{Ssec:LocalUnstabilityThirdOrder}

Set $\mu = 3$.
Then the map~(\ref{eq:f}) takes the form
\[
\left\{
\begin{array}{ccl}
\psi_1 & = & \psi + w,\\
w_1    & = & -3\psi -2w - \pi (\psi+w)^2 + \Order_3(\psi,w).
\end{array}
\right.
\]
The map $(\psi_3,w_3) = f^3(\psi,w)$ is close to the identity,
since
\[
\left\{
\begin{array}{ccl}
\psi_3 & = & \psi - \pi( 3 \psi^2 + 2 \psi w ) + \Order_3(\psi,w),\\
w_3    & = & w + \pi (6\psi^2+6\psi w + w^2) + \Order_3(\psi,w).
\end{array}
\right.
\]
We determine the local stability of the origin under the map $f^3$
(and so, under the map $f$) by applying Sim\'o's criterion~\cite{Simo1982}.
That criterion states that the origin is locally stable under an
analytic area preserving map of the form
\[
\left\{
\begin{array}{ccl}
x_1 & = & x + \Order_2(x,y),\\
y_1 & = & y + \Order_2(x,y),
\end{array}
\right.
\]
if and only if $G(x_1,y)$ has a strict extremum at the origin,
where $x_1 y + G(x_1,y)$ is the generating function of the map;
that is, $G(x_1,y)$ is a function determined by the implicit equations
\begin{equation}\label{eq:GeneratingFunctionEquations}
x_1 = x + \frac{\partial G}{\partial y}(x_1,y)
\qquad \mbox{and} \qquad
y_1 = y - \frac{\partial G}{\partial x_1}(x_1,y).
\end{equation}
Thus, we look for a function $G(w_3,\psi) = O_3(w_3,\psi)$ such that
\[
w_3 = w + \frac{\partial G}{\partial \psi}(w_3,\psi) \qquad \mbox{and}\qquad
\psi_3 = \psi - \frac{\partial G}{\partial w_3}(w_3,\psi).
\]
After a straightforward computation, we obtain that
\begin{eqnarray*}
G(w_3,\psi) & = & \pi (2\psi^3 + 3\psi^2 w_3 + \psi w_3^2)
+ \Order_4(\psi,w_3) \\
& = & \pi \psi (\psi + w_3) (2\psi + w_3) + \Order_4(\psi,w_3).
\end{eqnarray*}
This function has no strict extremum at the origin.
Therefore, the synchronous trajectory is unstable when $\mu = 3$.
We will give more details about the dynamics around the origin
near this third order resonance in Section~\ref{Ssec:HamiltonianThirdOrder}.

\subsection{Local stability in the second-order resonance}
\label{Ssec:SecondOrderResonance}

Set $\mu = 4$. Then the map~(\ref{eq:f}) takes the form
\[
\left\{
\begin{array}{ccl}
\psi_1 & = & \psi + w,\\
w_1    & = & -4\psi - 3 w + b(\psi_1),
\end{array}
\right.
\]
where $b(\psi_1) = 2\pi(\cos \psi_1 - 1) + 4(\psi_1 - \sin \psi_1)  =
\Order(\psi_1^2)$.

We consider the area preserving linear change of variables
\[
x = 2\psi + w,\qquad y = \psi,
\]
which brings the linear part of the map at the origin
into its Jordan normal form.
Indeed, the map becomes
\begin{equation}\label{eq:Map_At_2Resonance}
\left\{
\begin{array}{ccl}
x_1 & = & -x + b(x-y),\\
y_1 & = & x - y.
\end{array}
\right.
\end{equation}
Let $\sum_{j \ge 2} b_j u^j$ be the Taylor expansion of $b(u)$.
We know from Lemma 1.4 in~\cite{Simo1982} that if
$b_2 \neq 0$, $b_3 \neq 0$, and $2 b_3 + b_2^2 > 0$,
then the origin is locally stable under the
analytic map~(\ref{eq:Map_At_2Resonance}).

In our concrete map, $b_2 = -\pi$ and $b_3 = 2/3$.
Therefore, the synchronous trajectory is locally stable when $\mu = 4$.

This completes the proof of Theorem~\ref{thm:LocalStability}.

\section{Hamiltonian approximations}
\label{Sec:HamiltonianApproximations}

This section contains a semi-local study of the RICs around
the synchronous trajectory.
We will approximate them by the level curves of suitable
Hamiltonians in three different scenarios.
We will also prove the asymptotic
formulas~(\ref{eq:AsymptoticSaddleCenterBifurcation})
and~(\ref{eq:AsymptoticThirdOrderResonace}).

\subsection{Near the saddle-center bifurcation}
\label{Ssec:HamiltonianSaddleCenterBifurcation}

Let us study the size and shape of $\mathcal{D}_\mu \subset \mathcal{A}_\mu$
when $\mu \to 0^+$.

The elliptic fixed point $p_\rms = (\psi_\rms,0) = (0,0)$
and the hyperbolic fixed point $p_\rmh = (\psi_\rmh,0)$
of the map~(\ref{eq:f}) collapse as $\mu \to 0^+$.
The following rough quantitative estimates are typical for
saddle-center bifurcations.
The size of $\mathcal{A}_\mu$ in the $\psi$-coordinate
is $\Order(|\psi_\rms-\psi_\rmh|) = \Order(\mu)$.
Its size in the $w$-coordinate is
$\Order(|\psi_\rms-\psi_\rmh|) \times \Order(\mu^{1/2}) = \Order(\mu^{3/2})$,
since the angle between the eigenvectors of the matrix $M_\rmh$
is $\Order(\mu^{1/2})$.
Consequently,
$\mathcal{A}_\mu = \Order(\mu) \times \Order(\mu^{3/2}) = \Order(\mu^{5/2})$.
Finally, $|\mathcal{A}_\mu| \asymp |\mathcal{D}_\mu|$ when $\mu \to 0^+$.
We will confirm and refine these rough estimates.

Following the above comments, 
we scale the $\psi$-coordinate (respectively, $w$-coordinate) by
a factor of order $\mu$ (respectively, order $\mu^{3/2}$).
To be precise, we consider the change of scales
\begin{equation}\label{eq:Scaling_SaddleCenterBifurcation}
x = \psi/\mu,\qquad y = w/\mu^{3/2}.
\end{equation}

If $0 < \mu \ll 1$, then the map~(\ref{eq:f}) is transformed under this change
into a map $(x_1,y_1) = \tilde{f}(x,y)$ of the form
\begin{equation}\label{eq:ftilde_SaddleCenter}
\left\{
\begin{array}{ccl}
x_1 & = & x + \mu^{1/2} y, \\
y_1 & = & y - \mu^{1/2}(x + \pi x^2) + \Order(\mu).
\end{array}
\right.
\end{equation}
The map $\tilde{f}$ is close to the identity map
$\Identity(x,y) = (x,y)$, since
\[
\tilde{f} = \Identity + \mu^{1/2} \tilde{f}_1 + \Order(\mu),\qquad
\tilde{f}_1(x,y) = (y, -x - \pi x^2).
\]
The area-preserving character of $\tilde{f}$ implies that the term
$\mu^{1/2} \tilde{f}_1$ is a Hamiltonian vector field.
Concretely,
\[
\mu^{1/2} \tilde{f}_1 =
\mu^{1/2} \left(\frac{\partial \tilde{H}_1}{\partial y},
               -\frac{\partial \tilde{H}_1}{\partial x}
\right),
\]
where the Hamiltonian $\tilde{H}_1:\Rset^2 \to \Rset$ is given by
\begin{equation}\label{eq:H1_SaddleCenterBifurcation}
\tilde{H}_1(x,y) =
\frac{x^2 + y^2}{2} + \frac{\pi}{3} x^3 - \frac{1}{6\pi^2}.
\end{equation}
We have subtracted the constant $1/6\pi^2$ for convenience.
We denote by $\phi_H^t$ the $t$-time flow of a Hamiltonian $H$.
Then
\begin{equation}\label{eq:f_H1}
\tilde{f} =
\Identity + \mu^{1/2} \tilde{f}_1 + \Order(\mu) =
\phi^{\mu^{1/2}}_{\tilde{H}_1} + \Order(\mu) =
\phi^1_{\mu^{1/2} \tilde{H}_1} + \Order(\mu),
\end{equation}
since the Euler method has a second order
local truncation error.

\begin{figure}
\iffigures
\centering
\psfrag{qh}{$\tilde{q}_\rmh$}
\psfrag{qs}{$\tilde{q}_\rms$}
\psfrag{ar}{$\blacktriangleright$}
\psfrag{al}{$\blacktriangleleft$}
\includegraphics[height=2.4in]{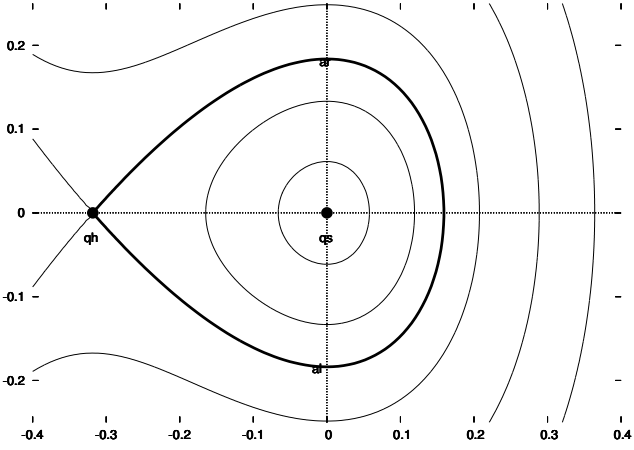}
\else
\vspace{2.4in}
\fi
\caption{The separatrix (thick line), some level curves (thin lines),
and the two equilibrium points (small circles) of the limit
Hamiltonian~(\ref{eq:H1_SaddleCenterBifurcation})
associated to the saddle-center bifurcation.
The black arrows show the Hamiltonian dynamics on the separatrix.}
\label{fig:H1_SaddleCenterBifurcation}
\end{figure}

Therefore, the integrable dynamics of the Hamiltonian $\tilde{H}_1$
approximates the dynamics of the map $\tilde{f}$ in the limit
$\mu \to 0^+$.
We say that $\tilde{H}_1$ is the \emph{limit Hamiltonian} associated
to the saddle-center bifurcation.
It has two equilibrium points.
Namely, the elliptic point $\tilde{q}_\rms = (0,0)$
and the saddle point $\tilde{q}_\rmh = (-1/\pi,0)$.
Besides, the unstable and stable invariant curves of the saddle point
coincide along one branch, giving rise to the separatrix
\[
\tilde{\mathcal{S}} =
\left\{
\tilde{q}=(x,y) \in \Rset^2 : \tilde{H}_1(\tilde{q}) = 0, \ x > -1/\pi
\right\}.
\]
Let $\tilde{\mathcal{R}}$ be the domain enclosed by
$\tilde{\mathcal{S}} \cup \{\tilde{q}_\rmh\}$;
that is,
\[
\tilde{\mathcal{R}} =
\left\{
\tilde{q}=(x,y) \in \Rset^2 : \tilde{H}_1(\tilde{q}) \le 0,\ x \ge -1/\pi
\right\}.
\]
The phase portrait of the limit Hamiltonian~(\ref{eq:H1_SaddleCenterBifurcation})
is sketched in Fig.~\ref{fig:H1_SaddleCenterBifurcation}.
Only the points inside the domain $\tilde{\mathcal{R}}$
give rise to bounded trajectories.
This implies that $\tilde{\mathcal{R}}$ is a good approximation of the scaled versions
of $\mathcal{A}_\mu$ and $\mathcal{D}_\mu$ when $0 < \mu \ll 1$.

The separatrix $\tilde{\mathcal{S}}$ is described by a homoclinic trajectory
$\tilde{q}(t) = (x(t),y(t))$ to the saddle point:
$\lim_{t \to \pm \infty} \tilde{q}(t) = \tilde{q}_\rmh$.
If we impose the initial condition $y(0) = 0$ and
solve the Hamiltonian equations on the separatrix $\tilde{\mathcal{S}}$,
we get that
\begin{equation}\label{eq:HomoclinicOrbit}
x(t) = \frac{3}{2\pi \cosh^2(t/2)}-\frac{1}{\pi},\qquad
y(t) = \frac{3 \sinh (t/2)}{2\pi \cosh^3(t/2)}.
\end{equation}
Therefore,
\[
|\tilde{\mathcal{R}}| = \oint_{\tilde{S}} y \rmd x =
\int_{-\infty}^{+\infty} y(t) x'(t) \rmd t =
\frac{6}{5\pi^2}.
\]
The first identity follows from Green's theorem,
the last one follows from the residue's theorem.
The asymptotic estimate~(\ref{eq:AsymptoticSaddleCenterBifurcation})
follows from the change of scales~(\ref{eq:Scaling_SaddleCenterBifurcation}).

\begin{figure}
\iffigures
\centering
\includegraphics[height=2.4in]{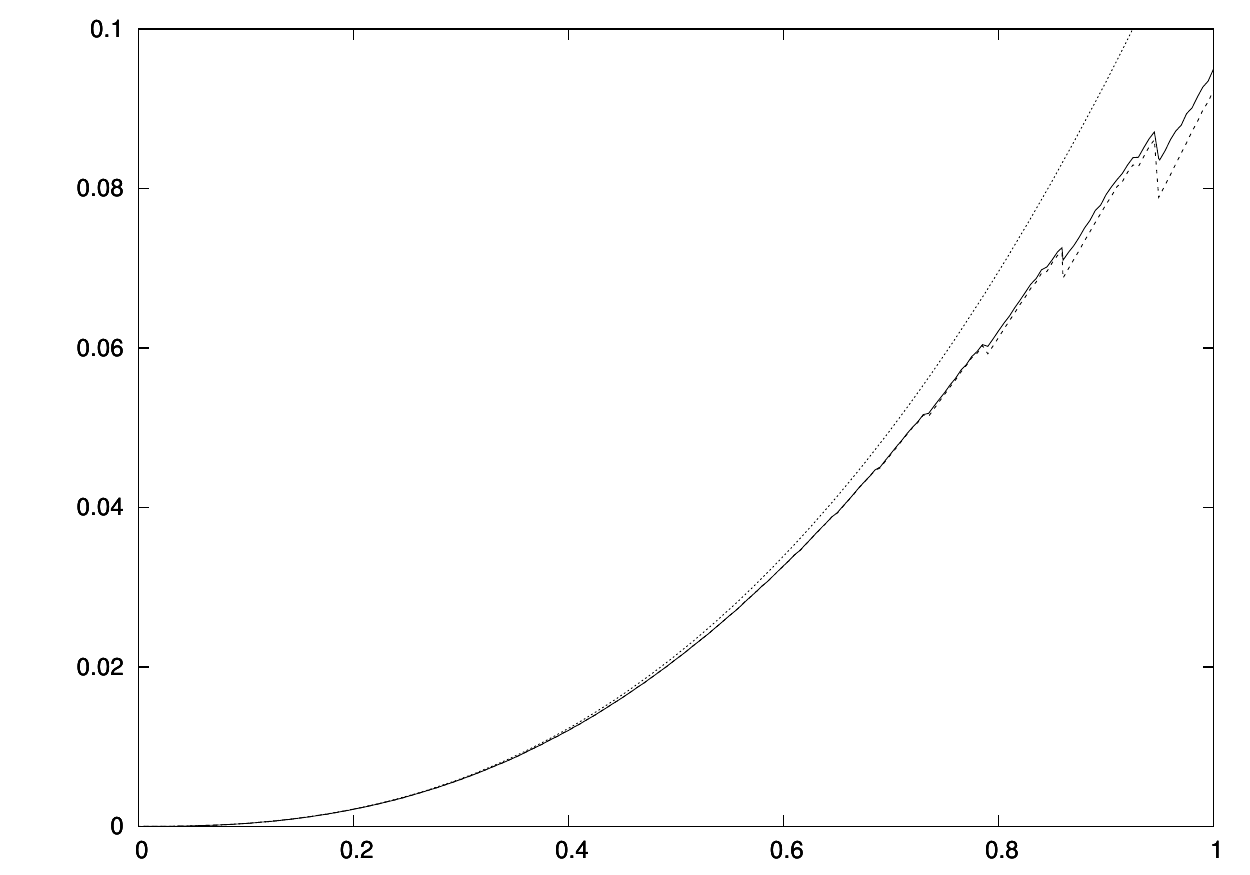}
\else
\vspace{2.4in}
\fi
\caption{$|\mathcal{A}_\mu|$ (continuous line),
$|\mathcal{D}_\mu|$ (dashed line), and the asymptotic
estimate $6\mu^{5/2}/5\pi^2$ (dotted line) versus $\mu$.}
\label{fig:AreaStabilityRegion_close_to_0}
\end{figure}

Fig.~\ref{fig:AreaStabilityRegion_close_to_0} shows a strong agreement
between the numerically computed values of $|\mathcal{A}_\mu|$
and $|\mathcal{D}_\mu|$ and their asymptotic
estimate $\mu \mapsto 6 \mu^{5/2}/5\pi^2$,
even for relatively big values of $\mu$.

Next, we compute better Hamiltonian approximations near
the saddle-center bifurcation.
The Lie's series method is the standard tool to find them,
however we will follow a generating function method that
fits perfectly with the map~(\ref{eq:ftilde_SaddleCenter}).

We recall the following formal result from~\cite[Remark 1]{Simo1982}.

Let $(x_1,y_1) = \tilde{f}(x,y)$ be an analytic area preserving map that is,
in some sense, close to the identity.
Let $\tilde{G}(x,y)$ be a function such that $x_1 y + \tilde{G}(x_1,y)$
is a generating function of $\tilde{f}$; that is,
it satisfies the implicit equations~(\ref{eq:GeneratingFunctionEquations}).
We note that $\| \tilde{G} \| \ll 1$.
The 1-time flow of the Hamiltonian $\tilde{H} = \sum_{j\ge 1} \hat{H}_j$
defined by
\begin{eqnarray*}
\hat{H}_1 & = & \tilde{G}, \\
\hat{H}_2 & = & \textstyle{\frac{1}{2}} \tilde{G}_x \tilde{G}_y, \\
\hat{H}_3 & = & \textstyle{\frac{1}{12}}(\tilde{G}_{xx} \tilde{G}_y^2 +
                4\tilde{G}_{xy} \tilde{G}_x \tilde{G}_y + \tilde{G}_{yy} \tilde{G}_x^2),\\
\hat{H}_4 & = & \textstyle{\frac{1}{12}}(\tilde{G}_{xxy} \tilde{G}_y +
                \tilde{G}_{xyy} \tilde{G}_x + \tilde{G}_{xx} \tilde{G}_{yy} +
                3\tilde{G}_{xy}^2)\tilde{G}_x \tilde{G}_y \\
          &   & \textstyle{\frac{1}{12}} \tilde{G}_{xy}(\tilde{G}_{xx} \tilde{G}_y^2
                + \tilde{G}_{yy} \tilde{G}_x^2),
\end{eqnarray*}
and so on, coincides with the map $\tilde{f}$ at a \emph{formal} level.
There is a small misprint in~\cite{Simo1982};
the numerical factors in the fourth term are $\frac{1}{12}$, not $\frac{1}{2}$.
Besides, we must skip the minus signs in front of the even terms that
appear in~\cite{Simo1982}, because our generating function has the
form $\tilde{G}(x_1,y)$, instead of $\tilde{G}(x,y_1)$.
Here, subindexes in $\tilde{G}$ mean partial derivatives.
We note that $\hat{H}_j = \Order(\| \tilde{G} \|^j)$.
The series $\tilde{H} = \sum_{j\ge 1} \hat{H}_j$ is generically divergent.
Otherwise, $\tilde{f} = \phi^1_{\tilde{H}}$ would be integrable,
which is an exceptional situation.

Now, we come back to our concrete problem.

We write the map~(\ref{eq:ftilde_SaddleCenter}) as
\[
\left\{
\begin{array}{ccl}
x_1 & = & x + \mu^{1/2} y, \\
y_1 & = & y - \mu^{1/2} c(x_1),
\end{array}
\right.
\]
where
\[
c(x_1) =
\frac{2\pi (1 - \cos(\mu x_1)) + \mu \sin(\mu x_1)}{\mu^2} =
x_1 + \pi x^2_1 + \Order(\mu^2).
\]
Then it is easy to check that the function
\[
\tilde{G}(x_1,y) =
\mu^{1/2} \left( y^2/2 + \int_{-1/\pi}^{x_1} c(s) \rmd s \right) =
\mu^{1/2} \tilde{H}_1(x_1,y) + \Order(\mu^{5/2}),
\]
satisfies the implicit equations~(\ref{eq:GeneratingFunctionEquations}),
where $\tilde{H}_1(x,y)$ is the limit
Hamiltonian~(\ref{eq:H1_SaddleCenterBifurcation}).
Let $\sum_{j\ge 1} \hat{H}_j$ be the formal series associated to this
generating function.
We note that $\tilde{G} = \Order(\mu^{1/2})$,
and so $\hat{H}_j = \Order(\mu^{j/2})$ for all $j \ge 1$.
Therefore,
if we retain just the first $n$ terms of the formal series,
we get an approximating Hamiltonian
$\tilde{H}^{[n]}(x,y;\mu) = \sum_{j = 1}^n \hat{H}_j(x,y;\mu)$ such that
\[
\tilde{f} = \phi^1_{\tilde{H}^{[n]}} + \Order(\mu^{(n+1)/2}).
\]
It is interesting to compare this formula with~(\ref{eq:f_H1}).
On the other hand, using the stronger estimate
$\tilde{G} = \mu^{1/2} \tilde{H}_1 + \Order(\mu^{5/2})$,
we get that $\hat{H}_j(x,y;\mu) = \mu^{j/2} \tilde{H}_j(x,y)$
for some function $\tilde{H}_j(x,y)$ that does not depend on $\mu$
for all $j \le 4$.
Indeed,
\begin{eqnarray*}
\tilde{H}_2(x,y) & = &
(\tilde{H}_1)_x (\tilde{H}_1)_y/2 = (x + \pi x^2)y/2,\\
\tilde{H}_3(x,y) & = &
\left(
(\tilde{H}_1)_{xx} (\tilde{H}_1)_y^2 + (\tilde{H}_1)_{yy} (\tilde{H}_1)_x^2
\right)/12 \\
& = & (x + \pi x^2)^2/12 + (1 + 2\pi x)y^2/12, \\
\tilde{H}_4(x,y) & = &
(\tilde{H}_1)_{xx} (\tilde{H}_1)_{yy} (\tilde{H}_1)_x (\tilde{H}_1)_y /12 \\
& = & (1 + 2\pi x)(x+\pi x^2)y/12,
\end{eqnarray*}
since all mixed partial derivatives of $\tilde{H}_1$ vanish.

\begin{figure}
\iffigures
\centering
\includegraphics[height=2.4in]{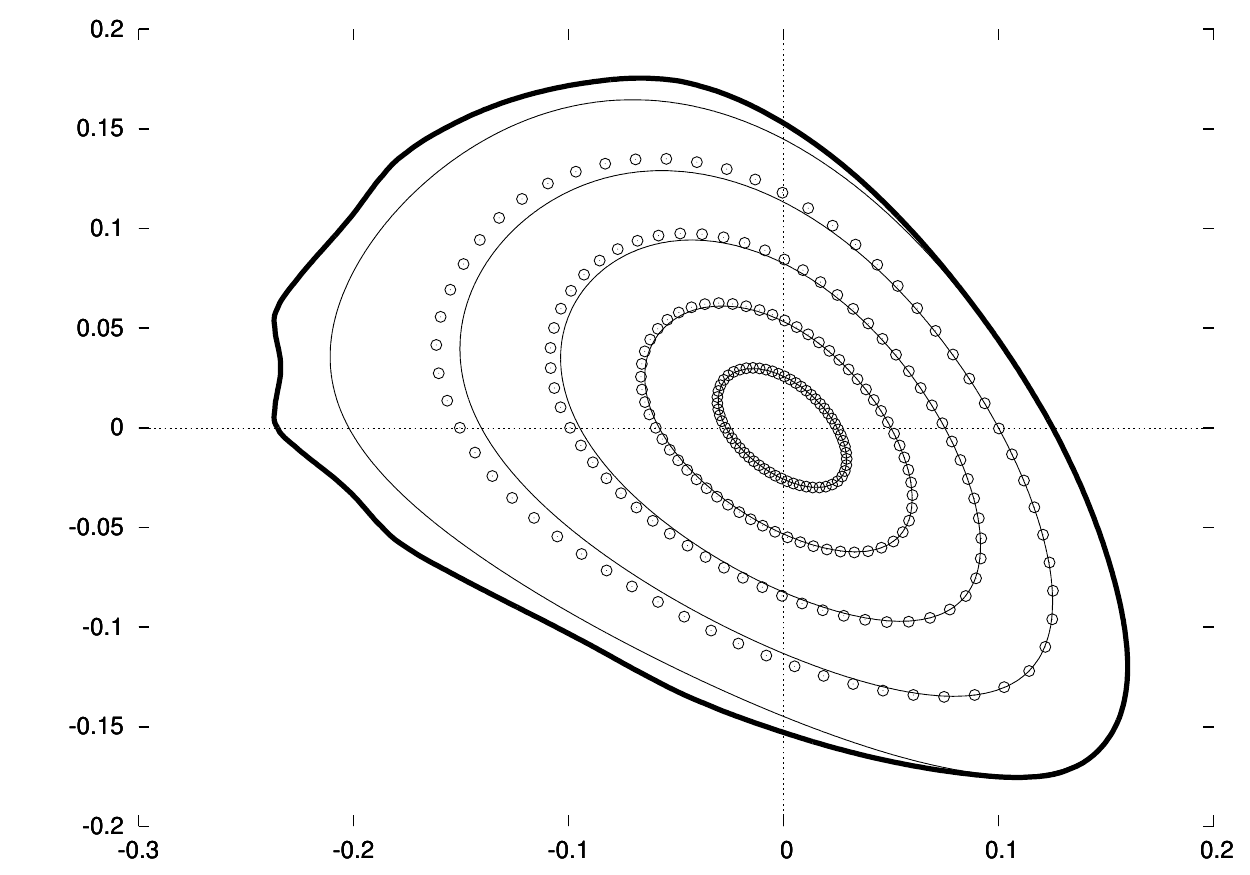}
\else
\vspace{2.4in}
\fi
\caption{Some level curves of the fourth order approximating Hamiltonian
associated to the saddle-center bifurcation (thin lines),
the LRIC (thick line), and some RICs of the microtron map
(small circles) for $\mu = 1$.}
\label{fig:LevelCurvesHamiltonian_mu_1}
\end{figure}

We have plotted some level curves of the approximating Hamiltonian
\[
H^{[4]}(\psi,w;\mu) =
\tilde{H}^{[4]}(\psi/\mu,w/\mu^{3/2};\mu) =
\sum_{j=1}^4 \mu^{j/2} \tilde{H}_j(\psi/\mu,w/\mu^{3/2})
\]
for $\mu = 1$ in Fig.~\ref{fig:LevelCurvesHamiltonian_mu_1}.
We have used the original coordinates $(\psi,w)$ for the sake of comparison,
instead of the scaled ones~(\ref{eq:Scaling_SaddleCenterBifurcation}).
We have also plotted some RICs of the microtron map~(\ref{eq:f}).
We stress that, in spite of the relatively big value of $\mu$,
the RICs fit surprisingly well with the level curves,
being this fitting better close to the elliptic fixed point.
The fitting improves when $\mu$ is smaller.
The LRIC is represented in a thick line.
The ``last'' level curve correspond to the value
$H^{[4]}(\psi,w) \simeq -0.00465$.

Finally, let us discuss briefly how to refine the asymptotic estimate
$|\mathcal{D}_\mu| \asymp 6 \mu^{5/2}/5\pi^2$,
although we do not carry out the associated computations.
They are too cumbersome.

If $0 < \mu \ll 1$, the Hamiltonian $\tilde{H}^{[n]}(x,y;\mu)$ has,
besides the elliptic point $(0,0)$,
a hyperbolic saddle point close to $\tilde{q}_\rmh=(-1/\pi,0)$ whose
unstable and stable invariant curves coincide along one branch.
Let $\tilde{\mathcal{R}}^{[n]}_\mu$ be the domain enclosed by
the corresponding separatrix.
Then $|\tilde{\mathcal{R}}^{[n]}_\mu| =
\sum_{j = 0}^{n-1} \alpha_j \mu^{j/2} + \Order(\mu^{n/2})$ for some
coefficients $\alpha_0,\ldots,\alpha_{n-1}$.
We have already seen that $\alpha_0 = 6/5\pi^2$.
In this way, we may get the refined asymptotic estimates
\begin{eqnarray*}
|\mathcal{D}_\mu| & = & \mu^{5/2} |\tilde{\mathcal{D}}_\mu | \asymp
\mu^{5/2} \left(|\tilde{\mathcal{R}}^{[n]}_\mu| + \Order(\mu^{n/2}) \right) \\
& \asymp &
6\mu^{5/2}/5\pi^2 + \cdots + \alpha_{n-1} \mu^{n/2 + 2} + \Order(\mu^{(n+5)/2}). 
\end{eqnarray*}
Anyway, we stress that $|\mathcal{D}_\mu|$ is not smooth in the parameter $\mu$.
Rather the opposite happens;
$|\mathcal{D}_\mu|$ has a fractal self-similar structure with infinitely
many jumps, although they become very small when $\mu \to 0^+$.
See~\cite{SimoTreschev1998,SimoVieiro2009}, Fig.~\ref{fig:AreaStabilityRegion_close_to_0},
and Fig.~\ref{fig:AreaStabilityRegion}.

\subsection{At the fourth order resonance}
\label{Ssec:HamiltonianFourthOrder}

Let $\tilde{f}$ be the transformed map in the
coordinates~(\ref{eq:CoordinatesFourthOrder}),
which were the suitable ones for $\mu = 2$.

The map $\tilde{f}^4$ is very close to the identity
in a neighborhood of the origin when $\mu = 2$.
To be precise, $\tilde{f}^4 = \Identity + O_3(x,y)$.
Hence, given any order $n \ge 4$,
there exists a unique approximating Hamiltonian of the form
\[
\tilde{H}^{[n]}(x,y) = \tilde{H}_4(x,y) + \cdots + \tilde{H}_n(x,y),
\]
being $\tilde{H}_j(x,y)$ a homogeneous polynomial of degree $j$, such that
\[
\tilde{f}^4 = \phi^1_{\tilde{H}^{[n]}} + O_n(x,y).
\]
A rather tedious computation, which we have the good taste to omit,
leads to the following formulas:
\begin{eqnarray*}
\tilde{H}_4(x,y) & = & -(x^4+y^4)/6 - \pi^2 x^2 y^2, \\
\tilde{H}_5(x,y) & = & -\pi^3 (x^4 y + x y^4) - \pi(x^3 y^2 + x^2 y^3)/3,\\
\tilde{H}_6(x,y) & = & \alpha (x^6 + y^6) + \beta (x^4 y^2 + x^2 y^4) + \gamma x^3 y^3,
\end{eqnarray*}
with $\alpha = 1/180 - \pi^4/3$, $\beta = -5\pi^2/12$,
$\gamma = 1/9 - 2\pi^4$.

All approximating Hamiltonians are symmetric:
\[
\tilde{H}^{[n]} (x,y) = \tilde{H}^{[n]} (y,x), \qquad \forall n \ge 4.
\]
This is a consequence of the fact that the reversor $r_0$ that we will
give in~(\ref{eq:Reversors}) becomes $\tilde{r}_0(x,y) = (y,x)$
in the coordinates~(\ref{eq:CoordinatesFourthOrder}).

\begin{figure}
\iffigures
\centering
\includegraphics[height=2.4in]{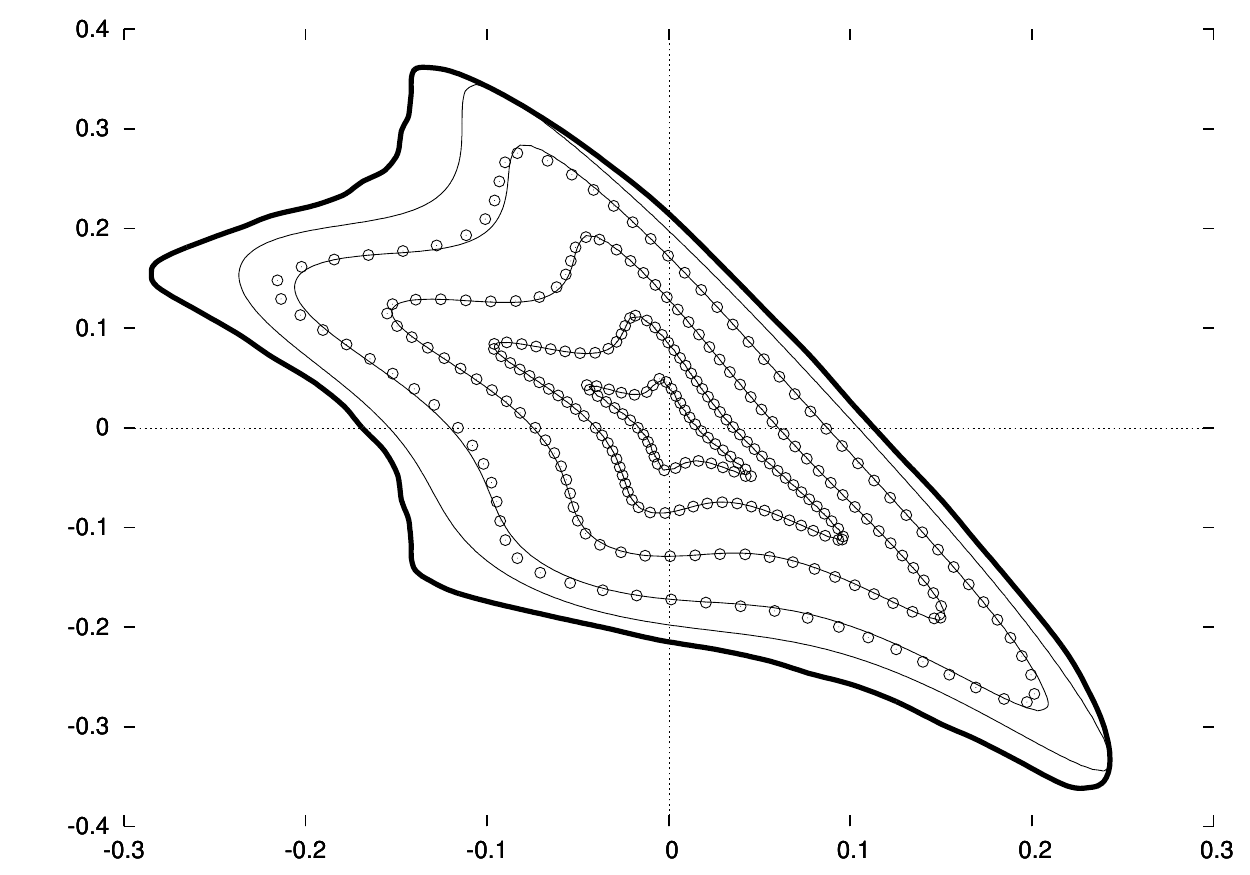}
\else
\vspace{2.4in}
\fi
\caption{Some level curves of the polynomial Hamiltonian
of degree six associated to the fourth order resonance (thin lines),
the LRIC (thick line), and some RICs of the microtron map (small circles)
for $\mu = 2$.}
\label{fig:LevelCurvesHamiltonian_mu_2}
\end{figure}

We have plotted some level curves of the Hamiltonian
\[
H^{[6]}(\psi,w) = \tilde{H}^{[6]}(\psi+w,\psi)
\]
in Fig.~\ref{fig:LevelCurvesHamiltonian_mu_2}.
We have used the original coordinates $(\psi,w)$.
The RICs of the microtron map fit quite well with the level curves.
Of course, the fitting improves close to the origin. 
The LRIC is represented in a thick line.
The ``last'' level curve correspond to the value $H^{[6]}(\psi,w) \simeq -0.0022$.

\subsection{Near the third order resonance}
\label{Ssec:HamiltonianThirdOrder}

\begin{figure}
\iffigures
\centering
\psfrag{qs}{\footnotesize $\tilde{q}_\rms$}
\psfrag{q1}{\footnotesize $\tilde{q}_1$}
\psfrag{q2}{\footnotesize $\tilde{q}_2$}
\psfrag{q3}{\footnotesize $\tilde{q}_3$}
\includegraphics[height=2.4in]{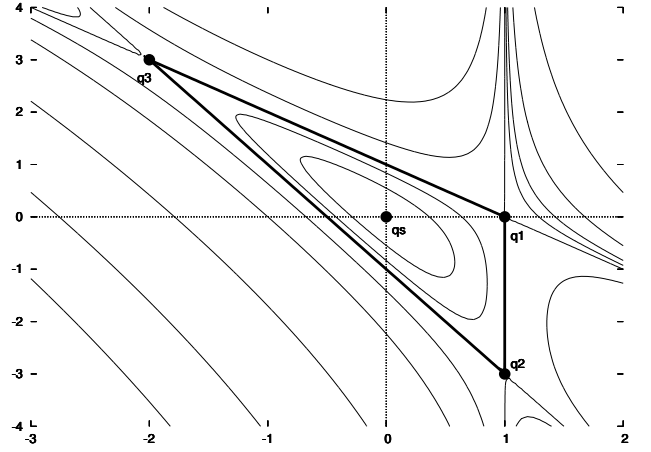}
\else
\vspace{2.4in}
\fi
\caption{The three separatrices (thick lines), some level curves (thin lines),
and the four equilibrium points (small circles) of the the limit
Hamiltonian~(\ref{eq:H1_ThirdOrderResonance})
associated to the third order resonance.}
\label{fig:H1_ThirdOrderResonance}
\end{figure}

Let us study the size and shape of the connected component
$\mathcal{D}_\mu$ of the acceptance when $\mu = 3 + \epsilon$
with $0 < |\epsilon| \ll 1$.
We will see that $\mathcal{D}_{3+\epsilon}$ is approximately
a triangle of vertices~(\ref{eq:TriangleVertices}),
which implies that the asymptotic
estimate~(\ref{eq:AsymptoticThirdOrderResonace}) holds.
We scale the variables according to these claims,
which suggest that all the interesting dynamics takes place in a
$O(\epsilon)$-neighborhood of the origin.
To be precise, we consider the change of scale
\begin{equation}\label{eq:Scaling_ThirdOrder}
x = \pi \psi/\epsilon,\qquad y = \pi w/\epsilon.
\end{equation}
If $\mu = 3+\epsilon$ with $|\epsilon| \ll 1$,
then the map~(\ref{eq:f}) is transformed under this change
into a map $(x_1,y_1) = \tilde{f}(x,y)$ of the form
\[
\left\{
\begin{array}{ccl}
x_1 & = & x + y, \\
y_1 & = & -3x - 2y - x_1(x_1+1)\epsilon + \Order(\epsilon^2).
\end{array}
\right.
\]
The map $\tilde{f}^3$ is close to the identity:
$\tilde{f}^3 = \Identity + \epsilon \tilde{f}_1 + \Order(\epsilon^2)$,
with
\[
\tilde{f}_1(x,y) = ( 3x + 2y - 3x^2 - 2xy, -6x - 3y + 6x^2 + 6xy + y^2 ).
\]
Following Section~\ref{Ssec:HamiltonianSaddleCenterBifurcation},
we realize that the first order term $\epsilon \tilde{f}_1$ is
a Hamiltonian vector field with Hamiltonian $\epsilon \tilde{H}_1$,
where
\begin{eqnarray}
\nonumber
\tilde{H}_1(x,y) & = & 3x^2 + 3xy + y^2 - 2 x^3 - 3 x^2 y - xy^2 - 1 \\
\label{eq:H1_ThirdOrderResonance}
& = & (1-x)(x+y-1)(2x+y+1).
\end{eqnarray}
Therefore, $\tilde{f}^3 = \phi^1_{\epsilon \tilde{H}_1} + \Order(\epsilon^2)$,
and the integrable dynamics of the Hamiltonian $\tilde{H}_1$ approximates
the dynamics of the map $\tilde{f}^3$ when $\mu \simeq 3$.
Thus, $\tilde{H}_1$ is the \emph{limit Hamiltonian} associated
to the third order resonance.
It has four equilibrium points.
One elliptic point: $\tilde{q}_\rms = (0,0)$, and three saddle points:
\[
\tilde{q}_1 = (1,0), \qquad
\tilde{q}_2 = (1,-3), \qquad
\tilde{q}_3 = (-2,3).
\]
Let $\tilde{\mathcal{R}}$ be the triangle whose vertices are
these three points.
Clearly, $|\tilde{\mathcal{R}}| = 9/2$.
Each side of $\tilde{\mathcal{R}}$ is both the stable invariant curve
of a saddle, and the unstable invariant curve of another saddle.
That is, $\partial \tilde{\mathcal{R}}$ is formed by the three saddles
and the three straight separatrices connecting them.
See Fig.~\ref{fig:H1_ThirdOrderResonance}.

\begin{figure}
\iffigures
\centering
\includegraphics[height=2.4in]{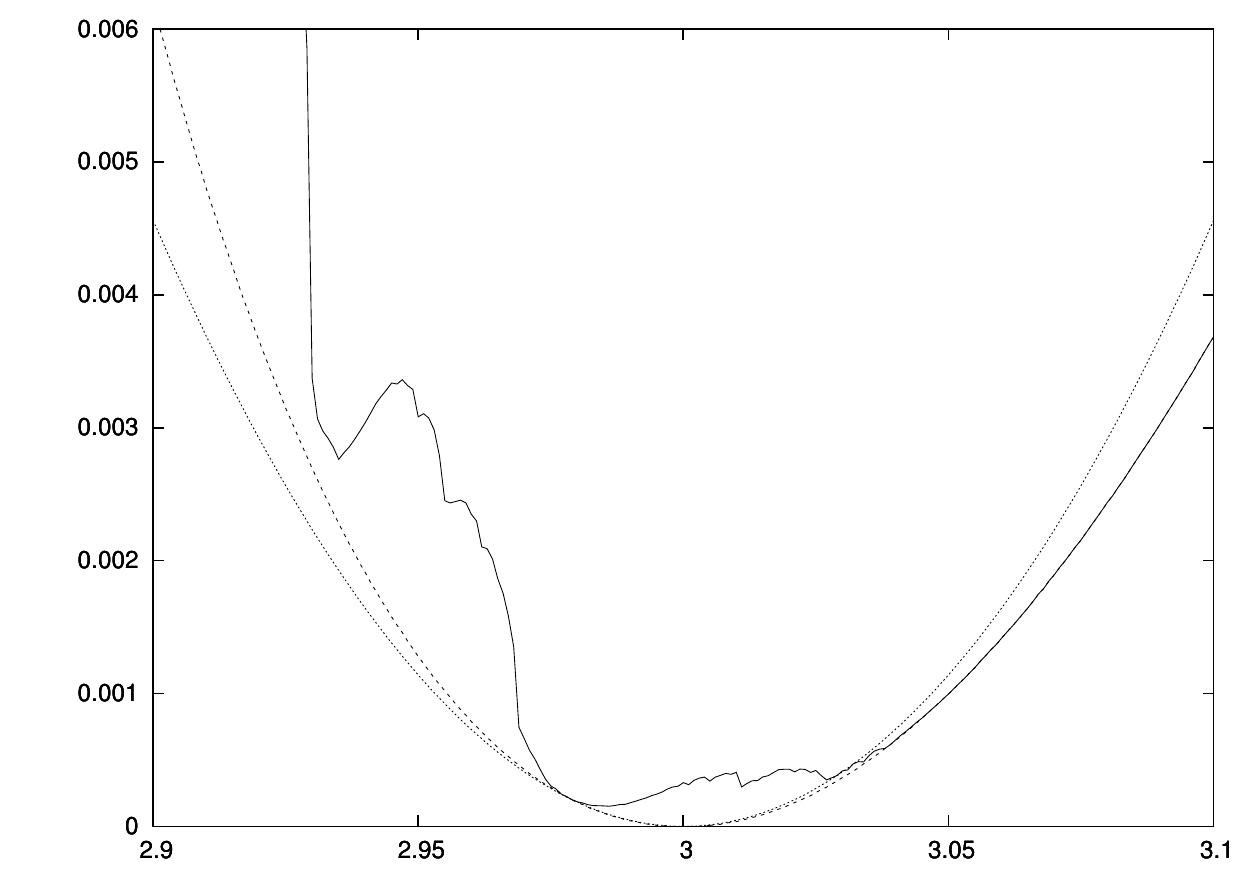}
\else
\vspace{2.4in}
\fi
\caption{$|\mathcal{A}_\mu|$ (continuous line),
$|\mathcal{D}_\mu|$ (dashed line), and the asymptotic
estimate $9(\mu-3)^2/2\pi^2$ (dotted line) versus
the parameter $\mu$.}
\label{fig:AreaStabilityRegion_close_to_3}
\end{figure}

Only the points inside the closed domain $\tilde{\mathcal{R}}$
give rise to bounded trajectories,
so $\tilde{\mathcal{R}}$ is a good approximation of the scaled version
of $\mathcal{D}_\mu$ for $\mu \simeq 3$.
Thus, we get the quadratic asymptotic estimate~(\ref{eq:AsymptoticThirdOrderResonace}).
We have displayed a comparison between the numerically computed area
$|\mathcal{D}_{3+\epsilon}|$ and its asymptotic estimate $9\epsilon^2/2\pi^2$
in Fig.~\ref{fig:AreaStabilityRegion_close_to_3}.
We note that $|\mathcal{A}_3| > 0$,
because of the traces of the $(1,3)$-periodic chain of elliptic islands
that was thrown away from the main connected component $\mathcal{D}_\mu$
at some value $\mu_\star \approx 2.85$.
We will visualize these traces in Fig.~\ref{fig:StabilitySection}.

\section{Global stability of the synchronous trajectory}
\label{Sec:GlobalStability}

This section has a global and experimental character.
We will numerically study the stability domain $\mathcal{A}$
and its connected component $\mathcal{D}$ in the range $0 < \mu < 4.6$.
We will also describe some ideas behind the algorithms.

\subsection{The reversors}
\label{Ssec:Reversors}

A map is \emph{reversible} when each orbit is related to its time
reverse orbit by a symmetry transformation, called a \emph{reversor}.
If a map is reversible, many of their periodic and homoclinic points
are located on certain \emph{symmetry lines}, and many of their
invariant objects are invariant under the reversors~\cite{Devaney1976}.
Two paradigmatic examples of such invariant objects are RICs around
elliptic points and stable and unstable invariant curves of hyperbolic points.
We will use these facts to simplify some computations.

The map~(\ref{eq:f}) can be written as the composition $f = r_1 \circ r_0$,
where $r_0,r_1 : \Tset\times\Tset  \to \Tset\times\Rset$ are the involutions
\begin{equation}\label{eq:Reversors}
r_0(\psi,w) = (\psi + w,-w),\qquad
r_1(\psi,w) = (\psi,\eta(\psi)-w),
\end{equation}
and $\eta(\psi) = 2\pi (\cos \psi - 1) - \mu \sin \psi$.
This means that the map $f$ is reversible with reversors $r_0$ and $r_1$.
See~\cite{LambRoberts1998}.
The symmetry lines of these reversors are their sets of fixed points;
that is,
\begin{eqnarray}\label{eq:FixedSet}
\FixedSet(r_0) & = & \left\{ (\psi,w) \in \Tset \times \Rset : w = 0 \right\}, \\
\nonumber
\FixedSet(r_1) & = &
\left\{ (\psi,w) \in \Tset \times \Rset : w = \eta(\psi)/2 \right\}.
\end{eqnarray}
All the RICs displayed in Fig.~\ref{fig:SummaryOfDynamics} are invariant under
the reversors $r_0$ and $r_1$, and so is the stability domain $\mathcal{A}$.
See Fig.~\ref{fig:StabilityRegionColor_FourthOrderResonance}.

In particular, if we consider the decomposition
\[
\mathcal{A} = \mathcal{A}^- \cup \mathcal{A}^+,\qquad
\mathcal{A}^\pm =
\mathcal{A} \cap \{ (\psi,w)\in \Tset \times \Rset : \pm w \ge 0 \},
\]
then $\mathcal{A}^\pm = r_0(\mathcal{A}^\mp)$.
This halves the computational effort to find $\mathcal{A}$.
The connected component $\mathcal{D}$ satisfies the same property.

The symmetry lines of any reversible map $f = r_1 \circ r_0$ have
many more useful properties.
Let us recall the characterization of \emph{symmetric periodic orbits}
(SPOs) given in~\cite{LambRoberts1998}.
SPOs are the periodic orbits that are invariant under both reversors $r_0$ and $r_1$.
An orbit of $f$ is an SPO if and only if it has exactly two points on
$\FixedSet(r_0) \cup \FixedSet(r_1)$, in which case it has a point
on each symmetry line if and only if it has odd period.
We have displayed in Fig.~\ref{fig:StabilityRegionColor_FourthOrderResonance}
a couple of $(1,4)$-SPOs.
Four is even, so one of these orbits has two points on $\FixedSet(r_0)$,
and the other orbit has two points on $\FixedSet(r_1)$.

\subsection{The rotation number}
\label{Ssec:RotationNumber}

The points inside the stability domain $\mathcal{A}$
rotate around the elliptic fixed point $p_\rms = (0,0)$ when $0 < \mu < 4$.
The points infinitesimally close to $p_\rms$ give $\theta/2\pi$ turns
per iteration, where $\theta$ is the angle defined in~(\ref{eq:RotationNumber}).
We recall that $\theta/2\pi$ is the \emph{rotation number} of the elliptic
point $p_\rms$, and we write $\rho(p_\rms) = \theta/2\pi$.
Next, we try to define the rotation number $\rho(p)$ for any point
$p \in \mathcal{A} \setminus\{ p_\rms \}$
following a standard approach~\cite{KatokHasselblatt1995}.

Given any $p = (\psi,w) \in \mathcal{A}$, let $\varphi$ be its ``argument''.
That is, $\varphi$ is the angle between the segment $[p_\rms,p]$ and the
semi-straight line $\{(\psi,0) : \psi> 0 \}$.
Analogously, let $\varphi_n \in \Rset$ be the ``argument'' of the
$n$-th iterate $p_n = f^n(p)$.
We consider these arguments on the universal cover $\Rset$,
not on $\Tset = \Rset/2\pi\Zset$.
Then we wonder whether the limit
\begin{equation}\label{eq:LimitRotationNumber}
\rho = \rho(p) := \frac{1}{2\pi} \lim_{n \to +\infty} \frac{\varphi_n - \varphi}{n}
\end{equation}
exists.
If so, we say that~(\ref{eq:LimitRotationNumber}) is the \emph{rotation number}
of the point $p$ under the map $f$ around the elliptic point $p_\rms$.
We may be tempted to use the crude numerical approximation
\[
\rho \approx \frac{\varphi_N - \varphi_0}{2\pi N}
\]
for big enough values of $N$,
but it has an $O(1/N)$-error in the most common situations.
Fortunately, it can be refined using the following algorithm.
See~\cite{SearaVillanueva2006,LuqueVillanueva2013} for details.

Given two integers $0 < P < Q$, we set $N = 2^Q$
and compute
\[
S^1_n = \sum_{j=1}^n (\varphi_j - \varphi_0),\quad
S^p_n = \sum_{j=1}^n S^{p-1}_j,\quad
\tilde{S}^p_q = {2^q +p \choose p+1}^{-1} S^p_{2^q},
\]
for $1 \le n \le N$, $1 \le p \le P$, and $1 \le q \le Q$.
Then, under some mild hypotheses, the refined approximation
\[
\rho \approx \Theta(P,Q) := \sum_{p=0}^P
(-1)^{p-1} \frac{2^{p(p+1)/2}}{\delta_p \delta_{P-p}} \tilde{S}^P_{Q-P+p},\quad
\delta_p = \prod_{j=1}^p (2^j -1),
\]
has an $\Order(1/N^{P+1})$-error.
Furthermore, we have the following empirical
bound of the error in the analytic setting:
\begin{equation}\label{eq:ErrorRotationNumber}
| \rho - \Theta(P,Q)| \lessapprox
2^{-(P+1)} |\Theta(P,Q) - \Theta(P,Q-1)|.
\end{equation}

\begin{figure*}
\iffigures
\centering
\includegraphics[height=4.4in]{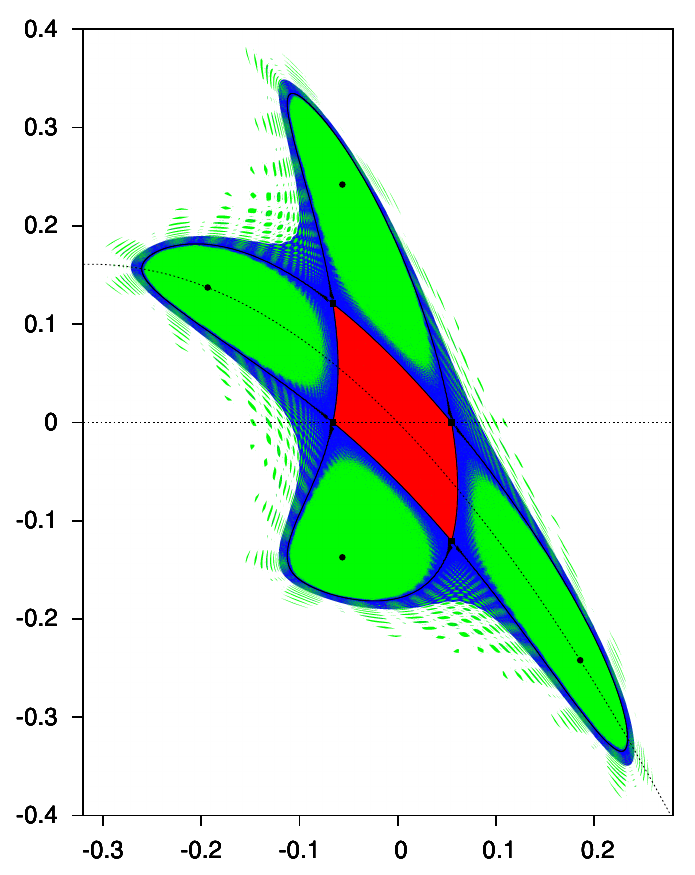}
\includegraphics[height=4.4in]{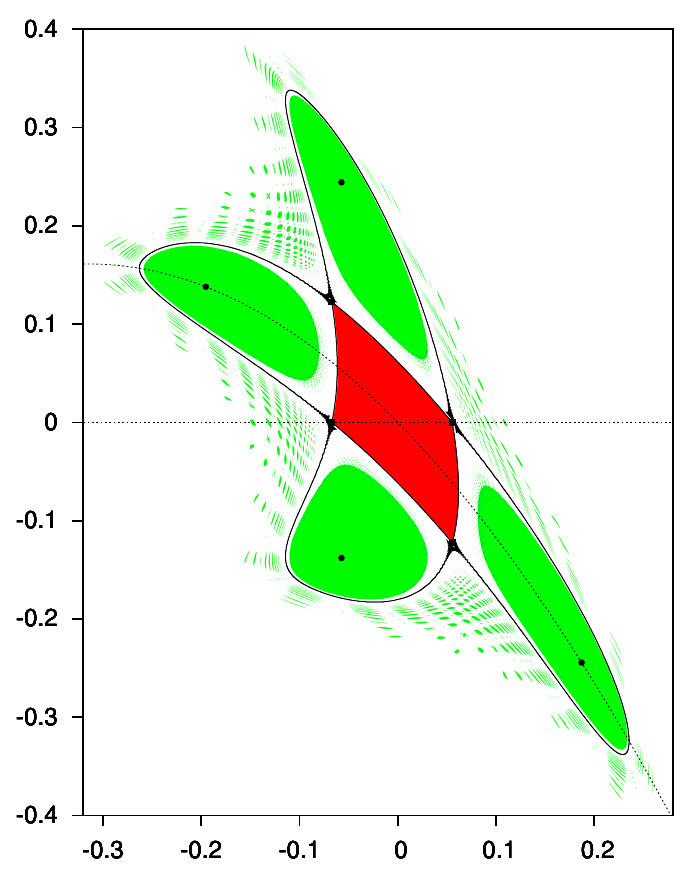}
\else
\vspace{4.5in}
\fi
\caption{The stability domain for $\mu = 2.037$ (left)
and $\mu = 2.038$ (right). Blue corresponds to chaotic seas,
green to periodic elliptic islands, and red to RICs.
The elliptic and hyperbolic $(1,4)$-SPO orbits are marked
with solid black circles and solid black squares, respectively. 
The symmetry lines $\FixedSet(r_0)$ and $\FixedSet(r_1)$
are displayed as dashed black lines.
Each SPO of even period has exactly two points
on a single symmetry line.
A short part of the stable and unstable invariant curves of the
hyperbolic $(1,4)$-SPO is drawn with continuous black lines.
These stable and unstable invariant curves split
(that is, they do not coincide),
but only the outer splitting can be seen at this scale.
See Section~\ref{Sec:InvariantCurves} for details.
In the electronic version one can magnify the plots.}
\label{fig:StabilityRegionColor_FourthOrderResonance}
\end{figure*}

We are interested in rotation numbers because they allow us to
distinguish the three main bounded dynamical behaviors in analytic
area-preserving
diffeomorphisms~\cite{KatokHasselblatt1995,SimoTreschev1998,SearaVillanueva2006}:
\begin{itemize}
\item
If $\mathcal{C}$ is a Moser-like RIC,
then the limit~(\ref{eq:LimitRotationNumber}) exists for all $p \in \mathcal{C}$,
and it does not depend on $p$, so we write $\rho(\mathcal{C})$.
Besides, $\rho(\mathcal{C})$ is, generically, a \emph{Diophantine number}
(that is, it is badly approximated by rational numbers),
and the previous refined algorithm works quite well.
\item
A $(m,n)$-\emph{periodic chain of elliptic islands} is an invariant
region with several connected components such that each
of them surrounds a $(m,n)$-periodic elliptic point.
The $(1,4)$-periodic chain around around $p_\rms$ is displayed
in Fig.~\ref{fig:StabilityRegionColor_FourthOrderResonance}
for two values of $\mu$.
It is formed by the four big green domains.
Each of them is mapped onto the next one in clockwise sense,
so $\rho(p) \equiv 1/4$ for any point $p$ inside them.
Similarly,
if $p$ is inside an $(m,n)$-periodic chain of elliptic islands,
the limit~(\ref{eq:LimitRotationNumber}) exists, $\rho(p) = m/n \in \Qset$,
and the refined algorithm also works quite well.
\item
A \emph{chaotic sea} (or \emph{Birkhoff instability zone})
is the region between two adjacent RICs minus the
stable elliptic islands.
If $p$ is inside a chaotic sea,
then the limit~(\ref{eq:LimitRotationNumber}) generically does not exist,
and the empirical bound~(\ref{eq:ErrorRotationNumber}) does not decrease
when $P$ and $Q$ increase.
\end{itemize}

\begin{figure*}
\iffigures
\centering
\includegraphics[height=3.25in]{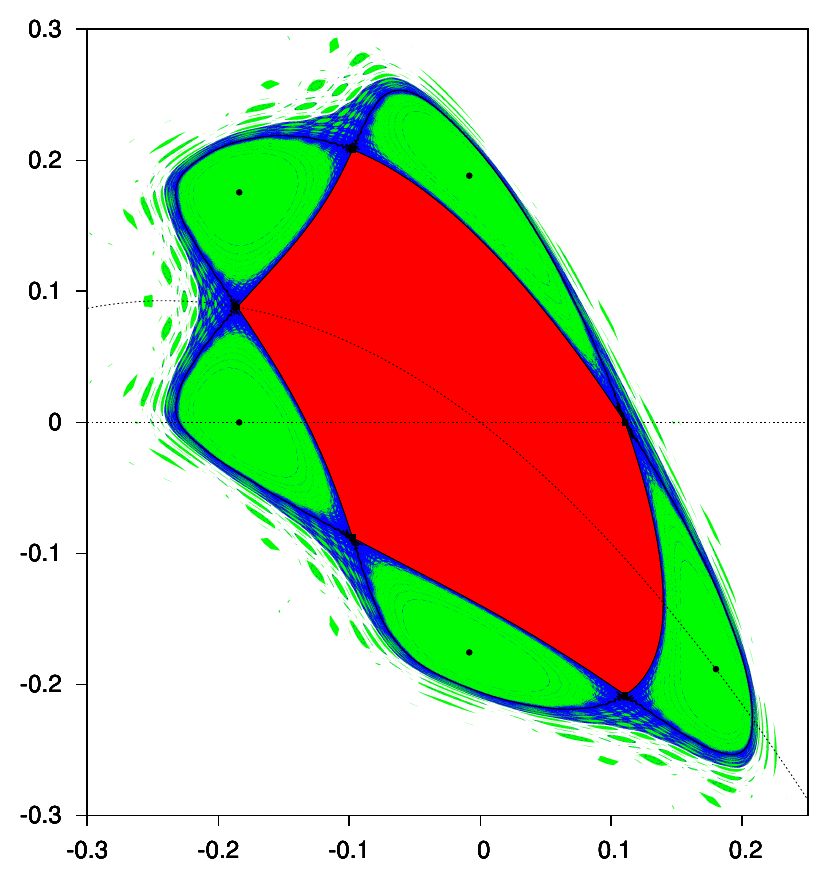}
\hspace{-0.15in}
\includegraphics[height=3.25in]{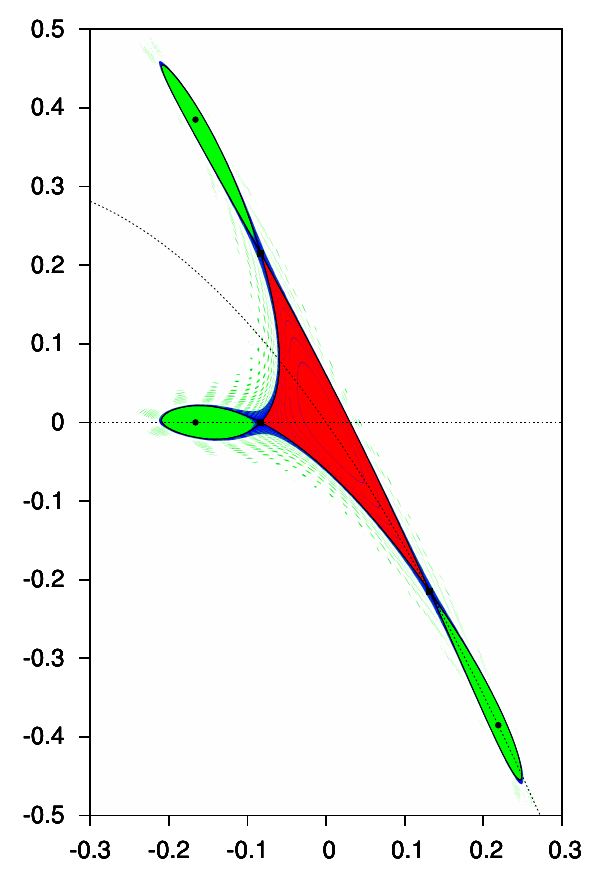}
\hspace{-0.15in}
\includegraphics[height=3.25in]{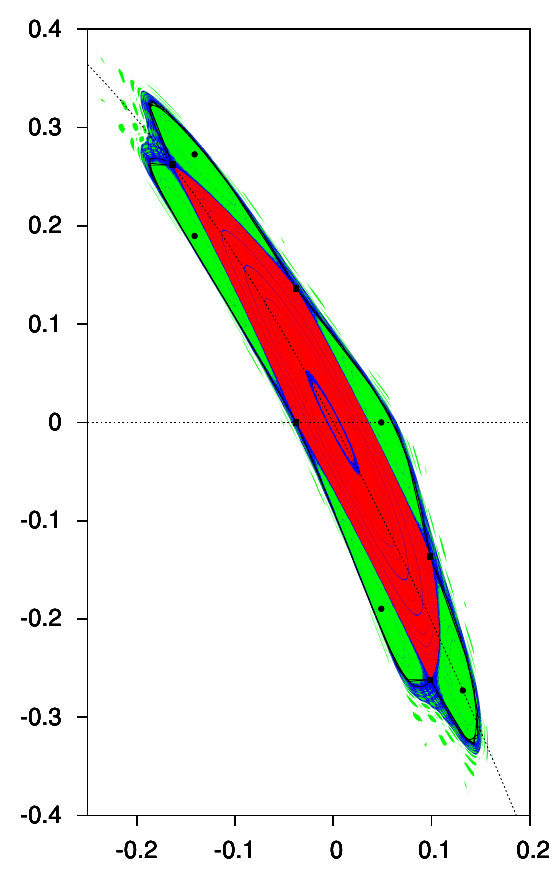}
\else
\vspace{3.25in}
\fi
\caption{The stability domain for $\mu = 1.539$ (left),
$\mu = 2.853$ (center), and $\mu = 3.735$ (right).
Colors, solid black circles, solid black squares, dashed black lines,
and continuous black lines have the same meaning as
in Figure~\ref{fig:StabilityRegionColor_FourthOrderResonance},
but here the SPOs are $(1,5)$-periodic (left),
$(1,3)$-periodic (center), and $(2,5)$-periodic (right).
Each SPO of odd period has exactly one point on each symmetry line.}
\label{fig:StabilityRegionColor_LowOrderResonances}
\end{figure*}

From now on, we paint the points inside chaotic seas in blue,
the ones inside elliptic islands in green,
and the ones on RICs in red.
Fig.~\ref{fig:StabilityRegionColor_FourthOrderResonance} and
Fig.~\ref{fig:StabilityRegionColor_LowOrderResonances} are samples
of that convention.
To be precise, given any point $p \in \mathcal{A}$,
we apply the refined algorithm with $P = 7$ and $Q = 15$ to compute
its rotation number.
If the empirical bound~(\ref{eq:ErrorRotationNumber}) is bigger than
the tolerance $\delta = 10^{-10}$, then we paint $p$ in blue.
Otherwise, we paint $p$ in red/green when $\rho(p)$
is irrational/rational,
which is decided by looking at its continued fraction;
see~\cite{Khinchin1964} for details.

Finally, we can understand the meaning of the small gaps that appear
in the graph of the function $\psi \mapsto \rho(\psi,0)$ displayed in
Fig.~\ref{fig:RotationNumber_mu_root}.
They simply correspond to the sections of some chaotic seas with
the symmetry line $\FixedSet(r_0) = \{ w = 0 \}$.
The rotation number is not well defined on those sections.

\subsection{The stability domain}
\label{Ssec:StabilityDomain}

We compute the stability domain of the microtron map~(\ref{eq:f})
using the \emph{orbit method} described in~\cite{Vieiro2009}.
Visual examples of several stability domains of the H\'enon map provided by
this method can be found in~\cite{SimoVieiro2009,SimoVieiro2011,SimoVieiro2012}.
Our figures show a strong resemblance with those ones.

\begin{figure*}
\iffigures
\centering
\includegraphics[height=2.5in]{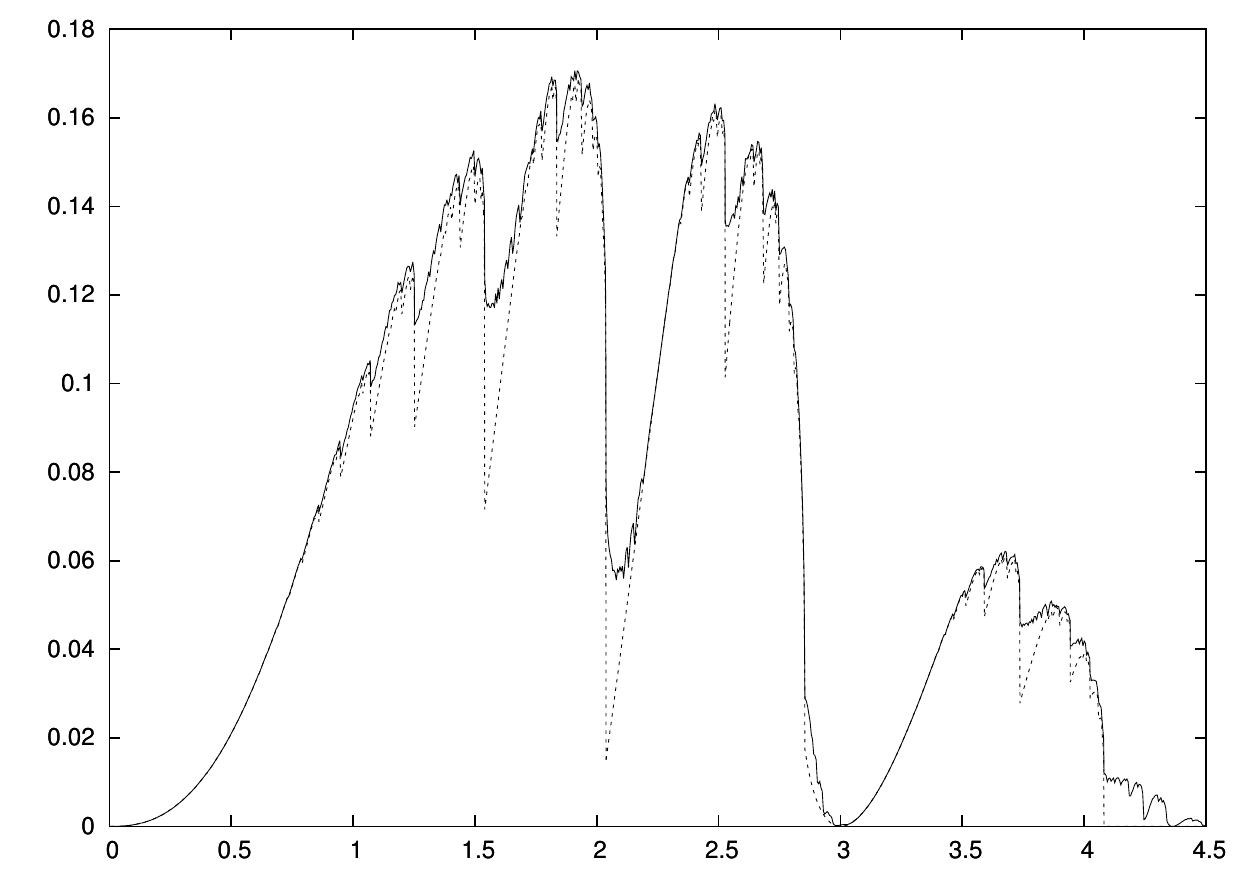}
\includegraphics[height=2.5in]{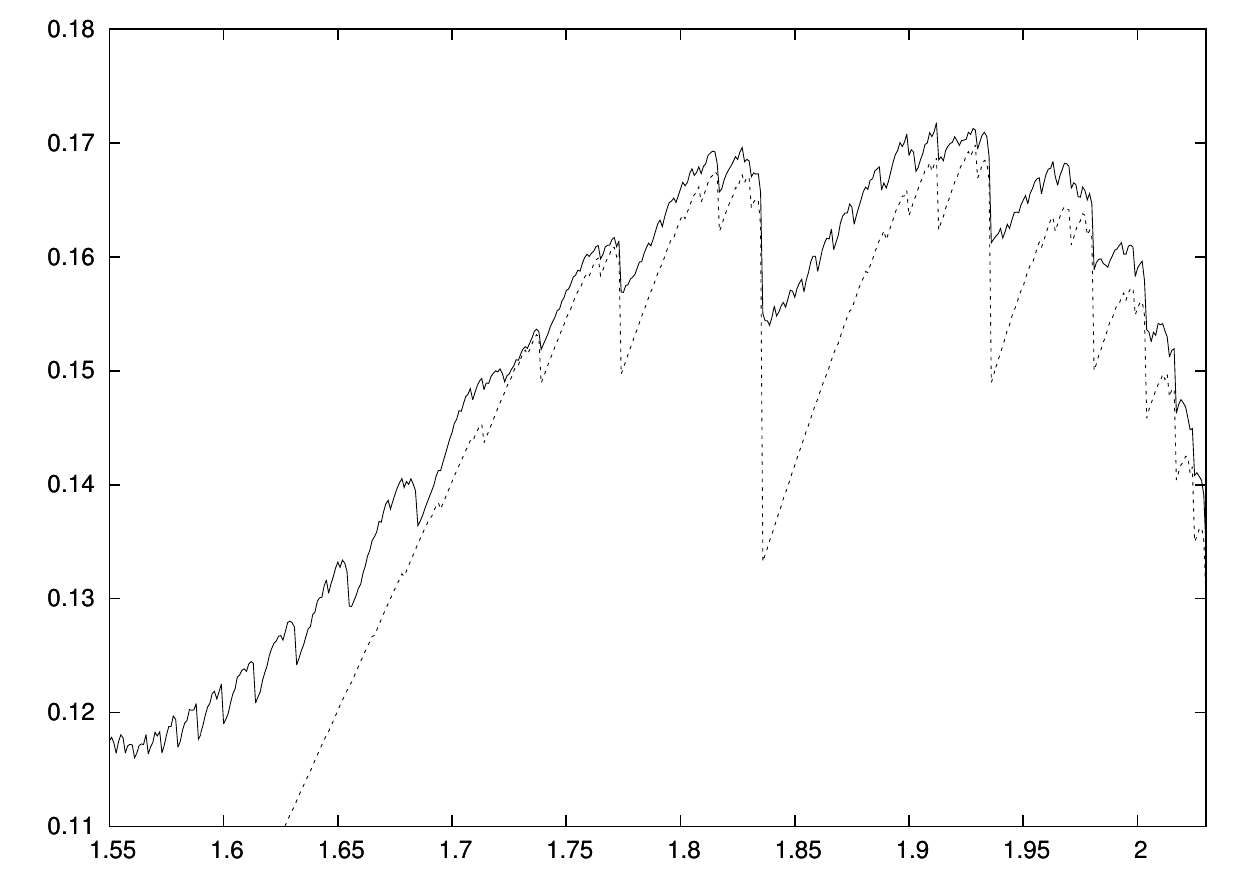}
\else
\vspace{2.5in}
\fi
\caption{Left: $|\mathcal{A}_\mu|$ (continuous line) and
$|\mathcal{D}_\mu|$ (dashed line) versus $\mu$.
Right: A zoom of the previous figure to visualize its self-similar structure.}
\label{fig:AreaStabilityRegion}
\end{figure*}

\begin{figure*}
\iffigures
\includegraphics[height=5.6in]{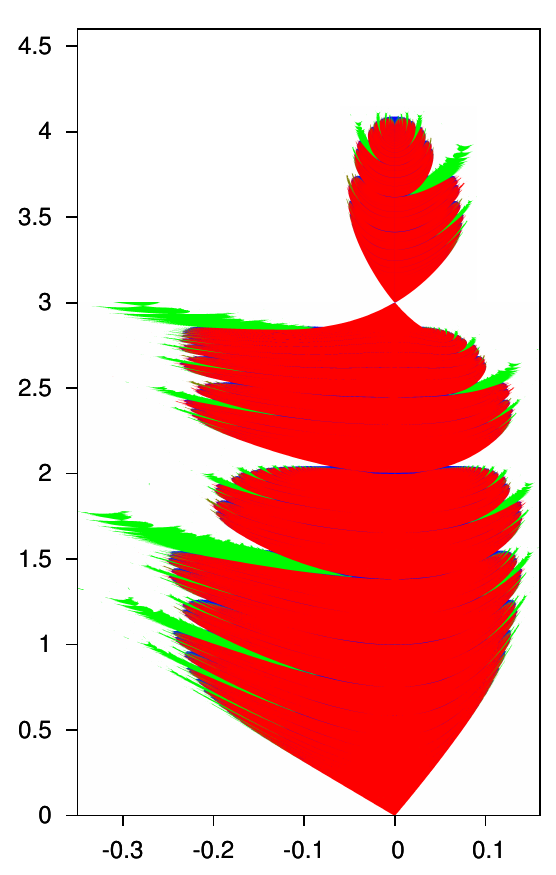}
\includegraphics[height=5.6in]{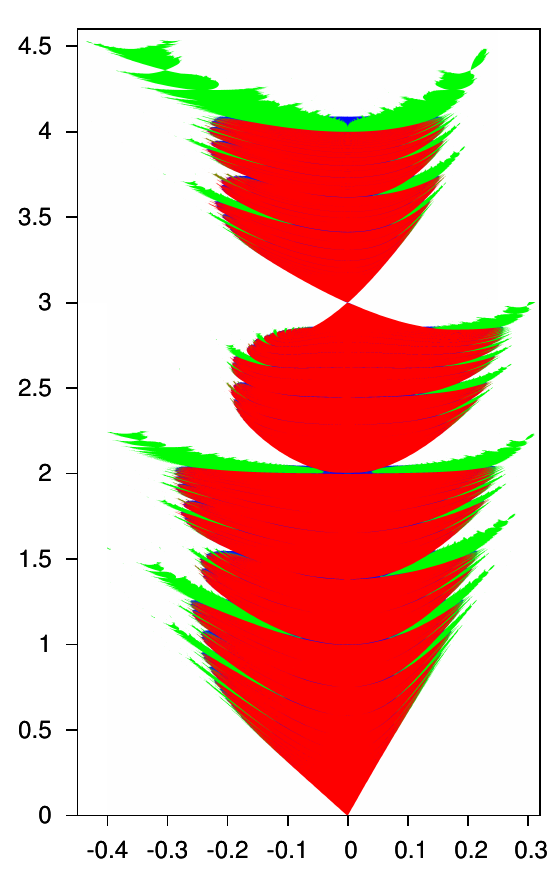}
\else
\vspace{5.7in}
\fi
\caption{The fractal sets $\mathcal{S}_0$ (left) and $\mathcal{S}_1$ (right)
in the $(\psi,\mu)$-plane. Compare with Fig.~2 in~\cite{Henon1969}.}
\label{fig:StabilitySection}
\end{figure*}

Let us describe our implementation of the orbit method.

First, we take a fine rectangular grid in a suitable rectangle
\[
[\psi_{\min},\psi_{\max}] \times [0,w_{\max}]
\]
that contains the upper half $\mathcal{A}^+$ of the stability domain.
Second, we paint in white all grid points such that some of their first
$1000$ iterates escape from the control region
\[
\left\{ (\psi,w) \in \Tset \times \Rset: |w| \le 1 \right\}.
\]
This is a fast step, since $1000$ is a relatively small number
for any modern computer.
Third, we consider the $n$-th iterates with $-10^7 \le n \le 10^7$
of all not-yet-white grid points adjacent to some already-white one.
If some of these iterates escapes from the control region,
we paint in white all grid points ``visited'' by
its corresponding unbounded orbit or by its
$r_0$-symmetric orbit.
We do not use the reversor $r_1$,
since $r_0$ is computationally cheaper and $r_1 = f \circ r_0$.
We repeat this process until no more points escape from the control region.
This is the hardest step.
Fourth, we determine the color of each non-white grid point by
computing its rotation number as explained before.
This step of the algorithm is not ``orbitally coherent''
in the sense of~\cite{Vieiro2009}, but we feel that it is the right choice,
since there is no clear way to choose the color in disputed cases.
Finally, we get the lower half of the stability domain from the
identity $\mathcal{A}^- = r_0 (\mathcal{A}^+)$.

We stress that all grid points in the interior (respectively, ``on'' the border)
of the stability domain have been iterated only 1000
(respectively, $2 \cdot 10^7$) times.
This is good, since the interior contains much more
grid points than the border.

Fig.~\ref{fig:StabilityRegionColor_FourthOrderResonance} shows a couple
of stability domains computed with this algorithm.
Each pixel is the center of a square with side $\ell = 1/4000$,
so the area of the stability domain is approximately equal to $\ell^2$
times the number of colored pixels.
We also count the number of pixels in the connected component
containing the origin following a standard algorithm in Computer Vision.
See, for instance,~\cite[pags.~72--75]{ShapiroStockman2002}.
In that way, we obtain that
$|\mathcal{A}_\mu| \approx 1.1657\cdot 10^{-1}$ and
$|\mathcal{D}_\mu| \approx 1.1029 \cdot 10^{-1}$ for $\mu = 2.037$,
but $|\mathcal{A}_\mu| \approx 7.6105 \cdot 10^{-2}$ and
$|\mathcal{D}_\mu| \approx 1.5067 \cdot 10^{-2}$ for $\mu = 2.038$.

These jumps in $|\mathcal{A}_\mu|$ and $|\mathcal{D}_\mu|$ have
a simple explanation.
We recall that the fourth order resonance takes place at $\mu = \mu_\bullet = 2$.
After that value is crossed,
four $(1,4)$-periodic elliptic islands surrounded by a chaotic sea
emanate from the elliptic point $p_\rms$.
This structure (elliptic islands plus chaotic sea) moves
away from $p_\rms$ as $\mu$ grows,
but remains inside $\mathcal{D}_\mu \subset \mathcal{A}_\mu$
while some RIC surrounds it.
However, the LRIC surrounding it disappears
at some value $\mu = \mu_\star \in (2.037,2.038)$.
After that value is crossed,
both the elliptic islands and the chaotic sea are thrown away
from the connected component $\mathcal{D}_\mu$,
although any elliptic island is, by definition, part
of $\mathcal{A}_\mu$.
Thus,
the jump in $|\mathcal{A}_\mu|$ only takes into account the loss of
the chaotic sea,
whereas the jump in $|\mathcal{D}_\mu|$ also takes into account the loss
of many elliptic islands.

Similar jumps take place for any periodic elliptic island,
although the greater is the order of the expelled island,
the smaller is the jump in both areas.
For instance,
we see in Fig.~\ref{fig:StabilityRegionColor_LowOrderResonances} the
stability domains for $\mu = 1.539$, $\mu = 2.853$, and $\mu = 3.735$,
which are parameter values smaller than, but very close to, the values
at which the periodic elliptic islands of orders three and five
are thrown away from the connected component $\mathcal{D}$.
Hence, the stability domains corresponding to $\mu = 1.540$, $\mu = 2.854$,
and $\mu = 3.736$ are significantly smaller,
because the more extern chaotic sea is lost. 

Fig.~\ref{fig:AreaStabilityRegion} shows $|\mathcal{A}_\mu|$
and $|\mathcal{D}_\mu|$ as a function of $\mu$.
The maximal value $|\mathcal{A}_\mu| \approx 1.718 \cdot 10^{-1}$ is
attained for $\mu \approx 1.912$,
$|\mathcal{D}_\mu| = 0$ for all $\mu \gtrsim 4.08$,
and $|\mathcal{A}_\mu| = 0$ for all $\mu \gtrsim 4.53$.
Obviously, $|\mathcal{D}_\mu| \le |\mathcal{A}_\mu|$,
since $\mathcal{D}_\mu \subset \mathcal{A}_\mu$.
We note that $|\mathcal{A}_\mu|$ displays many more jumps than $|\mathcal{D}_\mu|$,
because it is affected by \emph{secondary resonances}
(resonances inside the islands).
The graph of $|\mathcal{A}_\mu|$ has a self-similar fractal structure
caused by those secondary resonances, as the displayed magnification shows.

We also list in Table~\ref{tab:Summary} the exact value
\begin{equation}\label{eq:mu0}
\mu_\bullet := 2 - 2 \cos(2\pi m/n) \in [0,4]
\end{equation}
at which the elliptic fixed point $p_\rms$ becomes $(m,n)$-resonant,
so that the $(m,n)$-periodic chain of elliptic islands is created from $p_\rms$,
jointly with the numerically approximated value $\mu = \mu_\star$ at which
the $(m,n)$-periodic chain of elliptic islands is thrown away from
$\mathcal{D}$, for all $(m,n)$-resonances of order $n < 10$.
For instance, the $(2,7)$-periodic chain is thrown away at some value
$\mu_\star \in (2.526,2.527)$, which explains the seven ``holes'' delimited
by the unbounded orbit displayed in Fig.~\ref{fig:RootOfTheTwistCoefficient}
for $\mu_\rmr \approx 2.538$.

We observe that $\mu_\bullet < \mu_\star$, but in the $(1,3)$-resonance,
which comes as no surprise, since the twist coefficient is positive
if and only if $2.538 \approx \mu_\rmr < \mu < 3$,
see Subsection~\ref{Ssec:LocalStabilityElliptic}.

\begin{table}
\centering
\begin{tabular}{lll}
\hline
$(m,n)$ & Emanate at & Escape at $\mu_\star$ with \\
\hline
$(1,9)$ & $\mu_\bullet \approx 0.468$                            & $0.859 < \mu_\star < 0.860$ \\
$(1,8)$ & $\mu_\bullet = 2 - \sqrt{2} \simeq 0.586$              & $0.948 < \mu_\star < 0.949$ \\
$(1,7)$ & $\mu_\bullet \approx 0.753$                            & $1.071 < \mu_\star < 1.072$ \\
$(1,6)$ & $\mu_\bullet = 1$                                      & $1.251 < \mu_\star < 1.252$ \\
$(1,5)$ & $\mu_\bullet = \frac{1}{2}(5 - \sqrt{5}) \simeq 1.382$ & $1.539 < \mu_\star < 1.540$ \\
$(2,9)$ & $\mu_\bullet \approx 1.653$                            & $1.835 < \mu_\star < 1.836$ \\
$(1,4)$ & $\mu_\bullet = 2$                                      & $2.037 < \mu_\star < 2.038$ \\
$(2,7)$ & $\mu_\bullet \approx 2.445$                            & $2.526 < \mu_\star < 2.527$ \\
$(1,3)$ & $\mu_\bullet = 3$                                      & $2.853 < \mu_\star < 2.854$ \\
$(3,8)$ & $\mu_\bullet = 2 + \sqrt{2} \simeq 3.414$              & $3.589 < \mu_\star < 3.590$ \\
$(2,5)$ & $\mu_\bullet = \frac{1}{2}(5 + \sqrt{5}) \simeq 3.618$ & $3.735 < \mu_\star < 3.736$ \\
$(3,7)$ & $\mu_\bullet \approx 3.802$                            & $3.942 < \mu_\star < 3.943$ \\
$(4,9)$ & $\mu_\bullet \approx 3.879$                            & $4.023 < \mu_\star < 4.024$ \\
$(1,2)$ & $\mu_\bullet = 4$                                      & $4.080 < \mu_\star < 4.081$ \\
\hline
\end{tabular}
\caption{\label{tab:Summary}Exact values $\mu_\bullet$ at which
the main resonances emanate from $p_\rms$,
and approximated values $\mu_\star$ at which they escape from $\mathcal{D}$.}
\end{table}

We end with a couple of warnings regarding
the accuracy of our pictures.

On the one hand, we know that the red domains displayed
in Fig.~\ref{fig:StabilityRegionColor_FourthOrderResonance},
and in the left picture of Fig.~\ref{fig:StabilityRegionColor_LowOrderResonances}
are not completely filled with RICs.
Indeed, KAM theory implies that the set of all RICs has a
complicated Cantorian structure, whose gaps are filled with resonances.
Let us explain why we do not see those resonances.
We restrict our explanation to the case of
Fig.~\ref{fig:StabilityRegionColor_FourthOrderResonance},
since the phenomenon is the same in both cases.
We do not see any resonance inside the red zone displayed in
Fig.~\ref{fig:StabilityRegionColor_FourthOrderResonance} because
\[
\rho_\star := \rho(p_\rms) =
\theta/2\pi = \acos(1-\mu_\star/2)/2\pi \approx 0.25302,
\]
when $\mu = \mu_\star \in (2.037,2.038)$, see Table~\ref{tab:Relations}.
Therefore, all missed resonances have orders $n \ge 83$,
since $21/83$ is the rational number in the interval $(1/4,\rho_\star)$
with the smallest denominator.
Such resonances are too small to be detected with our pixel resolution.
We recall that a generic $(m,n)$-resonance in an
$\Order(\eta)$-neighborhood of an elliptic fixed point of an analytic
area preserving map has an $\Order(\eta^{n/4})$-size~\cite{SimoVieiro2009}.

On the other hand, the border of the connected component $\mathcal{D}_\mu$
should be a red curve (the LRIC), which is missing in the left picture
of Fig.~\ref{fig:StabilityRegionColor_FourthOrderResonance} and in the three
pictures of Fig.~\ref{fig:StabilityRegionColor_LowOrderResonances}.
The reason is, once more, that the LRIC is too thin
to be detected with our pixel resolution.
In the same way, probably there are some RICs in the middle of the blue zone,
which would mean that that blue zone is composed by several chaotic seas,
instead of by a single big chaotic sea.

\subsection{The sections with the symmetry lines}
\label{Ssec:Sections}

We consider the sections of the stability domain with the symmetry
lines~(\ref{eq:FixedSet}) for two reasons.
First, the stability domain is symmetric with respect to these lines.
Second,
we recall that each SPO has exactly two points on these lines,
and such SPOs are the basic invariant objects that organize
the dynamics inside their resonance.

We look at Fig.~\ref{fig:StabilityRegionColor_FourthOrderResonance}
to understand how SPOs organize the resonant dynamics.
We see an elliptic $(1,4)$-SPO with two points on $\FixedSet(r_1)$ that
organizes the stable dynamics inside the four big green elliptic islands,
and a hyperbolic $(1,4)$-SPO with two points on $\FixedSet(r_0)$
whose stable and unstable invariant curves delimit almost perfectly
the red region ``filled'' with RICs.

We do not deal with each one-dimensional section as a separate object,
but we gather them into the two-dimensional sets:
\begin{eqnarray*}
\mathcal{S}_0 & = &
\left\{ (\psi,\mu) \in \Tset \times (0,+\infty) :
        (\psi,0) \in \mathcal{A}_\mu \right\}, \\
\mathcal{S}_1 & = &
\left\{ (\psi,\mu) \in \Tset \times (0,+\infty) :
        (\psi,\eta(\psi)/2) \in \mathcal{A}_\mu \right\},
\end{eqnarray*}
in order to visualize their evolution in the parameter $\mu$.

These sets are represented in Fig.~\ref{fig:StabilitySection}
with the usual color codes.
We see that the connected component $\mathcal{D}_\mu$
collapses to the elliptic fixed point at $\mu = 3$,
and undergoes its major loss at some $\mu = \mu_\star \in(2.037,2.038)$.
The collapse is associated to the local instability of the third
order resonance.
The loss takes place when the fourth order resonance is thrown away.
These are the most relevant phenomena regarding the size of
the stability domain in generic families of area-preserving maps.

We appreciate a clear self-similar fractal structure in the
sets $\mathcal{S}_0$ and $\mathcal{S}_1$.
One can magnify their pictures to appreciate several details.
Let us describe the main ones.

Resonances emanate from the elliptic fixed point $p_\rms$ in the form
of thin blue and green tongues, since they are sections of chaotic seas
surrounding elliptic islands.
These tongues begin at the points $(\psi,\mu) = (0,\mu_\bullet)$,
where $\mu_\bullet$ is defined as in~(\ref{eq:mu0}) for any rational
number $m/n \in (0,1/2)$, so there are infinitely many of such tongues.
Besides, they look symmetric with respect to the vertical line
$\{\psi = 0\}$ in a neighborhood of that line. 
The resonances, but the $(1,4)$-one, have a small width
when they are close to $p_\rms$.
Hence, most of the tongues become visible only at some distance
of $\{ \psi = 0 \}$.

The elliptic islands inside a resonance grow in size when they move away
from $p_\rms$.
If the resonance has even order,
then the green part becomes the biggest part of the tongue in the set $\mathcal{S}_1$,
but not in the set $\mathcal{S}_0$.
This means that each elliptic SPO with even period has two points on
$\FixedSet(r_1)$, but none on $\FixedSet(r_0)$.
Fig.~\ref{fig:StabilityRegionColor_FourthOrderResonance} is a sample
of that empirical claim.
On the contrary, if the resonance has odd order,
then the green part becomes the biggest part of just one branch of the tongue,
in both sets $\mathcal{S}_0$ and $\mathcal{S}_1$.
The other branch remains blue.
This is a completely expected behavior, since we know that all SPOs
with odd period have just one point on each symmetry line.

We warn that red and green regions should not be in direct contact,
but we have again difficulties to detect the blue chaotic seas
between them with our current pixel resolution.

Resonances are thrown away from $\mathcal{D}$,
which is the reason for the saw-like border of
$\mathcal{S}_0$ and $\mathcal{S}_1$.
Once a green tongue is separated from the main red body in
the $\psi$-direction, it is no longer delimited by blue portions,
because chaotic seas do not form part of $\mathcal{A}$
when their resonances are thrown away.
See Fig.~\ref{fig:StabilityRegionColor_FourthOrderResonance}.
Besides, the separated part of any green tongue shows a shape similar
to the shape of the whole set, which is due to the secondary resonances.
Indeed, we see the two above-mentioned phenomena
(major loss in the stability domain and collapse of the
stability domain to a point) along many of these tongues.

Such phenomena are clearly visible in the tongues of $\mathcal{S}_1$
that cover the range $4.08 < \mu < 4.53$.
The fixed point $p_\rms$ is already globally unstable in that range,
but there is still a locally stable elliptic two-periodic orbit on
the symmetry line $\FixedSet(r_1)$.

H\'enon studied some sets similar to $\mathcal{S}_0$ and
$\mathcal{S}_1$ for the H\'enon map more than
forty years ago~\cite{Henon1969}.
His computations already show many of the above-described phenomena,
in spite of the limitations of the computers in that time.
Such limitations were recently overcome
in~\cite{Migueletal2013}.

\section{On the invariant curves of some hyperbolic points}
\label{Sec:InvariantCurves}

The stable and unstable invariant curves of hyperbolic
(fixed or periodic) points organize the dynamics of area
preserving maps in several ways.
Let us just mention two of them.

On the one hand,
some of these invariant curves are approximate boundaries of
stability domains~\cite{Giovannazzi1996},
in which case the area of the lobes between them is equal to the flux
through certain closed curves composed by arcs of invariant
curves~\cite{Meiss1992}.
These lobes have an exponentially small area in the analytic
case~\cite{FontichSimo1990},
so that the above-mentioned closed curves become \emph{partial barriers}
of the dynamics~\cite{MacKayMP1984,Meiss1992}.
If the map is entire,
then the stable and unstable invariant curves never coincide~\cite{Ushiki1980},
so these partial barriers are never complete barriers.

\begin{figure}
\iffigures
\centering
\includegraphics[height=2.4in]{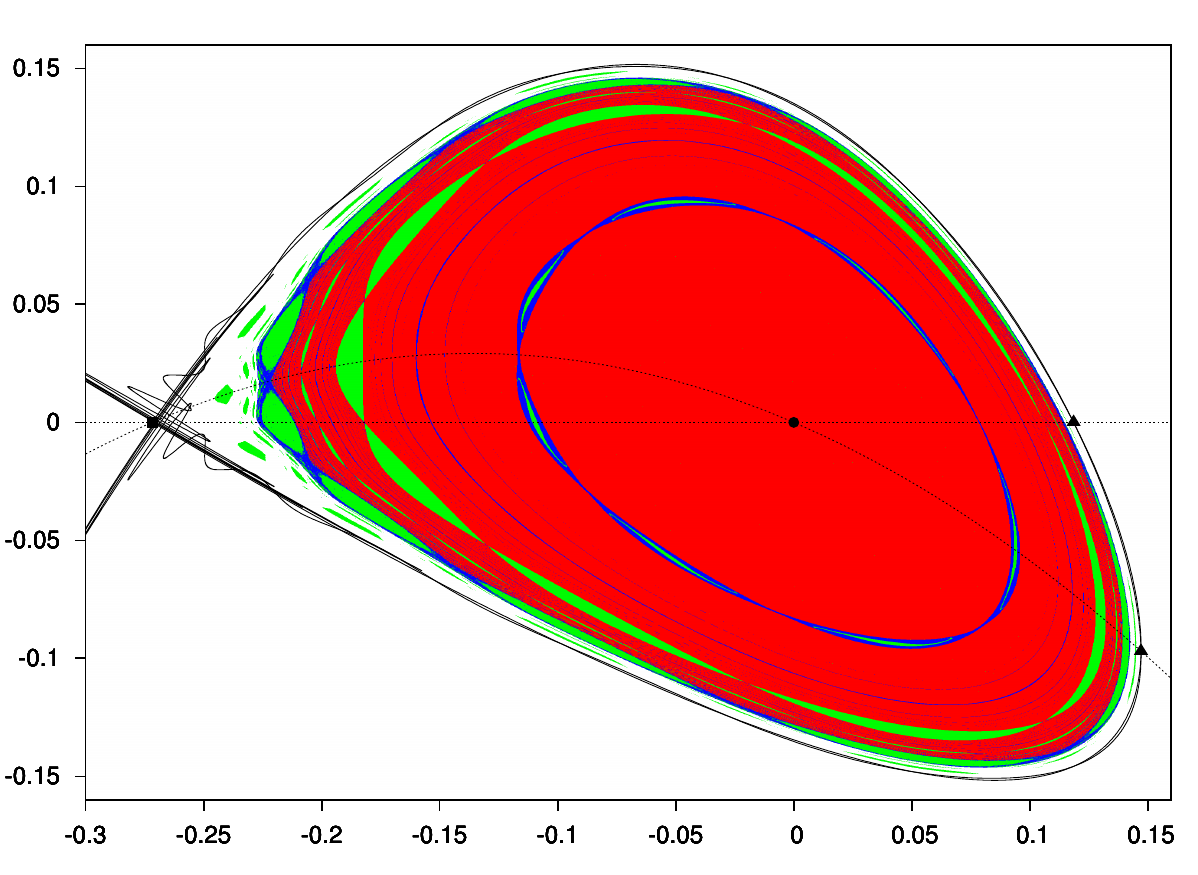}
\else
\vspace{2.4in}
\fi
\caption{The stability domain for $\mu = 0.859$.
Blue corresponds to chaotic seas,
green to periodic elliptic islands, and red to RICs.
The elliptic fixed point $p_\rms = (0,0)$ and
the hyperbolic fixed point $p_\rmh = (-2\phi_\rms,0)$ are marked
with a solid black circle and a solid black square, respectively. 
The symmetry lines $\FixedSet(r_0)$ and $\FixedSet(r_1)$
are displayed as dashed black lines.
A short part of the stable and unstable invariant curves
of $p_\rmh$ is shown as continuous black lines.
The primary intersections of these invariant curves with the
symmetry lines are marked with two solid black triangles.}
\label{fig:StabilityRegion_mu_0p859}
\end{figure}

On the other hand,
if the unstable invariant curve of some
periodic orbit intersects the stable invariant curve of
another periodic orbit,
then there can be no RICs between both periodic orbits.
This \emph{obstruction criterion} was established
in~\cite{OlveraSimo1987}.

Next, we discuss these ideas in the setting of map~(\ref{eq:f}).

\begin{figure}
\iffigures
\centering
\includegraphics[height=2.4in]{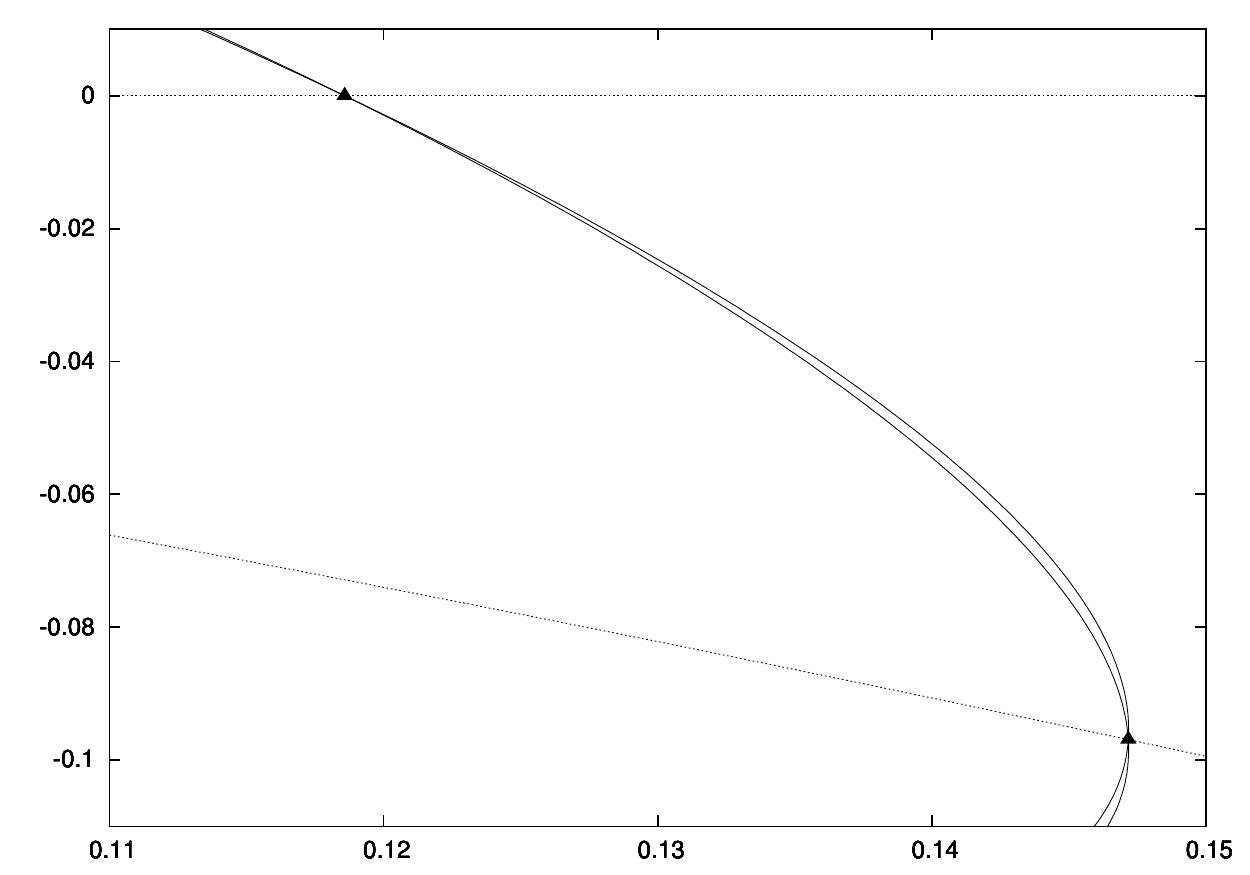}
\else
\vspace{2.4in}
\fi
\caption{A zoom of Figure~\ref{fig:StabilityRegion_mu_0p859},
but without the stability domain.
The area of the lobe $\mathcal{L}$ delimited by the separatrices
between the two primary homoclinic points marked with
solid black triangles is
$|\mathcal{L}| \approx 3.808194826948494 \times 10^{-5}$.}
\label{fig:Zoom_mu_0p859}
\end{figure}

\subsection{Singular splitting near the saddle-center bifurcation}
\label{Ssec:SeparatricesSaddleCenterBifurcation}

We saw in Section~\ref{Ssec:HamiltonianSaddleCenterBifurcation} that
the map~(\ref{eq:f}) is approximated,
after the rescaling~(\ref{eq:Scaling_SaddleCenterBifurcation}),
by the $\mu^{1/2}$-time flow of the Hamiltonian~(\ref{eq:H1_SaddleCenterBifurcation})
when $0 < \mu \ll 1$.
Besides, the Hamiltonian~(\ref{eq:H1_SaddleCenterBifurcation})
has a separatrix that encloses a region which resembles the
stability domain of the map when $0 < \mu \ll 1$.
Compare the stability domain displayed in Fig.~\ref{fig:AfterParabolic0}
with the phase portrait of the Hamiltonian~(\ref{eq:H1_SaddleCenterBifurcation})
sketched in Fig.~\ref{fig:H1_SaddleCenterBifurcation}.
The separatrix is described by the homoclinic trajectory~(\ref{eq:HomoclinicOrbit}),
which is analytic in a complex strip of width $d_0 = \pi$.

Nevertheless,
the stable and unstable invariant curves of the saddle point
$p_\rmh = (-2\phi_\rms,0)$ of our map~(\ref{eq:f}) do not coincide,
since the map is entire.
This result goes back to Ushiki~\cite{Ushiki1980}.
We have displayed the stability domain and the separatrices of
the saddle point $p_\rmh = (-2\phi_\rms,0)$ for $\mu = 0.859$ in
Fig.~\ref{fig:StabilityRegion_mu_0p859}.
We check that the separatrices enclose the stability domain.
The reversibility of our map implies that the separatrices
have a primary homoclinic point on each symmetry line~(\ref{eq:FixedSet}).
Let $\mathcal{L}$ be the region (such region is called \emph{lobe})
delimited by the pieces of the separatrices between these
two primary homoclinic points.
For instance, we display the lobe $\mathcal{L}$ for $\mu = 0.859$
 in figure~\ref{fig:Zoom_mu_0p859}.
In that case, the lobe area is
$|\mathcal{L}| \approx 3.808194826948494 \times 10^{-5}$.

Fontich and Sim\'o~\cite{FontichSimo1990} proved that
the splitting of the separatrices for any close to the identity analytic
area preserving map is exponentially small in the characteristic exponent
$h$ of the saddle point.
To be precise, they established that the splitting size is smaller
than $\Order(\rme^{-2\pi d/h})$ for any $0 < d < d_0 = \pi$.
Here, $d_0$ is the width of the analyticity strip of the
homoclinic solution of the limit Hamiltonian.
Since in our case $d_0 = \pi$, we get the upper bound
$|\mathcal{L}| \le \Order(\rme^{-c/h})$ for any $0 < c < 2\pi^2$.
We recall that $\mu = 2(\cosh h -1) = h^2 + \Order(h^4)$, so
$h \asymp \sqrt{\mu}$ as $\mu \to 0^+$.

Ten years later,
Gelfreich~\cite{Gelfreich2000} derived an asymptotic formula
for the splitting angle between the separatrices in analytic
saddle-center bifurcations, although he did not provide a
complete proof.
Gelfreich's formula, once adapted to our map,
says that $|\mathcal{L}| \asymp a_0 \rme^{-2\pi^2/h}$
as $h \to 0^+$ for some constant $a_0 \in \Rset$.

Our numerical experiments strongly suggest that there exist
some asymptotic coefficients $a_n \in \Rset$, $n \ge 0$, such that
\begin{equation}\label{eq:AsymptoticFormula}
|\mathcal{L}| \asymp \rme^{-2\pi^2/h} \sum_{n \ge 0} a_n h^{2n},
\qquad (h \to 0).
\end{equation}
This fits perfectly with both Fontich-Sim\'o's upper bound,
and Gelfreich's asymptotic formula.
Our refined asymptotic formula~(\ref{eq:AsymptoticFormula})
means that if we retain only finitely many terms of the right-hand side,
then the error will be of the order of the first discarded term.
Such refined asymptotic formulas in singular splitting problems
were first presented in~\cite{GelfreichLS1994} for the Standard map,
and first proved in~\cite{MartinSS2011} for the perturbed McMillan map.

Besides, we have numerically seen that the first asymptotic
coefficient in formula~(\ref{eq:AsymptoticFormula}) is non-zero:
\[
a_0 \approx 1.42098502709189813726617259727 \times 10^5,
\]
whereas the second asymptotic coefficient vanishes: $a_1 = 0$,
so the approximation $|\mathcal{L}| \approx a_0 \rme^{-2\pi^2/h}$
has an $\Order(h^4)$ relative error.
We have also checked that the asymptotic series
$\sum_{n \ge 0} a_n h^{2n}$ is divergent,
but its Borel transform $\sum_{n\ge 0} a_n h^{2n}/(2n)!$
has radius of convergence $2\pi^2$.
This is a typical behaviour for many other maps,
see~\cite{DelshamsRamirez1999,Ramirez2005,GelfreichSimo2008,Migueletal2013}.

Let us consider the closed curve formed by the unstable
invariant curve from the saddle point $p_\rmh$ to the primary
homoclinic point on some fixed symmetry line plus the
stable invariant curve from that primary point to $p_\rmh$.
This closed curve encloses a planar domain $\mathcal{R}$
slightly bigger than the stability domain $\mathcal{A}$,
see Fig.~\ref{fig:StabilityRegion_mu_0p859}.
The key observation is that this closed curve is an
\emph{effective barrier} when $0 < \mu \ll 1$.
The term \emph{effective} means that the flux through this closed curve
is so small that it looks like a true barrier for a very big number
of iterates of the map.
For instance, if we set $\mu=0.2$, then $h \approx 0.44357$,
\[
|\mathcal{R}| > |\mathcal{A}| \approx 2.1455 \times 10^{-3},\quad
|\mathcal{L}| \approx a_0 \rme^{-2\pi^2/h} \approx 6.7000 \times 10^{-15}.
\]
Besides, we know that the lobe area $|\mathcal{L}|$
is an exact measure of the flux through $\partial \mathcal{R}$
after one iteration of the map, see~\cite{MacKayMP1984,Meiss1992}.
This means that after $10^9$ iterates of the map $f$,
less than three thousandths parts of the points
inside $\mathcal{R}$ have escaped.
Thus, one may approximate the stability domain $\mathcal{A}$
by the region $\mathcal{R}$ in many practical situations.

Finally, we note that the numerical computation of any exponentially
small splitting quantity (angle, area, or distance)
gets complicated by problems of precision, stability, and time.
In order to overcome them,
Sim\'o proposed to use a multiple-precision arithmetic,
to expand the invariant curves up to high order, 
and to take advantage of the reversor~\cite{Simo1990}.
These ideas have been used
in~\cite{DelshamsRamirez1999,Ramirez2005,GelfreichSimo2008}.
We have also used them.

\subsection{Singular splitting near the third-order resonance}
\label{Ssec:SeparatricesThirdOrderResonance}

We saw in Section~\ref{Ssec:HamiltonianThirdOrder} that the third power
of the map~(\ref{eq:f}) is approximated,
after the rescaling~(\ref{eq:Scaling_ThirdOrder}), by the
$\epsilon$-time flow of the Hamiltonian~(\ref{eq:H1_ThirdOrderResonance})
when $\mu = 3 + \epsilon$ with $0 < |\epsilon| \ll 1$.
Besides, the Hamiltonian~(\ref{eq:H1_ThirdOrderResonance})
has three saddle points whose invariant curves coincide giving rise
to the triangle sketched in Fig.~\ref{fig:H1_ThirdOrderResonance}.

If $\mu \simeq 3$, then the stability domain of the map~(\ref{eq:f}) has a
central part with a triangular shape, that contains many RICs, and three
``sheets'', that contain points with rotation number equal to $1/3$,
attached to the vertices of that ``triangle''.
The vertices of this ``triangle'' correspond to hyperbolic
three-periodic points whose stable and unstable invariant curves
do not coincide.
There are two different splitting phenomena in this setting.
Namely,
the inner splitting (associated to the invariant curves that
enclose the ``triangle'') and
the outer splitting (associated to the invariant curves that
enclose the ``sheets'').
Each splitting should be studied separately.
The inner one is generically much smaller than the outer
one~\cite{SimoVieiro2009}.

We have displayed the stability domain for $\mu = 2.853$
in the central picture
of Fig.~\ref{fig:StabilityRegionColor_LowOrderResonances}.
The red part is the ``triangle'', the green parts are the three ``sheets'',
and the continuous  black lines are the invariant curves of the hyperbolic
three-periodic points.
We stress that,
although the value of $|\epsilon| = |\mu - 3|$ is not very small, 
the inner splitting can not be detected even after a
big magnification of our picture.
This suggest that the inner splitting is exponentially small in
$|\epsilon| = |\mu -3|$.
G.~Moutsinas~\cite{Moutsinas2016} has studied the inner splitting
in analytic area-preserving maps close to the third-order resonance.
He deduced, under a generic assumption on the third-order Birkhoff
normal form around the elliptic fixed point at the exact third-order
resonance, that the inner splitting is exponentially small
in the characteristic exponent of the third iterate of the map
at the hyperbolic three-periodic points.
To be more precise, he found that the Lazutkin homoclinic invariant
associated to some distinguished heteroclinic orbits has a
refined asymptotic formula of the form~(\ref{eq:AsymptoticFormula}),
but now $h$ is the characteristic exponent of the map $f^3$
at the three-periodic points instead of the characteristic exponent
of the map $f$ at the origin.

On the contrary, the outer splitting in the central picture
of Figure~\ref{fig:StabilityRegionColor_LowOrderResonances}
can be perceived after a suitable magnification of a small
neighborhood of a hyperbolic three-periodic point.
This visual inspection fits with the results given
in~\cite[Section~6.1]{SimoVieiro2009},
where it is established that the outer splitting associated to a
generic third-order resonance does not tend to zero as we
approach the resonance.
That is, the outer splitting is $\Order(1)$.

We can extract two practical consequences of these results.

First, set $\mu \simeq 3$ and
let $\mathcal{R}^{(1,3)}_{\rm inner}$ and $\mathcal{R}^{(1,3)}_{\rm outer}$
be the regions enclosed by suitable parts of the stable and
unstable invariant curves of the hyperbolic three-periodic points
such that $\mathcal{R}^{(1,3)}_{\rm inner}$ contains the triangular shaped
part of $\mathcal{A}$ containing many RICs and
$\mathcal{R}^{(1,3)}_{\rm outer}$ contains all the points with rotation number
equal to $1/3$.
Then the flux through the effective barrier
$\partial \mathcal{R}^{(1,3)}_{\rm inner}$ is much smaller
than the flux through $\partial \mathcal{R}^{(1,3)}_{\rm outer}$.

Second, let $\mu_\star \in [2.853,2.854]$ be the value
at which the third-order resonance is thrown away from $\mathcal{D}$.
Then $\mathcal{R}^{(1,3)}_{\rm inner}$ is a really good approximation
of the connected component $\mathcal{D}$ when $\mu \gtrsim \mu_\star$.

\subsection{Singular splitting near high-order resonances}
\label{Ssec:SeparatricesHiguerOrderResonance}

The singular splitting near resonances of order $n \ge 4$
shares several qualitative and quantitative features
with the singular splitting near the saddle-center bifurcation
and near the third-order resonance.
Let us explain this.

Let $\mu_\bullet$ and $\mu_\star$ be the values at which
the $(m,n)$-resonance emanates from $p_\rms$ and is thrown
away from $\mathcal{D}$, respectively.
Let $\mathcal{R}^{(m,n)}_{\rm inner}$
(respectively, $\mathcal{R}^{(m,n)}_{\rm outer}$)
be the region enclosed by suitable parts of the inner
(respectively, enclosed between suitable parts of the inner and outer)
branches of the stable and unstable invariant curves of
the hyperbolic $(m,n)$-periodic points.
The inner region usually looks like a red ``polygon'' with $n$
curved sides, because it is almost completely foliated by RICs.
The outer region contains the $(m,n)$-periodic chain of elliptic islands,
and it also contains part of its surrounding chaotic sea before
the $(m,n)$-resonance is thrown away.
See Figures~\ref{fig:StabilityRegionColor_FourthOrderResonance}
and~\ref{fig:StabilityRegionColor_LowOrderResonances} for several
pictures about the resonances
\[
(m,n) = \{(1,4), (1,5), (1,3), (2,5)\}.
\]

Since the flux through the borders of the inner and outer regions
can be geometrically interpreted as the area of certain
lobes~\cite{MacKayMP1984,Meiss1992},
we obtain the following information about the inner and outer flux.
The inner flux is smaller than the outer flux, and both of them
are exponentially small in $|\mu-\mu_\bullet|$~\cite{SimoVieiro2009}.
The inner region is a good approximation of the connected component
$\mathcal{D}$ when $\mu \gtrsim \mu_\star$.
See, for instance, the right picture in
Figure~\ref{fig:StabilityRegionColor_FourthOrderResonance}.
The inner and outer regions are not completely contained
in the stability domain when $\mu > \mu_\star$, since there is a small,
but not zero, flux through their borders~\cite{Ushiki1980}.

\subsection{On the obstruction criterion for the existence of RICs}
\label{Ssec:Obstructions}

\begin{figure}
\iffigures
\centering
\includegraphics[height=2.4in]{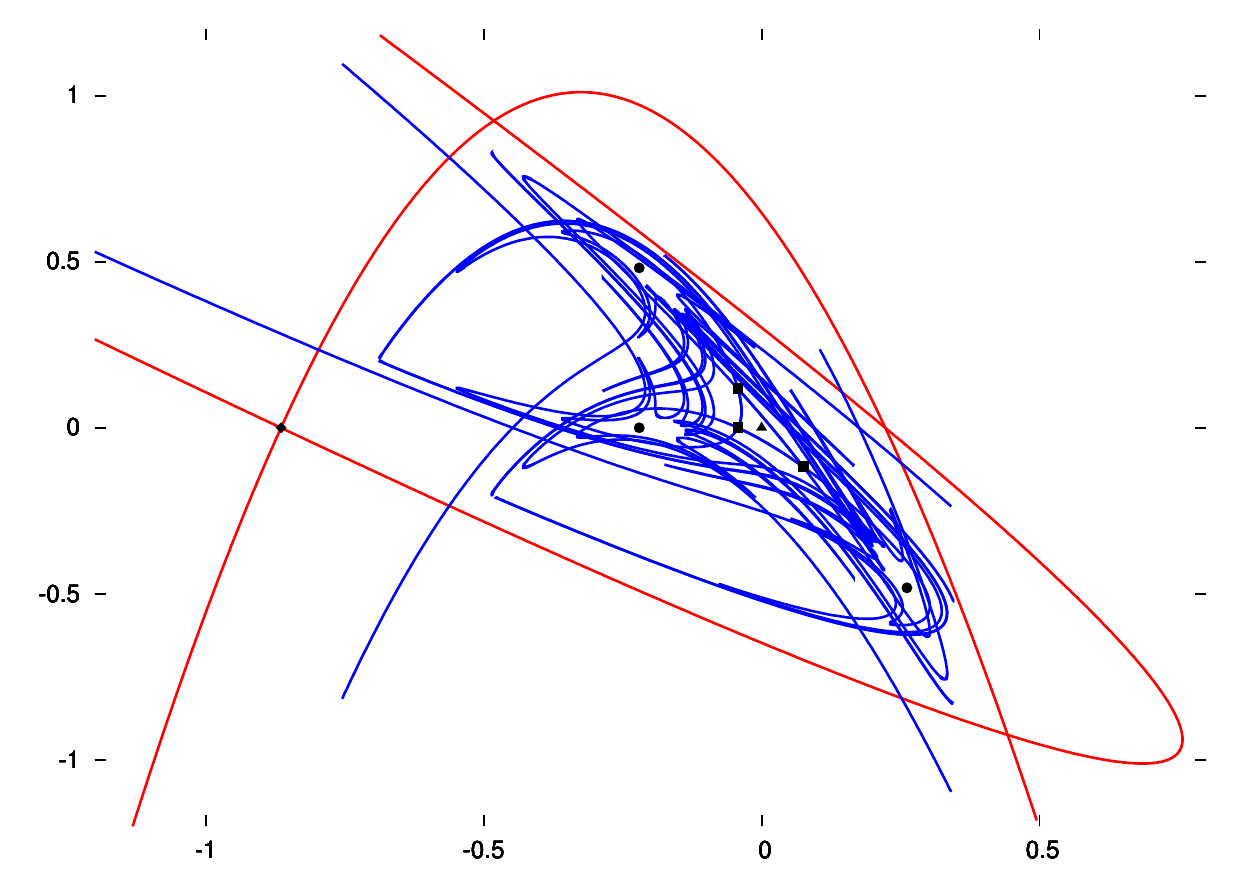}
\else
\vspace{2.4in}
\fi
\caption{
The stable and unstable invariant curves of the hyperbolic fixed point $p_{\rm h}$
and the hyperbolic $(1,3)$-periodic orbit intersect transversally for $\mu = 2.9$.
The hyperbolic/elliptic fixed point is marked with a black rhombus/triangle.
The hyperbolic/elliptic 3-periodic points are marked with black squares/circles.
The stable and unstable invariant curves of the hyperbolic fixed point
(respectively, 3-periodic points) are displayed in red (respectively, in blue).}
\label{fig:ObstructionCriterion_mu_2p9}
\end{figure}

We recall the obstruction criterion for the existence of RICs
stated in~\cite{OlveraSimo1987}.
If an area-preserving twist diffeomorphism on the annulus $\Tset \times \Rset$
has two hyperbolic periodic orbits of rotation numbers $m_1/n_1 < m_2/n_2$
such that their stable and unstable invariant curves intersect transversally,
then the map has no RIC with a rotation number $\rho \in [m_1/n_1,m_2/n_2]$.

Let $f$ be the map~(\ref{eq:f}).
The point $p_{\rm h} = ( -2\phi_{\rm s}, 0)$ is a hyperbolic fixed point
or, equivalently, a hyperbolic $(0,1)$-periodic point.
If the invariant curves of $p_{\rm h}$ intersect
transversally the invariant curves of a hyperbolic $(m,n)$-periodic orbit
of the map, then the map has no RIC with rotation number $\rho \in [0,m/n]$ and
the $(m,n)$-resonance has already escaped from the connected component $\mathcal{D}$.
Therefore, we should expect that the exact value $\mu = \mu_\star$
at which the $(m,n)$-resonance escapes from $\mathcal{D}$ coincides with the
bifurcation value at which the invariant curves of the hyperbolic
$(m,n)$-periodic orbit have their first contact with the invariant curves
of $p_{\rm h}$.

Let us present a concrete application of this idea.
We have already seen that the $(1,3)$-resonance escapes from $\mathcal{D}$
at some value $\mu = \mu_\star \in (2.853,2.854)$ by means of the
brute force method used in subsection~\ref{Ssec:StabilityDomain}.
Next, we study this escape with the obstruction criterion.

We have drawn the stable and unstable invariant curves of
$p_{\rm h}$ and the hyperbolic $(1,3)$-periodic orbit for
several values of the parameter $\mu$ in the interval $(2.8,3)$.
For instance, the $(1,3)$-resonance has already escaped
from $\mathcal{D}$ when $\mu = 2.9$, since the invariant
curves drawn in Figure~\ref{fig:ObstructionCriterion_mu_2p9}
intersect transversally.
We have obtained similar pictures for $\mu \in \{2.89, 2.88, 2.87\}$,
but the closer we are to $\mu = \mu_\star$,
the bigger part of the invariant curves we have to draw in order to
find intersections.
In particular,
pictures for $\mu = 2.88$ and $\mu = 2.87$ are far from pretty.
We have not found an intersection for $\mu = 2.86$,
because the computation and visualization of such a big part
of the invariant curves is not an easy task.

\section{Conclusions}\label{Sec:Conclusions}

We have studied the stability of longitudinal beam motion in RTMs.
Namely, we have analyzed the stability domain $\mathcal{A}$
(and its central connected component $\mathcal{D}$)
of the area-preserving map that describes the phase oscillations
using standard Dynamical Systems tools. 
We have found the range of values of the synchronous phase $\phi_\rms$
for which $\mathcal{A}$ and $\mathcal{D}$ exist.
We have studied their structure and calculated their area as a function of
$\phi_\rms$.

The knowledge of $\mathcal{A}$,
called longitudinal acceptance in the theory of particle accelerators,
is of much importance for the optimization of the beam motion in RTMs.
Indeed, the adjustment of machine parameters for the efficient acceleration
of the beam during its commissioning consists in matching the domain
in the phase space occupied by the particles emitted by an injector
(often an electron gun) to the acceptance for a given value of $\phi_\rms$.
The optimal beam matching allows to minimize beam losses and undesired excess
of strayed radiation produced by the accelerator and, so, maximize
the output beam current without increasing the current at the injection.
For the adjustment to be most efficient,
the acceptance area must be maximal and the shape of the phase domain
of the injected beam must fit the acceptance shape.
Therefore, our detailed analysis of the acceptance geometry could be useful.

Let us comment on two ``empirical'' rules used in particle
accelerators~\cite{Rand1984}.
The first rule claims that the values of $\phi_\rms$ for which an accelerator
can operate are contained in the interval of linear stability
of the synchronous trajectory.
The second rule states that the optimal values of $\phi_\rms$
are close to the middle point of such interval.
In our study, we have checked that both rules are in fact quite precise
within our RTM model, where the interval of linear stability is $(0,\phi_\rmp)$,
with
\[
\phi_\rmp := \arctan(2/\pi) \approx 32.5^{\circ}.
\]
First,
we have numerically seen that $|\mathcal{D}| > 0$ for $0 < \phi_\rms < 33^{\circ}$,
except for the value
\[
\phi_\rmu := \arctan(3/2\pi) \approx 25.5^{\circ}
\]
that corresponds to the third order resonance.
Second,
we have found that the acceptance area reaches its maximal
value $|\mathcal{A}| \approx 0.17$ at $\mu \approx 1.912$,
which roughly corresponds to
\[
\phi_\rms \approx \arctan(1.912/2\pi) \approx 16.9^{\circ} \approx \phi_\rmp/2.
\]
In fact, $|\mathcal{A}|$ is sufficiently large for a rather wide
range of the values of $\mu$.
For instance, $|\mathcal{A}|\ge 0.1$ if 
\[
\mu \in [1.027,1.071]\cup[1.079,2.037] \cup [2.245,2.827].
\]
Third, we have studied $\mathcal{A}$ and $\mathcal{D}$
in the vicinity of resonant values.
In particular, we have checked that $|\mathcal{D}| = 0$ at
the third order resonance $\phi_\rms = \phi_\rmu$,
whereas it reduces significantly, till $|\mathcal{D}| \approx 0.02$,
at the fourth order resonance
\[
\phi_\rms = \arctan(1/\pi) \approx 17.7^{\circ}.
\]
Other resonances do not lead to so sharp decreases of $|\mathcal{D}|$.
This data is quite important because one of the criteria of choosing
the design value of $\phi_\rms$, or the working point of
the machine, is to avoid values close to resonant ones.
Otherwise even a small natural drift of machine parameters may lead
to $\phi_\rms$ approaching one of the dangerous resonant values
and consequently to excessive beam losses. In this respect the asymptotic
formulas~(\ref{eq:AsymptoticSaddleCenterBifurcation})
and~(\ref{eq:AsymptoticThirdOrderResonace})
are of much interest. Another important aspect of the acceptance structure
are the elliptic islands and chaotic seas like the ones displayed in
Figs.~\ref{fig:StabilityRegionColor_FourthOrderResonance}--\ref{fig:StabilityRegionColor_LowOrderResonances}. For instance,
each drastic change in $|\mathcal{D}|$ is associated
to the escape of a chain of elliptic islands from $\mathcal{D}$.

The sizes of $\mathcal{A}$ and $\mathcal{D}$ along the
$\psi$ and $w$ axes are also important for the beam matching.
If the beam is previously bunched around $\phi_\rms$ then the bunch
length in $\psi$ must be shorter than the corresponding size of
the acceptance and the energy dispersion around $E_{n,s}$ measured
in terms of $w$ ---see~(\ref{eq:psi-w-defn})--- must be smaller
than its size in this variable.
For example, let us consider the case $\mu = 2$.
Then $\psi_{\rm max}-\psi_{\rm min} = 0.28$ for $w=0$,
and $w_{\rm max} -w_{\rm min} = 0.4$ for $\psi = 0$.
See Fig.~\ref{fig:FourthOrderResonance}.
The bunches at the
injection should fit these sizes in order to avoid beam losses during
the acceleration (in practice, to minimize beam losses).
The latter means that 
\[
\left| \frac{E_0 - E_{0,s}}{\Delta_{\rm s}} \right| < 0.064.
\]

Small accelerators do not have buncher and the beam is produced by
an electron gun which emits particles continuously,
so they occupy the whole	 interval $[0,2\pi]$ in the phase variable $\psi$
at the AS entrance.
Our results provide an estimate of the beam capture efficiency $\epsilon$;
that is, the fraction of the initial beam that is successfully accelerated.
For $\mu = 2$ this fraction is
\[
\epsilon = \frac{\psi_{\rm max}-\psi_{\rm min}}{2\pi} = 0.04.
\]
The numerical computations show that $\epsilon \le 0.13$ for all $\mu$
in our RTM model.
See Remark~\ref{rem:PhaseDeviation} in Section~\ref{Sec:MainResult}.

Here, we understand stability in a mathematical sense.
That is, we are dealing with perpetual stability,
although only $2 \cdot 10^7$ turns were considered in our numerical
computations of $\mathcal{A}$.
On the contrary, the number of turns made by each particle is
typically of just a few tens in real RTMs.
For instance, the number of turns is roughly 90 in the RTM machine of the
MAMI complex at the Institute for Nuclear Physics in Mainz, which
is nowadays the largest RTM facility in operation~\cite{Jankowiak_etal2008}.
Thus, the physical acceptance and the true capture efficiency
are larger than the mathematical stability domain and the estimates
of $\epsilon$ given above, respectively.
In fact, the difference may not be that large due to the instabilities of
trajectories after just a few iterations.

We have studied the case of the multiplicity increase factor $k = 1$,
see Section~\ref{Sec:RTMModel}.
The general case $k \in \Nset$ can be analyzed in a similar way,
and the stability domains turn out to be smaller.

Our detailed description of the acceptance is an essential widening
and improvement of results reported in~\cite{Melekhin1975}. To
the best of our knowledge our results give the first complete characterization
of the stability domain of the full non-linear model of the beam longitudinal
motion in RTMs.
We would like to emphasize that the obtained results
show the importance of the non-linear effects in the RTM beam dynamics.
This feature was well known from the experience of operation of
this type of electron accelerators.
Many properties of the RTM map are similar to those of the H\'enon map.
Also, as it was pointed out in~\cite{Melekhin1972}, similar maps appear
in the theory of anharmonic oscillator and optical theory of open
resonators.

In our study we assumed that the acceleration
gap is of zero length and that the velocity of particle is equal to
the speed of light already at the injection. The latter is not the
case for compact RTMs with the injection from a standard electron
gun. It would be useful to develop an approach in which these conditions
are relaxed. A finite-size accelerating gap can be taken into account
by introducing a transit-time factor~\cite{Wiedemann2003}.
The non-relativistic dynamics,
which is of practical importance for the RTM design,
can be considered by a corresponding modification of map~(\ref{eq:f}).
The concept of generalized synchronous particle was introduced and
phase oscillations and a corresponding map were studied in~\cite{Kubyshin_etal2008}.
It was shown that the phase of the synchronous particles changes (slips)
from turn to turn at the first orbits and this effect should be taken
into account in choosing the phase of electrons at the injection.

We have modeled the magnetic field in the RTM end magnets
by a simplified hard-edge distribution,
without taking into account neither the fringe field effect~\cite{Wiedemann2003}
nor more complicated field profiles~\cite{Vladimirov_etal2014}.
Studying the longitudinal dynamics in these cases is also of interest
for the RTM beam physics.

Finally, let us note that in our study the phase oscillations of particles
of the beam were considered as independent of the transverse oscillations,
vertical and horizontal. This approximation is valid if the amplitudes
of these oscillations are small.
A coupling between all three oscillations should be included
in a more precise and detailed analysis.

\end{document}